\documentclass[12pt]{amsart}
\usepackage{amscd}     
\usepackage{amssymb}
\usepackage{amsmath, amsthm, graphics}
\usepackage{xypic}     
\LaTeXdiagrams        
\usepackage[all]{xy}
\xyoption{2cell} \UseAllTwocells \xyoption{frame} \CompileMatrices
\allowdisplaybreaks[3]
\usepackage{amsfonts}
\usepackage{bm}
\usepackage[normalem]{ulem}
\usepackage{hyperref}
\usepackage{tikz}
\usepackage{mathtools}
\usepackage{verbatim} 
\usepackage{upgreek}

\usepackage[none]{hyphenat}

\usepackage{latexsym}
\usepackage{epsfig}
\usepackage{amsfonts}
\usepackage{enumerate}
\usepackage{times}
\usepackage{bbm}

\usepackage[T1]{fontenc}

\newtheorem{theorem}{Theorem}[subsection]

\newtheorem{lemma}[theorem]{Lemma}

\newtheorem{problem}[theorem]{Problem}

\newtheorem{proposition}[theorem]{Proposition}

\newtheorem{definition}[theorem]{Definition}

\theoremstyle{remark}

\theoremstyle{remark}

\numberwithin{equation}{section}

\newcommand{\C}{\mathbb{C}}
\newcommand{\R}{\mathbb{R}}

\def\<{\left\langle}
\def\>{\right\rangle}

\newcommand{\ebar}{\overline{\varepsilon}}

\def\b1{{\mathbf 1}}

\begin{document}

\title[Gaussian asymptotics of Jack measures from enumeration of ribbon paths]{Gaussian asymptotics of Jack measures on partitions\\ from weighted enumeration of ribbon paths}

\author[A. Moll]{Alexander Moll}
\address{Department of Mathematics\\ 
University of Massachusetts Boston\\ 
Boston, MA USA}
\email{alexander.moll@umb.edu}

\begin{abstract} In this paper we determine two asymptotic results for Jack measures $M(v^{\textnormal{out}}, v^{\textnormal{in}})$, a measure on partitions defined by two specializations $v^{\textnormal{out}}, v^{\textnormal{in}}$ of Jack polynomials proposed by Borodin-Olshanski in [European J. Combin. 26.6 (2005): 795-834].  Assuming $v^{\textnormal{out}} = v^{\textnormal{in}}$, we derive limit shapes and Gaussian fluctuations for the anisotropic profiles of these random partitions in three asymptotic regimes associated to vanishing, fixed, and diverging values of the Jack parameter.  To do so, we introduce a generalization of Motzkin paths we call ``ribbon paths'', show for arbitrary $v^{\textnormal{out}}, v^{\textnormal{in}}$ that certain Jack measure joint cumulants $\textcolor{black}{\kappa_n}$ are weighted sums of connected ribbon paths on $n$ sites with $n-1+g$ pairings, and derive our two results from the contributions of $(n,g)=(1,0)$ and $(2,0)$, respectively.  Our analysis makes use of Nazarov-Sklyanin's spectral theory for Jack polynomials.  As a consequence, we give new proofs of several results for Schur measures, Plancherel measures, and Jack-Plancherel measures.  In addition, we relate our weighted sums of ribbon paths to the weighted sums of ribbon graphs of maps on non-oriented real surfaces recently introduced by Chapuy-Do{\l}{\k{e}}ga.
\end{abstract}

  \maketitle
  
\setcounter{tocdepth}{1}

\tableofcontents

\setcounter{section}{0}

\section{Introduction and statement of results} \label{SECIntro}

\noindent In his work in multivariate statistics, Jack \cite{Jack} introduced a family of orthogonal polynomials now known as \textit{Jack polynomials}.  Since then, these polynomials have been extensively studied, rediscovered in a variety of contexts, and serve as a valuable tool in applications.  In addition, the theory of Jack polynomials is still actively developing, driven in a large part by several conjectures posed in the pioneering works by Stanley \cite{Stanley} and Goulden-Jackson \cite{GouldenJackson1996} which remain unresolved.\\
\\
\noindent In probability theory, Jack polynomials arise in the study of exactly solvable measures on partitions whose \textcolor{black}{parts} define stochastic point processes with long-range correlations.  In \cite{Ke4}, Kerov used Jack polynomials to define deformations of the Plancherel measures on partitions \cite{BaikDeiftSuidanBOOK, RomikBOOK} now known as \textit{Jack-Plancherel measures}.  These measures have been analyzed in combinatorics and probability \cite{Fulman, HoOb, MatsumotoJack, MatsumotoJack2, DoFe, DoSni, BoGoGu, GuionnetHuang} and are the simplest case of the measures on partitions of Nekrasov-Okounkov \cite{NekOk} whose work inspired tremendous research activity in the last two decades.  In particular, Nekrasov-Pestun-Shatashvili \cite{NekPesSha} and Do{\l}{\k{e}}ga-{\'S}niady \cite{DoSni} found asymptotic regimes in which the Jack-Plancherel measures form piecewise-linear limit shapes with infinitely-many local extrema, verifying one of the many influential predictions of Nekrasov-Shatashvili \cite{NekShat}.\\
\\
\noindent In this paper, we study \textit{Jack measures} on partitions proposed by Borodin-Olshanski in \textsection 3 of \cite{BoOl1} and by Olshanski in \textsection 4.4 of \cite{Ol0}.  Jack measures $M(v^{\textnormal{out}}, v^{\textnormal{in}})$ depend on two specializations $v^{\textnormal{out}}, v^{\textnormal{in}}$ of Jack polynomials and are a generalization of both Okounkov's Schur measures \cite{Ok1} and Kerov's Jack-Plancherel measures \cite{Ke4}.  We review the definition of Jack polynomials in \textsection [\ref{SUBSECJackPolynomials}] and define Jack measures in \textsection [\ref{SUBSECJackMeasures}].  Our main results are that in three different asymptotic regimes, 
\begin{enumerate}
\item The profile of a typical partition $\lambda$ sampled from a Jack measure forms a limit shape, and
\item The fluctuations of random profiles around the limit shapes are asymptotically Gaussian.
\end{enumerate}
\noindent In particular, we extend results from \cite{NekPesSha, DoSni} to show that piecewise-linear limit shapes form for arbitrary Jack measures and coincide with the \textit{dispersive action profiles} introduced by the author \cite{Moll1} and G\'{e}rard-Kappeler \cite{GerardKappeler2019}.  We state (1) and (2) in Theorem [\ref{Theorem1LLN}] and Theorem [\ref{Theorem2CLT}] below.\\
\\
\noindent To prove our results, we introduce a generalization of Motzkin paths we call \textit{ribbon paths} and determine the statistics of Jack measures from the weighted enumeration of ribbon paths.  Precisely, we show that certain Jack measure $n$th joint cumulants are weighted sums of connected ribbon paths on $n$ sites and $n-1+g$ pairings and derive our two results from the contributions of $(n,g)=(1,0)$ and $(2,0)$, respectively.  In \textsection [\ref{APPENDIXRibbonPaths}], we relate ribbon paths to the ribbon graphs on non-oriented surfaces in Chapuy-Do{\l}{\k{e}}ga \cite{ChapuyDolega2020}.  In light of this connection, our ribbon paths are of independent interest in enumerative geometry beyond their present application in probability.

\subsection{Jack polynomials} \label{SUBSECJackPolynomials} \noindent Throughout, we use the parametrization of Jack polynomials \cite{Jack} used in the study of random partitions by Nekrasov-Okounkov \cite{NekOk}: Jack polynomials $P_{\lambda}(\rho_1, \rho_2, \ldots | \ebar, \hbar)$ are polynomials in variables $\rho_1, \rho_2, \ldots$ indexed by partitions $\lambda$ and which depend on two parameters $\ebar \in \mathbb{R}$ and $\hbar>0$.  In \cite{NekOk}, for any $\ebar \in \mathbb{R}$ and $\hbar>0$, one chooses real parameters $\varepsilon_2, \varepsilon_1$ satisfying \begin{equation} \label{OmegaVariables} \varepsilon_2 < 0 < \varepsilon_1 \end{equation} and so that $\ebar \in \mathbb{R}$ and $\hbar>0$ are recovered from $\varepsilon_2, \varepsilon_1$ by \begin{eqnarray} \label{Ebar} \ebar &=& \varepsilon_1 + \varepsilon_2 \\ \label{Hbar} \hbar &=& - \varepsilon_1 \varepsilon_2. \end{eqnarray} \noindent The conventions for Jack polynomials in Stanley \cite{Stanley} and \textsection VI.10 of Macdonald \cite{Mac} are given by \begin{eqnarray} \label{Convention1} p_k &=& \rho_k / (- \varepsilon_2) \\
\noindent \label{Convention2} \alpha &=& {\varepsilon_1}/({- \varepsilon_2})\end{eqnarray}
\noindent where $p_k$ are the power sum symmetric functions and $\alpha$ is the Jack parameter in \cite{Mac,Stanley}.
\pagebreak

\noindent We now review the definition of Jack polynomials as eigenfunctions of a self-adjoint operator (\ref{OperatorIntro}).  Recall that a {partition} $\lambda$ is a weakly-decreasing sequence $0 \leq \cdots \leq \lambda_2 \leq \lambda_1$ of non-negative integers $\lambda_i$ so that $| \lambda | := \sum_{i=1}^{\infty}  \lambda_i < \infty.$  If $|\lambda|=d$, we say $\lambda$ is a partition of size $d$ with parts $\lambda_i$.  \textcolor{black}{For any partition $\lambda$, its transpose $\lambda'$ is the partition with parts $\lambda_j' = \# \{ i \ : \ \lambda_i \geq j \}$ and $l(\lambda) = \lambda_1'$ is the length of $\lambda$.} For any partition $\mu$, let $N_k'[\mu] = \# \{ i  :  \mu_i = k\}$ be the number of parts of $\mu$ of size $k$.  A partition $\mu$ is determined by these non-negative integers $(N_1'[\mu], N_2'[\mu], \ldots ) \in \mathbb{Z}_{\geq 0}^{\infty}$.  As a consequence, since the ring of polynomials $\mathbb{C}[\rho_1, \rho_2, \ldots]$ in variables $\rho_1, \rho_2, \ldots$ has a vector space basis $\rho_1^{d_1} \rho_2^{d_2} \cdots$ indexed by sequences $(d_1, d_2, \ldots) \in \mathbb{Z}_{\geq 0}^{\infty}$, one can conclude that this vector space basis is in fact indexed by partitions $\mu$ with $N_k' [ \mu] = d_k$.  In this case, we write $\rho_{\mu}:= \rho_1^{d_1} \rho_2^{d_2} \cdots$.  Note $\rho_{\mu}$ is a rescaling of the power sum $p_{\mu}$ in \cite{Mac} via (\ref{Convention1}).\\ 
\\
\noindent Next, fix $\hbar >0$ and equip $\mathbb{C}[\rho_1, \rho_2, \ldots]$ with the inner product $\langle \cdot , \cdot \rangle_{\hbar}$ depending on $\hbar>0$ defined by declaring the vector space basis $\rho_{\mu} = \rho_1^{d_1} \rho_2^{d_2} \cdots$ to be orthogonal with norm \begin{equation} \label{NormIntro}  || \rho_{\mu}||_{\hbar}^2 = \prod_{k=1}^{\infty} (\hbar k)^{d_k} d_k!. \end{equation}
\noindent In particular, $||\rho_k||^2_{\hbar} = \hbar k$.  After (\ref{Convention1}), (\ref{Convention2}), one has $||p_k||^2_{\hbar} = \alpha k$ and so $\langle \cdot , \cdot \rangle_{\hbar}$ coincides with the $\alpha$-deformed Hall inner product in \cite{Mac, Stanley}.  Write $\mathcal{F}$ for the Hilbert space completion of $\mathbb{C}[\rho_1, \rho_2, \ldots]$ with respect to $\langle \cdot, \cdot \rangle_{\hbar}$.  For any \textcolor{black}{positive integer} $k$, define unbounded operators on $\mathcal{F}$ by
\begin{eqnarray} \label{CreationOperator} \widehat{\rho}_k &:=& \rho_k \\ \label{AnnihilationOperator} \widehat{\rho}_{-k} & :=& \hbar k \frac{\partial}{\partial \rho_k}. \end{eqnarray}
\noindent These operators are mutual adjoints $\widehat{\rho}_{\pm k}^{\dagger} = \widehat{\rho}_{\mp k}$.  Set $\widehat{\rho}_0 = 0$, fix $\ebar \in \mathbb{R}$, and consider the operator
\begin{equation} \label{OperatorIntro} \widehat{{H}} ( \ebar, \hbar) =  \sum_{j_1, j_2 = 0}^{\infty} \widehat{\rho}_{j_1} \widehat{\rho}_{j_2 - j_1} \widehat{\rho}_{-j_2} +  \ebar \sum_{j=0}^{\infty} j \widehat{\rho}_{j} \widehat{\rho}_{-j} .\end{equation} 

\noindent This operator $\widehat{H}(\ebar, \hbar)$ depends explicitly on $\ebar \in \mathbb{R}$ in (\ref{OperatorIntro}) and implicitly on $\hbar>0$ through the $\widehat{\rho}_{-k}$ defined in (\ref{AnnihilationOperator}), hence one expects any eigenfunctions of $\widehat{H}(\ebar, \hbar)$ to depend on both $\ebar$ and $\hbar$ as well.  One can check that (\ref{OperatorIntro}) is self-adjoint on $(\mathcal{F}, \langle \cdot, \cdot \rangle_{\hbar})$ for the inner product \textnormal{(\ref{NormIntro})} using $\widehat{\rho}_{\pm k}^{\dagger} = \widehat{\rho}_{\mp k}$.  One can also check that (\ref{OperatorIntro}) commutes with the degree operator $\sum_{k=1}^{\infty} \widehat{\rho}_{k} \widehat{\rho}_{-k}$ which has finite-dimensional eigenspaces $\mathcal{F}_d$ spanned by $\rho_{\mu}$ with $\mu$ of size $|\mu|=d$.  As a consequence, (\ref{OperatorIntro}) has discrete spectrum with polynomial eigenfunctions.  Since the basis $\rho_{\mu}$ of $\mathcal{F}$ is indexed by partitions $\mu$, the polynomial eigenfunctions of (\ref{OperatorIntro}) must also be indexed by partitions $\lambda$.  The resulting eigenfunctions are the multivariate orthogonal polynomials of Jack \cite{Jack}.

\begin{definition} \label{DefinitionJacks} For $\hbar >0$, $\ebar \in \mathbb{R}$, and partitions $\lambda$, the Jack polynomials $P_{\lambda}(\rho_1, \rho_2, \ldots | \ebar, \hbar)$ are defined up to overall scalars as the polynomial eigenfunctions of the self-adjoint operator \textnormal{(\ref{OperatorIntro})}. \end{definition} 

\noindent The equivalence of Definition [\ref{DefinitionJacks}] and the definition of Jack polynomials in Macdonald \cite{Mac} is a result of {\textcolor{black}{Stanley \cite{Stanley}, Polychronakos \cite{Poly1995}, and Awata-Matsuo-Odake-Shiraishi \cite{AwMtOdSh}}}.  Below, we will consider the normalization of Jack polynomials with respect to the inner product $\langle \cdot, \cdot \rangle_{\hbar}$ defined by \begin{equation} \label{NormalizedJacks} P_{\lambda}^{\textnormal{norm}}(\rho_1, \rho_2, \ldots | \ebar, \hbar) := \frac{P_{\lambda}(\rho_1, \rho_2, \ldots | \ebar, \hbar)}{ || P_{\lambda} ( \rho_1, \rho_2, \ldots  | \ebar, \hbar)||_{\hbar}}. \end{equation}

\pagebreak

\subsection{Jack measures} \label{SUBSECJackMeasures} \noindent We now define Jack measures $M(v^{\textnormal{out}}, v^{\textnormal{in}})$ on partitions $\lambda$ in terms of two specializations $v^{\textnormal{out}}$ and $v^{\textnormal{in}}$ of Jack polynomials $P_{\lambda}(\rho_1, \rho_2, \ldots | \ebar, \hbar)$ satisfying certain conditions.  Recall that for any sequence $V_1, V_2, \ldots \in \mathbb{C}$, the \textit{specialization $v$ of $\rho_1, \rho_2, \ldots$ at $V_1, V_2, \ldots$} is the evaluation map $v:\C[\rho_1, \rho_2, \ldots ] \rightarrow \mathbb{C}$ which sends $Q(\rho_1, \rho_2, \ldots)$ to $Q(V_1, V_2, \ldots) \in \C$.  \textcolor{black}{Throughout this paper} we write $v = \{V_k\}_{k=1}^{\infty}$ to denote the specialization $v$ of $\rho_1, \rho_2, \ldots$ at $V_1, V_2, \ldots \in \C$.\begin{definition} \label{DefinitionJackMeasures} Fix $\ebar \in \mathbb{R}$, $\hbar >0$, and two specializations $v^{\textnormal{out}} = \{V_k^{\textnormal{out}}\}_{k=1}^{\infty}$ and $v^{\textnormal{in}} = \{V_k^{\textnormal{in}}\}_{k=1}^{\infty}$ satisfying the regularity condition
\begin{equation} \label{RegularityCondition} \sum_{k=1}^{\infty} \frac{\overline{V_k^{\textnormal{out}}} V_k^{\textnormal{in}}}{k} < \infty \end{equation}  \noindent and, for each partition $\lambda$, the non-negativity condition \begin{equation} \label{NonNegativityCondition} \overline{P^{\textnormal{norm}}_{\lambda} (V_1^{\textnormal{out}}, V_2^{\textnormal{out}}, \ldots | \ebar, \hbar) } P^{\textnormal{norm}}_{\lambda} (V_1^{\textnormal{in}}, V_2^{\textnormal{in}}, \ldots | \ebar, \hbar)  \geq 0  \end{equation} \noindent where $P_{\lambda}^{\textnormal{norm}}$ is the normalized Jack polynomial defined in \textnormal{(\ref{NormalizedJacks})}.  \noindent Then the Jack measure $M(v^{\textnormal{out}}, v^{\textnormal{in}})$ is the probability measure on the discrete set of all partitions $\lambda$ with density \begin{equation} \label{JackMeasureLaw} \textnormal{Prob}(\lambda) =  \overline{P^{\textnormal{norm}}_{\lambda} (V_1^{\textnormal{out}}, V_2^{\textnormal{out}}, \ldots | \ebar, \hbar) } P^{\textnormal{norm}}_{\lambda} (V_1^{\textnormal{in}}, V_2^{\textnormal{in}}, \ldots | \ebar, \hbar) \cdot \textnormal{exp} \Big ( -  \frac{1}{\hbar} \sum_{k=1}^{\infty} \frac{\overline{V_k^{\textnormal{out}}} V_k^{\textnormal{in}}}{k} \Big ) .\end{equation} \end{definition} 
\noindent Jack measures are probability measures on the set of all partitions $\lambda$, not on partitions of a fixed size.  For specializations $v^{\textnormal{out}}, v^{\textnormal{in}}$ satisfying (\ref{RegularityCondition}) and (\ref{NonNegativityCondition}), the formula (\ref{JackMeasureLaw}) defines a probability measure since $\sum_{\lambda} \textnormal{Prob}(\lambda) =1$ reduces to the Cauchy identity in Proposition 2.1 of Stanley \cite{Stanley}: \begin{equation} \label{CauchyIdentity}  \sum_{\lambda} \overline{P^{\textnormal{norm}}_{\lambda} (V_1^{\textnormal{out}}, V_2^{\textnormal{out}}, \ldots | \ebar, \hbar) } P^{\textnormal{norm}}_{\lambda} (V_1^{\textnormal{in}}, V_2^{\textnormal{in}}, \ldots | \ebar, \hbar)  = \textnormal{exp} \Big (  \frac{1}{\hbar} \sum_{k=1}^{\infty} \frac{\overline{V_k^{\textnormal{out}}} V_k^{\textnormal{in}}}{k} \Big ) . \end{equation} 

\noindent At $\ebar = 0$, formulas (\ref{Ebar}), (\ref{Hbar}), and (\ref{Convention2}) give $\alpha=1$, in which case the Jack polynomials (\ref{NormalizedJacks}) are Schur polynomials and the Jack measures are the Schur measures of Okounkov \cite{Ok1}.  In this paper, we focus on Jack measures $M(v,v)$ defined by a specialization $v = \{V_k\}_{k=1}^{\infty}$ satisfying \begin{equation} \label{DiagonalRegularityCondition} \sum_{k=1}^{\infty} \frac{|V_k|^2}{k} < \infty. \end{equation}

\noindent For such $v$, the Jack measure $M(v,v)$ has $v^{\textnormal{out}}  =  v^{\textnormal{in}} = v$ hence the regularity condition (\ref{RegularityCondition}) holds.  When $v^{\textnormal{out}} = v^{\textnormal{in}}$, the non-negativity conditions (\ref{NonNegativityCondition}) are always satisfied and the law (\ref{JackMeasureLaw}) is \noindent  \begin{equation} \label{DiagonalJackMeasureLaw} \textnormal{Prob}(\lambda) =  \Big | P_{\lambda}^{\textnormal{norm}}(V_1, V_2, \ldots | \ebar, \hbar)  \Big |^2 \cdot \textnormal{exp} \Big ( -  \frac{1}{\hbar} \sum_{k=1}^{\infty} \frac{|V_k|^2}{k} \Big ). \end{equation}

\noindent The assumption $v^{\textnormal{out}} = v^{\textnormal{in}}=v$ covers the case of Poissonized Jack-Plancherel measures studied in \cite{Ke4, Fulman, MatsumotoJack, DoFe, DoSni, BoGoGu, GuionnetHuang}.  For Plancherel specializations $v_{\textnormal{PL}} = \{1,0,0,\ldots\}$, Jack measures $M(v_{\textnormal{PL}}, v_{\textnormal{PL}})$ are mixtures of the Jack-Plancherel measures $M_d^{\textnormal{PL}}$ on partitions of size $d$ \textcolor{black}{with law} \textcolor{black}{{\small  \begin{equation} \label{JackPlancherelDoubleHookLaw} \textnormal{Prob}^{(d)}_{\textnormal{PL}} (\lambda) = \hbar^d d! \prod_{i=1}^{l(\lambda)} \prod_{j=1}^{\lambda_i} \frac{1}{ ( - \varepsilon_2 (\lambda_j' - i) +  \varepsilon_1 (\lambda_i - j) + \varepsilon_1) (  - \varepsilon_2 (\lambda_j' - i) +  \varepsilon_1 (\lambda_i - j) - \varepsilon_2)} \end{equation}}}\noindent by a Poisson distribution in $d$ of intensity $\tfrac{1}{\hbar}$.  \textcolor{black}{We analyze them in \textsection [\ref{SUBSECJackPlancherelLimitShapes}], \textsection [\ref{SUBSECJackPlancherelGaussianFluctuations}], \textsection [\ref{SUBSECJackPlancherelMeasureComments}], \textsection [\ref{APPENDIXdePoissonization}]. If $\ebar = 0$ then $- \varepsilon_2 = \varepsilon_1 = {\hbar}^{1/2}$ and (\ref{JackPlancherelDoubleHookLaw}) are Plancherel measures of symmetric groups; see \textsection [\ref{SUBSECPlancherelMeasureComments}].}

\subsection{Problem: asymptotics of random partition profiles as $\hbar \rightarrow 0$} \label{SUBSECProfilesProblem} In this section we pose the asymptotic problem for Jack measures addressed in this paper.  We begin with an exact computation.

\begin{proposition} \label{PropositionSizeExpectation} For $\ebar \in \mathbb{R}$, $\hbar>0$, and $v = \{V_k\}_{k=1}^{\infty}$ satisfying \textnormal{(\ref{DiagonalRegularityCondition})}, the size $| \lambda| = \sum_{i=1}^{\infty} \lambda_i$ of $\lambda$ sampled from the Jack measure $M(v,v)$ with law \textnormal{(\ref{DiagonalJackMeasureLaw})} is a random variable with expected value
\begin{equation} \label{SizeExpectationIntro} \mathbb{E} [ | \lambda| ] = \frac{1}{\hbar} \cdot \sum_{k=1}^{\infty} |V_k|^2 \end{equation} directly proportional to $\sum_{k=1}^{\infty} |V_k|^2$, inversely proportional to $\hbar$, and independent of $\ebar$. \end{proposition}

\noindent Proposition [\ref{PropositionSizeExpectation}] is verified in \textsection [\ref{SUBSECProofExpectedSizeProposition}].  As a consequence of (\ref{SizeExpectationIntro}), if $v = \{V_k\}_{k=1}^{\infty}$ satisfies $\sum_{k=1}^{\infty} \frac{|V_k|^2}{k} < \infty$ in (\ref{DiagonalRegularityCondition}) but has $\sum_{k=1}^{\infty} |V_k|^2 = \infty$, then $M(v,v)$ is well-defined but the expected value of $|\lambda|$ is infinite, i.e. the Jack measure $M(v,v)$ has heavy tails.  If $\sum_{k=1}^{\infty} |V_k|^2<\infty$, then as $\hbar \rightarrow 0$ with $v$ fixed, the typical size of the random partition $\lambda$ in (\ref{SizeExpectationIntro}) diverges at rate $1/\hbar$.  Equivalently, if we draw a partition $\lambda$ as a Young diagram with $| \lambda|$ identical $1 \times 1$ boxes, the typical area contained in the resulting random Young diagram diverges at rate $1/\hbar$.  This asymptotic behavior suggests that if we draw $\lambda$ instead as a \textit{rescaled} Young diagram with $|\lambda|$ identical boxes each of area proportional to $\hbar$, then the typical area under these $\hbar$-rescaled Young diagrams will be independent of $\hbar$ and thus one has hope to see \textcolor{black}{global} $\hbar \rightarrow 0$ asymptotics.  Kerov's insight in \cite{Ke4} was that in drawing such $\hbar$-rescaled Young diagrams, we have the freedom to draw boxes as \textit{rectangles} $R(r_2, r_1)$ of side lengths $-r_2 \sqrt{2} \times r_1 \sqrt{2}$ for any $r_2 < 0 < r_1$ so long as $\hbar = - r_1 r_2$.  He found that if we {choose} $(r_2, r_1) = (\varepsilon_2, \varepsilon_1)$ in (\ref{OmegaVariables}), (\ref{Ebar}), (\ref{Hbar}) so that $\ebar$ determines the {anisotropy} of these boxes of area $2 \hbar$, then many \textcolor{black}{quantities in the special function theory} of Jack polynomials $P_{\lambda}(\rho_1, \rho_2, \ldots | \ebar, \hbar)$ are reflected in the resulting ``anisotropic Young diagrams'' \cite{Ke4}.  We now formalize this discussion following Kerov \cite{Ke4}.

 \begin{figure}[htb]
\centering
\includegraphics[width=0.6 \textwidth]{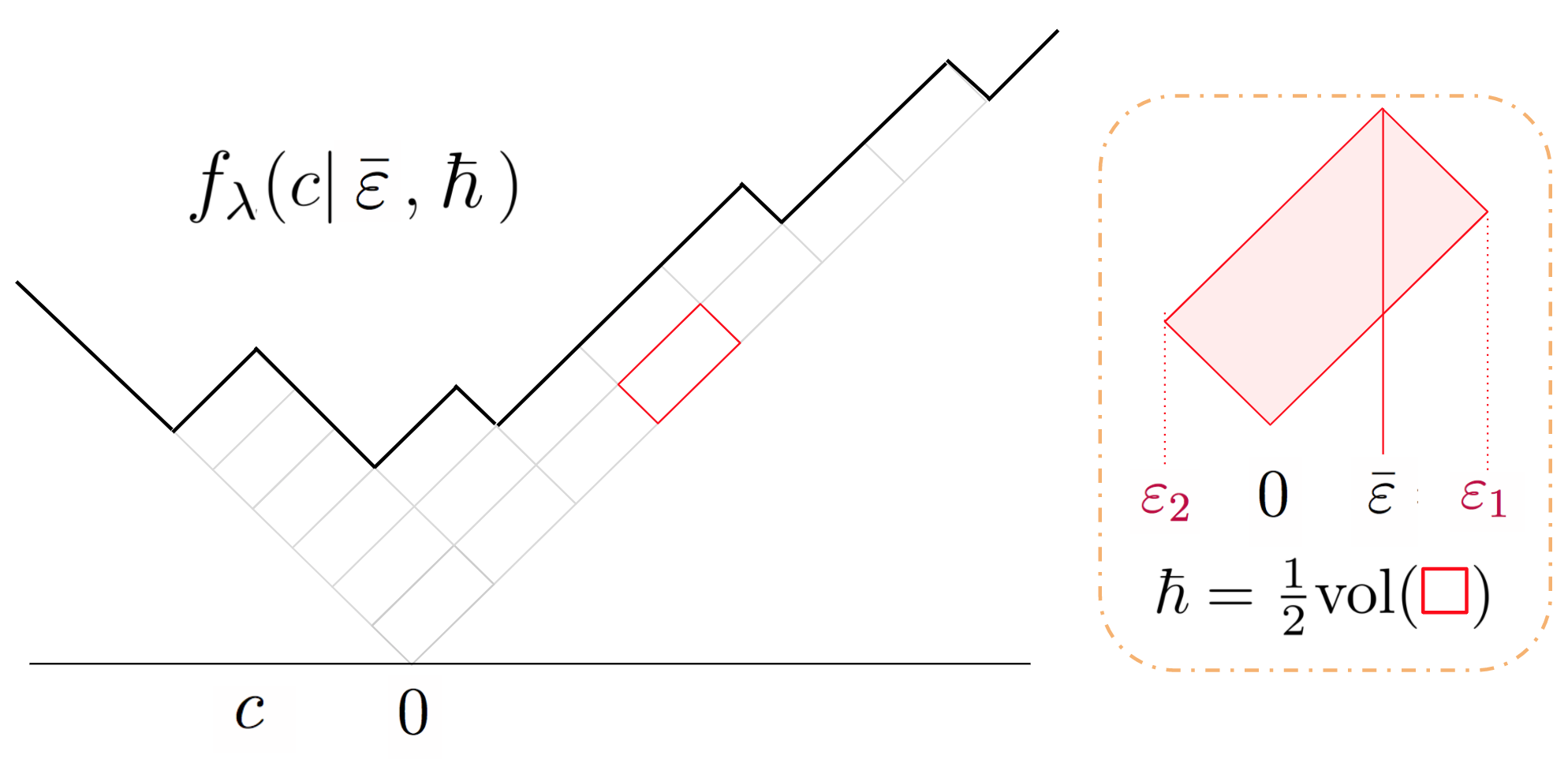}

 \vspace*{-5mm}
\caption{Graph of anisotropic partition profile $f_{\lambda}(c | \ebar, \hbar)$ of the partition $\lambda$ of size $|\lambda|=17$ with row lengths $\cdots \leq 0 \leq 1 \leq 1 \leq 1 \leq 2 \leq 5 \leq 7$ determined by the Young diagram of $\lambda$ drawn with rectangles $R(\varepsilon_2, \varepsilon_1)$ of positive anisotropy $\ebar>0$.}
\label{Anisotropic2020}
\end{figure}

\begin{definition} \label{DEFAnisotropicProfile} \noindent \textnormal{[Kerov \cite{Ke4}]} For $\ebar \in \mathbb{R}$, $\hbar>0$, and any partition $\lambda$, the anisotropic partition profile $f_{\lambda}(c | \ebar, \hbar)$ is the piecewise-linear function of $c \in \mathbb{R}$ with slopes $\pm 1$ so that the region \begin{equation} \label{TheRegion} \Omega_{\lambda}(\ebar, \hbar) =\{ (c,y) \in \R^2 \ : \ |c  | \leq y \leq f_{\lambda}(c| \ebar, \hbar) \} \end{equation} \noindent is a disjoint union of $|\lambda| = \sum_{i=1}^{\infty} \lambda_i$ rectangles $R(\varepsilon_2, \varepsilon_1)$ of side lengths $-\varepsilon_2 \sqrt{2}$, $\varepsilon_1 \sqrt{2}$ arranged in rows of positive slope indexed by $i=1,2,3, \ldots$ starting from the right so the $i$th row has $\lambda_i$ rectangles.  Here $\varepsilon_2, \varepsilon_1$ are determined from $\ebar, \hbar$ by \textnormal{(\ref{OmegaVariables}), (\ref{Ebar}), (\ref{Hbar})}.  See \textnormal{Figure [\ref{Anisotropic2020}]}.
 \end{definition}

\noindent With Definition [\ref{DEFAnisotropicProfile}], our discussion of Proposition [\ref{PropositionSizeExpectation}] above may be recapitulated as follows.  Consider the area of the region $\Omega_{\lambda}( \ebar, \hbar)$ enclosed by $f_{\lambda}(c | \ebar, \hbar)$ above $|c|$ as defined in (\ref{TheRegion}).  By Definition [\ref{DEFAnisotropicProfile}], this is the area $2\hbar$ of each box $R(\varepsilon_2, \varepsilon_1$) multiplied by the number of boxes $|\lambda|$:\begin{equation} \label{AreaDuh} \textnormal{Area} [\Omega_{\lambda}(\ebar, \hbar)]= 2 \hbar | \lambda|. \end{equation} \noindent By Proposition [\ref{PropositionSizeExpectation}], formulas (\ref{SizeExpectationIntro}) and (\ref{AreaDuh}) imply that \begin{equation} \label{AreaExpectationFixed} \mathbb{E} \Big [ \textnormal{Area} [\Omega_{\lambda}(\ebar, \hbar)] \Big ] = 2 \sum_{k=1}^{\infty} |V_k|^2 .\end{equation} \noindent In conclusion, the expected {area} of $\Omega_{\lambda}( \ebar, \hbar)$ is \textit{independent of} $\hbar$.  On the other hand, the number of available configurations of the {boundary} $\partial \Omega_{\lambda}( \ebar, \hbar)$ dramatically \textit{increases as} $\hbar \rightarrow 0$ since the condition that $\Omega_{\lambda}( \ebar, \hbar)$ be tiled by rectangles of area $2 \hbar$ becomes less and less restrictive.  The non-trivial portion of the boundary $\partial \Omega_{\lambda}( \ebar, \hbar)$ is the profile $f_{\lambda}( c | \ebar, \hbar)$, which suggests the following: 
\noindent \begin{problem} \label{ProblemStatement} {For Jack measures $M(v,v)$ on partitions $\lambda$ defined by $v = \{V_k\}_{k=1}^{\infty}$ so that \textnormal{(\ref{AreaExpectationFixed})} is finite, what is the typical statistical behavior of the random anisotropic profile $f_{\lambda}(c | \ebar, \hbar)$ as $\hbar \rightarrow 0$?} \end{problem}

\subsection{Results: limit shapes and Gaussian fluctuations as $\hbar \rightarrow 0$} \label{SUBSECResults} Our main contribution in this paper is to solve Problem [\ref{ProblemStatement}] for the random profiles $f_{\lambda}( c | \ebar, \hbar)$ assuming $c$ is independent of $\hbar$ and that $|V_k| \leq A r^k$ for some $A>0$ and $0<r<1$.  First, we prove that the random $f_{\lambda}(c | \ebar, \hbar)$ \textcolor{black}{forms} a \textit{limit shape} $\mathbf{f}(c | v; \ebar)$ as $\hbar \rightarrow 0$.  Second, we show that the fluctuations of the random profile $f_{\lambda}(c | \ebar, \hbar)$ around the limit shape $\mathbf{f}(c | v; \ebar)$ occur at the scale $\hbar^{1/2}$ and that the rescaled differences converge to a \textit{Gaussian field} $\mathbf{G}(c | v; \ebar)$ as $\hbar \rightarrow 0$: \begin{equation} \label{InformalPresentationResults} f_{\lambda} (c | \ebar, \hbar) \sim \mathbf{f}(c | v; \ebar) + \hbar^{1/2} \mathbf{G}(c | v; \ebar). \end{equation} 
 \noindent We now formulate these results precisely in Theorem [\ref{Theorem1LLN}] and Theorem [\ref{Theorem2CLT}] \textcolor{black}{below}.
 
 \begin{theorem} \label{Theorem1LLN} \textnormal{[First Main Result: Limit Shapes]} Fix $\ebar \in \mathbb{R}$, $\hbar >0$, and $v = \{V_k\}_{k=1}^{\infty}$ a specialization satisfying $|V_k| \leq A r^k$ for some $A>0$ and $0<r<1$.  Let $f_{\lambda}( c | \ebar, \hbar)$ be the random anisotropic partition profile sampled from the Jack measure $M(v,v)$.  There exists a deterministic function $\mathbf{f}(c | v; \ebar)$ of $c \in \mathbb{R}$ so that for any $p \in \mathbb{N}$, as $\hbar \rightarrow 0$ one has \textcolor{black}{convergence in probability}
 \begin{equation} \label{LLNStatement} \int_{-\infty}^{+\infty}  c^{p} \cdot  \tfrac{1}{2} f_{\lambda}''(c |\ebar, \hbar)  dc \rightarrow \int_{-\infty}^{+\infty} c^{p} \cdot \tfrac{1}{2} \mathbf{f}''(c | v; \ebar) dc\end{equation}
  
 \noindent of random variables to constants.  Moreover, we describe the limit shape $\mathbf{f}(c |  v; \ebar)$ explicitly:
  \begin{itemize}
 \item \textnormal{Regime I:} \ \ If \ $\textcolor{black}{\ebar<0}$, $\mathbf{f}(c | v; \ebar)$ is the ``dispersive action profile'' in \textnormal{Definition [\ref{DefinitionDispersiveActionProfiles}]} below which was originally introduced in \textnormal{\cite{Moll1}}.  $\mathbf{f}(c | v; \ebar)$ is piecewise-linear in $c \in \mathbb{R}$ with slopes $\pm 1$ and, for generic $v$, has infinitely-many local extrema accumulating only at $-\infty$.  In addition, the limit shape result \textnormal{(\ref{LLNStatement})} still holds if $v = v(\ebar)$ depends on $\ebar$ for fixed $\ebar<0$.\\ 
 
 \item \textnormal{Regime II:} \ If \ $\ebar = 0$, $\mathbf{f}(c | v; 0)$ is the ``convex action profile'' in \textnormal{Definition [\ref{DefinitionConvexActionProfiles}]} below originally discovered by Okounkov in \textnormal{\cite{Ok2}}.  $\mathbf{f}( c | v; 0)$ is convex in $c \in \mathbb{R}$ \textcolor{black}{and $\mathbf{f} (c | v; 0)  - |c|$ has compact support}.  In addition, if $\hbar, \ebar \rightarrow 0$ at a comparable rate $\frac{\ebar^2}{\hbar} = \gamma \in \mathbb{R}$, one still has \textnormal{(\ref{LLNStatement})} on $\mathbf{f}(c | v; 0)$.  In particular, this limit shape is independent of $\gamma$.\\ 
 
 \item \textnormal{Regime III:} If \ $\textcolor{black}{\ebar >0}$, the limit shape $\mathbf{f}(c | v; \ebar)$ is determined by the dispersive action profile $\mathbf{f}(c | -v ; - \ebar)$ in Regime I for the Jack measure $M(-v,-v)$ by the non-local relation \textnormal{(\ref{NonLocalRelationOneThree})}.
 \end{itemize}  \end{theorem}

  \noindent We prove Theorem [\ref{Theorem1LLN}] in \textsection [\ref{SECProof1}] and discuss previous limit shape results in \textsection [\ref{SECComments}].  Note that in terms of the Jack parameter $\alpha$, by (\ref{Ebar}), (\ref{Hbar}), (\ref{Convention2}), Regimes I, II, III in Theorem [\ref{Theorem1LLN}] correspond to $\textcolor{black}{\alpha \rightarrow 0}$, $\alpha$ fixed, and $\textcolor{black}{\alpha \rightarrow \infty}$, respectively.  In Regime II the rate is $\gamma = (\sqrt{\alpha} - \frac{1}{\sqrt{\alpha}})^2$.
  
  \begin{figure}[htb]
\centering
\includegraphics[width=0.95 \textwidth]{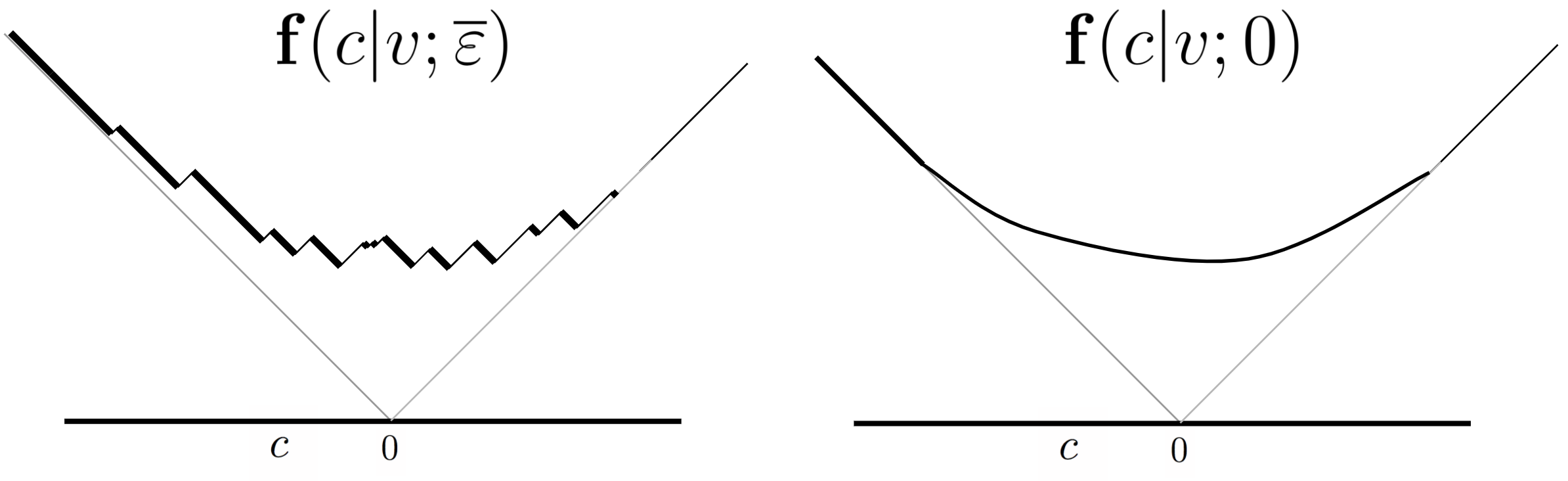}
\caption{Limit shapes of random partitions arising in Theorem [\ref{Theorem1LLN}]: the dispersive action profile $\textbf{f}( c | v; \ebar)$ with $\textcolor{black}{\ebar <0}$ and the convex action profile $\textbf{f}( c| v; 0)$.}
\label{LLNFIG2020}
\end{figure}

\noindent We now state our second main result: the random profile fluctuations are asymptotically Gaussian.
\begin{theorem} \label{Theorem2CLT} \textnormal{[Second Main Result: Gaussian Fluctuations]} Fix $\ebar \in \mathbb{R}$, $\hbar >0$, and $v = \{V_k\}_{k=1}^{\infty}$ a specialization satisfying $|V_k| \leq A r^k$ for some $A>0$ and $0<r<1$.  Let $f_{\lambda}( c | \ebar, \hbar)$ be the random anisotropic partition profile sampled from the Jack measure $M(v,v)$ \textcolor{black}{and $\mathbf{f}(c | v; \ebar)$ the limit shapes from \textnormal{Theorem [\ref{Theorem1LLN}]}}.  There exists a mean zero Gaussian random distribution $\mathbf{G}(c | v; \ebar)$ on $c \in \mathbb{R}$ so that for any $p \in \mathbb{N}$, as $\hbar \rightarrow 0$ one has weak joint convergence to Gaussian random variables:
 \begin{equation} \label{CLTStatement} \int_{-\infty}^{+\infty} c^{p}   \Bigg ( \frac{ \tfrac{1}{2} f''_{\lambda}( c | \ebar, \hbar) - \tfrac{1}{2} \mathbf{f}''(c | v; \ebar) }{ \hbar^{1/2}} \Bigg ) dc \rightarrow \int_{-\infty}^{+\infty}c^{p} \cdot \tfrac{1}{2} \mathbf{G}'' ( c | v, \ebar) dc.\end{equation}
 \noindent Moreover, we calculate the covariance of $\mathbf{G}(c | v; \ebar)$ from the limit shapes $\mathbf{f}(c | \textcolor{black}{\widetilde{v}}; \ebar)$ \textcolor{black}{for $\widetilde{v}$ near $v$} in \textnormal{(\ref{SecondCovarianceFormula})} and describe $\textcolor{black}{\mathbf{G}}(c| v; \ebar)$ in terms of the Gaussian random specialization $\boldsymbol{\varphi} = \{\boldsymbol{\varphi}_k\}_{k=1}^{\infty}$ with $\boldsymbol{\varphi}_k \in \mathbb{C}$ independent mean $0$ rotation invariant complex Gaussians with $\mathbb{E}[ | \boldsymbol{\varphi}_k|^2] =k$.  \textcolor{black}{Note that $\boldsymbol{\varphi}_k$ are independent of both $v$ and $\ebar$ and are Fourier modes of the fractional Gaussian field \textnormal{(\ref{FGFRelevantHere})}.}
 \begin{itemize}
 \item \textnormal{Regime I:} \ \  If \ $\textcolor{black}{\ebar <0}$, let $\textcolor{black}{\mathbf{f}_{\ebar}}$ be the map $\textcolor{black}{\mathbf{f}_{\ebar} : \widetilde{v}} \mapsto \mathbf{f}(c | \textcolor{black}{\widetilde{v}}; \ebar)$.  Then $\textcolor{black}{\mathbf{G}}(c | v; \ebar) = \tfrac{1}{2} \big ( d \textcolor{black}{\mathbf{f}_{\ebar}} |_v\big )_* \boldsymbol{\varphi}$ is the push-forward of $\tfrac{1}{2}\boldsymbol{\varphi}$ along the \textcolor{black}{linear map given by the} differential $d \textcolor{black}{\mathbf{f}_{\ebar}}|_{v}$ \textcolor{black}{of $\mathbf{f}_{\ebar}$ at $\widetilde{v} = v$}.\\
 
  \item \textnormal{Regime II:} \ If \ $\ebar = 0$, let $\textcolor{black}{\mathbf{f}_{0}}$ be the map $\textcolor{black}{\mathbf{f}_{0} : \widetilde{v}} \mapsto \mathbf{f}(c | \textcolor{black}{\widetilde{v}}; 0)$.  Then $\textcolor{black}{\mathbf{G}}(c | v; 0) = \tfrac{1}{2} \big ( d \textcolor{black}{\mathbf{f}_{0}} |_v\big )_* \boldsymbol{\varphi}$.  Moreover, $\mathbf{G}(c| v;0)$ coincides with the Gaussian field discovered in the asymptotics of Borodin's biorthogonal ensembles \textnormal{\cite{Bo1}} by Breuer-Duits \textnormal{\cite{BreuerDuits}}.  In addition, if $\hbar, \ebar \rightarrow 0$ and $\ebar \rightarrow 0$ at a comparable rate $\frac{\ebar^2}{\hbar} = \gamma \in \mathbb{R}$, one still has \textnormal{(\ref{CLTStatement})} except with $\mathbf{G}(c | v; 0)$ replaced by $\mathbf{G}(c | v; 0) + \sqrt{\gamma} \cdot \mathbf{X}(c|v;0)$, where $\mathbf{X}(c| v;0)$ is the deterministic function of $c \in \mathbb{R}$ in \textnormal{(\ref{FormulaForRegimeIIGaussianFieldMeanShift})}.  In particular, the covariance is independent of $\gamma$.\\
 \item \textnormal{Regime III:} If \ $\textcolor{black}{\ebar >0}$, the Gaussian process $\mathbf{G}(c | v; \ebar)$ is determined by the Gaussian process $\mathbf{G}(c | -v; - \ebar)$ in Regime I for the Jack measure $M(v,v)$ by the non-local relation \textnormal{(\ref{NonLocalRelationForCLT})}.
 \end{itemize}
\end{theorem}

\noindent We prove Theorem [\ref{Theorem2CLT}] in \textsection [\ref{SECProof2}] and discuss previous Gaussian fluctuation results in \textsection [\ref{SECComments}].\\

\pagebreak

\subsection{Polynomial expansion of transition measure joint cumulants} To prove our two results, we study the random variables $T^{\uparrow}_{\ell}(\ebar, \hbar)|_{\lambda}$ indexed by $\ell \in \mathbb{Z}_{\geq 0}$ and defined through $u \in \mathbb{C} \setminus \mathbb{R}$ by \begin{equation} \label{TransitionStatistics} \sum_{\ell=0}^{\infty} u^{- \ell-1} T^{\uparrow}_{\ell}(\ebar, \hbar) \Big |_{\lambda} = \textnormal{exp} \Bigg ( \int_{- \infty}^{+\infty} \log \Bigg [ \frac{1}{u-c} \Bigg ] \cdot \tfrac{1}{2} f_{\lambda}'' ( c | \ebar, \hbar) \Bigg ) \end{equation} 
where $f_{\lambda}( c | \ebar, \hbar)$ is the random anisotropic profile of shape $\lambda$ sampled from $M(v^{\textnormal{out}}, v^{\textnormal{in}})$.  In Kerov's theory of profiles \cite{Ke1, Ke4}, $T_{\ell}^{\uparrow}(\ebar, \hbar)|_{\lambda} = \int_{-\infty}^{+\infty} c^{\ell} d \tau^{\uparrow}_{\lambda}(c |  \ebar, \hbar)$ where $\tau^{\uparrow}_{\lambda}(c |  \ebar, \hbar)$ is the transition measure of the profile $f_{\lambda}( c | \ebar, \hbar)$.  In Theorem [\ref{Theorem3AOE}] below, for \textcolor{black}{arbitrary $v^{\textnormal{out}}$ and $v^{\textnormal{in}}$,} $n=1,2,3,\ldots$, $\vec{\ell}=(\ell_1, \ldots, \ell_n) \in \mathbb{Z}_{\geq 0}^{n}$, and $g,m=0,1,2,\ldots$, we prove there exist $W_{n,g,m}(\ell_1, \ldots, \ell_n | v^{\textnormal{out}}, v^{\textnormal{in}})$ finite so for $|| \vec{\ell}||_1 = \ell_1 + \cdots + \ell_n$, the $n$th joint cumulants of $T_{\ell}^{\uparrow}(\ebar, \hbar)|_{\lambda}$ are polynomials in $\ebar$ and $\hbar$ 
 \begin{equation} \label{IntroTheorem3AOE} \kappa_n \big ( T_{\ell_1}^{\uparrow} (\ebar, \hbar), \ldots, T_{\ell_n}^{\uparrow} (\ebar, \hbar) \big ) = \sum_{g=0}^{\tfrac{1}{2} || \vec{\ell}||_1 -n+1} \sum_{m=0}^{|| \vec{\ell}||_1}  W_{n,g,m}(\ell_1, \ldots, \ell_n | v^{\textnormal{out}}, v^{\textnormal{in}})  \hbar^{n-1+g} \textcolor{black}{\ebar^m}.  \end{equation}
\noindent \textcolor{black}{The expansion} (\ref{IntroTheorem3AOE}) implies that \textcolor{black}{in Regimes I, II, and III} as $\hbar \rightarrow 0$ \textcolor{black}{with $\ebar$ fixed or $\ebar \sim {\hbar}^{1/2}$,} \begin{equation} \label{KeyRelation} \kappa_n \big ( T^{\uparrow}_{\ell_1} ( \ebar, \hbar), \cdots, T_{\ell_n}^{\uparrow}(\ebar, \hbar) \big ) \sim \hbar^{n-1}. \end{equation} 

\noindent In \textsection [\ref{SECProof1}] and \textsection [\ref{SECProof2}], we show that (\ref{KeyRelation}) implies that \textcolor{black}{the random variables $T_{\ell}^{\uparrow}(\ebar, \hbar)|_{\lambda}$ satisfy a weak law of large numbers and central limit theorem in Regimes I, II, and III} \textcolor{black}{for Jack measures with arbitrary $v^{\textnormal{out}}, v^{\textnormal{in}}$.}  However, to characterize these limits as in our Theorems [\ref{Theorem1LLN}] and [\ref{Theorem2CLT}], we still have to calculate $W_{n,g,m}$ for $(n,g)=(1,0)$ and $(2,0)$.  \textcolor{black}{For this, we use the assumption $v^{\textnormal{out}} = v^{\textnormal{in}} = v$.}

\subsection{Weighted enumeration of ribbon paths} \textcolor{black}{Our main technical result}, Theorem [\ref{Theorem3AOE}] below, is \textcolor{black}{a refinement of} (\ref{IntroTheorem3AOE}).  \textcolor{black}{In particular, if} $v^{\textnormal{out}} = \{V^{\textnormal{out}}_k\}_{k=1}^{\infty}$, $v^{\textnormal{in}} = \{V^{\textnormal{in}}_k\}_{k=1}^{\infty}$ have $V_k^{\textnormal{out}} = V_k^{\textnormal{in}} = 0$ except for finitely-many $k=1,2,\ldots$, \textcolor{black}{Theorem [\ref{Theorem3AOE}] implies that} $W_{n,g,m}(\ell_1, \ldots, \ell_n | v^{\textnormal{out}}, v^{\textnormal{in}})$ in (\ref{IntroTheorem3AOE}) are themselves polynomials in $\overline{V_k^{\textnormal{out}}}$ and $V_k^{\textnormal{in}}$ with non-negative integer coefficients: \begin{equation} \label{LocatingRibbonPathsIntro} W_{n,g,m}(\ell_1, \ldots, \ell_n | v^{\textnormal{out}}, v^{\textnormal{in}}) = \sum_{| \mu^{\textnormal{out}}| = | \mu^{\textnormal{in}}| } C_{n,g,m}(\ell_1, \ldots, \ell_n | \mu^{\textnormal{out}}, \mu^{\textnormal{in}} ) \overline{V_{\mu^{\textnormal{out}}}^{\textnormal{out}}} V_{\mu^{\textnormal{in}}}^{\textnormal{in}} \end{equation} \noindent where $\mu^{\textnormal{out}}, \mu^{\textnormal{in}}$ are partitions of the same size and $V_{\mu} = \prod_{k=1}^{\infty} V_k^{N_k'[\mu]}$ for $N_k'[\mu] = \# \{i \ : \ \mu_i = k\}$ as in \textsection [\ref{SUBSECJackPolynomials}].  Our core innovation is to realize $C_{n,g,m}(\ell_1, \ldots, \ell_n | \mu^{\textnormal{out}}, \mu^{\textnormal{in}} ) \in \mathbb{Z}_{\geq 0}$ in (\ref{LocatingRibbonPathsIntro}) as weighted sums of new combinatorial objects we call ``ribbon paths'' with weights in $\mathbb{Z}_+$.  In Theorem [\ref{Theorem3AOE}] below, we prove $C_{n,g,m}(\ell_1, \ldots, \ell_n | \mu^{\textnormal{out}}, \mu^{\textnormal{in}} )$ are weighted sums of connected ribbon paths on $n$ sites of lengths $\ell_1, \ldots, \ell_n$ with $m$ slides, $n-g+1$ pairings, and unpaired jump profiles $\mu^{\textnormal{out}}, \mu^{\textnormal{in}}$.

\subsection{Organization of the paper} \label{SUBSECOrganization} In \textsection [\ref{SECSzegoConvex}] we review the definition of convex action profiles $\mathbf{f}( c | v ;0)$ from \cite{Moll1}.  In \textsection [\ref{SECSlidingDispersive}] we review the definition of dispersive action profiles $\mathbf{f}(c | v; \ebar)$ from \cite{Moll1}.  These convex and dispersive action profiles both appear as limit shapes in our Theorem [\ref{Theorem1LLN}].  In \textsection [\ref{SECRibbonAnisotropic}], we derive our main technical result for Jack measures in Theorem [\ref{Theorem3AOE}], the polynomial expansions (\ref{IntroTheorem3AOE}), (\ref{LocatingRibbonPathsIntro}) of joint cumulants over ribbon paths.  In \textsection [\ref{SECProof1}] we prove our limit shape result as stated in Theorem [\ref{Theorem1LLN}].  In \textsection [\ref{SECProof2}] we prove our Gaussian fluctuation result as stated in Theorem [\ref{Theorem2CLT}].  In \textsection [\ref{SECComments}] we discuss Jack measures and ribbon paths and compare our main results to previous results.  In Appendix \textsection [\ref{APPENDIXdePoissonization}], we show that our methods recover several results for the Jack-Plancherel measures \textcolor{black}{\textnormal{(\ref{JackPlancherelDoubleHookLaw})}}.  In Appendix \textsection [\ref{APPENDIXRibbonPaths}], we relate our ribbon paths to the ribbon graphs in Chapuy-Do{\l}{\k{e}}ga \cite{ChapuyDolega2020}.




\section{Szeg\H{o} paths and convex action profiles} \label{SECSzegoConvex} 

\noindent In this section we define convex action profiles $\mathbf{f}( c | v ;0)$ associated to a specialization $v = \{V_k\}_{k=1}^{\infty}$.  Convex action profiles were defined in \textsection 5 of \cite{Moll1}, originally discovered as limit shapes in \cite{Ok2}, and are limit shapes of Jack measures in Regime II of Theorem [\ref{Theorem1LLN}] as we prove in \textsection [\ref{SECProof1}].  Our presentation of convex action profiles $\mathbf{f}( c | v; 0)$ below is based on the notion of ``Szeg\H{o} paths.''  In \textsection [\ref{SUBSECSzegoPaths}] we define Szeg\H{o} paths $\gamma$ as a generalization of Catalan paths.  In \textsection [\ref{SUBSECWeightedEnumerationSzegoPaths}], for any two specializations $v^{\textnormal{out}}, v^{\textnormal{in}}$, we pose a $(v^{\textnormal{out}}, v^{\textnormal{in}})$-weighted enumeration problem for Szeg\H{o} paths of length $\ell$ in Problem [\ref{ProblemSzegoPathEnumeration}].  In \textsection [\ref{SUBSECConvexActionProfiles}], we describe the solution $W_{1,0,0}( \ell | v,v)$ of this problem in the case $v^{\textnormal{out}}= v^{\textnormal{in}}=v$ in Proposition [\ref{PropositionSzegoPaths}] and define convex action profiles $\mathbf{f}(c | v; 0)$ in terms of $W_{1,0,0}( \ell | v,v)$.  In \textsection [\ref{SUBSECCatalanVKLS}], for Plancherel specializations $v_{\textnormal{PL}}$, we show that the $(v_{\textnormal{PL}}, v_{\textnormal{PL}})$-weighted enumeration of Szeg\H{o} paths is the usual count of Catalan paths \cite{StanleyCatalanBook} and the associated convex action profile $\mathbf{f}(c | v_{\textnormal{PL}}, 0)$ is the convex profile of Vershik-Kerov \cite{KeVe} and Logan-Shepp \cite{LoSh}.

\subsection{Szeg\H{o} paths} \label{SUBSECSzegoPaths} The following notation for lattice paths is used in the remainder of the paper.  Write $\mathbb{Z}_{\pm} = \{\pm 1, \pm 2, \ldots\}$, $\mathbb{Z}_{\bullet} = \{0,1,2,\ldots\}$, and $\mathbb{Z}_{\bullet}^2$ for the quarter lattice in $\mathbb{R}^2$.  A lattice path $\mathbf{\gamma}$ in $\mathbb{Z}_{\bullet}^2$ of length $\ell$ is an ordered sequence of $\ell+1$ vertices $\textit{\textbf{y}}_0, \ldots, \textit{\textbf{y}}_{\ell}$ in $\mathbb{Z}_{\bullet}^2$.  Every lattice path determines $\ell$ steps $\textit{\textbf{e}}_1, \ldots,\textit{\textbf{e}}_{\ell} \in \mathbb{Z}^2$ by the equation $\textit{\textbf{y}}_{i+1} = \textit{\textbf{y}}_i + \textit{\textbf{e}}_i$.  Write $\mathbf{E}[\gamma]$ for the set of steps $\textit{\textbf{e}}_i$ in $\gamma$.  We now define a special class of lattice paths in $\mathbb{Z}_{\bullet}^2$ which we call \textit{Szeg\H{o} paths}.
\begin{definition} \label{DefinitionSzegoPaths} A Szeg\H{o} path $\mathbf{\gamma}$ of length $\ell$ is a lattice path $\textbf{y}_0, \textbf{y}_1, \ldots, \textbf{y}_{\ell}$ in $\mathbb{Z}_{\bullet}^2$ with
\begin{itemize}
\item Boundary conditions: $\textbf{y}_0 = (0, 0)$ and $\textbf{y}_{\ell} = (0,\ell)$ 
\item Step types: $\textbf{y}_{i+1} = \textbf{y}_{i} + \textbf{\textit{e}}_i$ with steps $\textit{\textbf{e}}_i = (1,k)$ for some non-zero integer $k \in \mathbb{Z}_- \cup \mathbb{Z}_+$.  
\end{itemize}
\noindent We say that the step $\textit{\textbf{e}}_i = (1,k)$ is a jump of degree $k$ and write $\deg (\textit{\textbf{e}}_i)=k$. \end{definition}

\begin{figure}[htb]
\centering
\includegraphics[width=0.6 \textwidth]{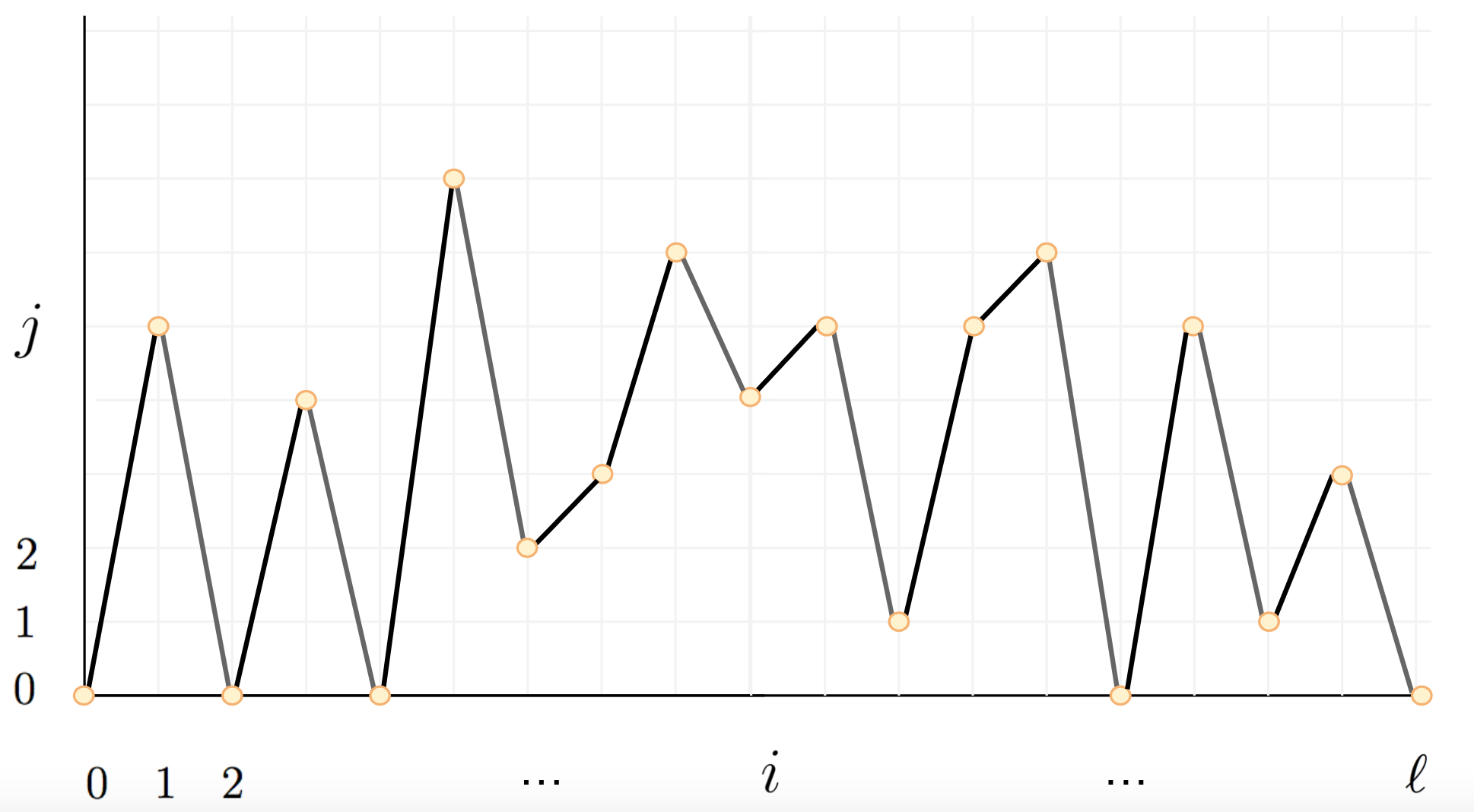}
\caption{A Szeg\H{o} path of length $\ell=18$.}
\label{SzegoPath2020}
\end{figure}

\noindent In a Szeg\H{o} path, all vertices lie in the first quadrant and are of the form $\textit{\textbf{y}}_i = (i, j_i)$ for some $j_i \in \mathbb{Z}_{\bullet}$.  We depict a Szeg\H{o} path in Figure [\ref{SzegoPath2020}].  By definition, any jump $\textit{\textbf{e}}$ in $\gamma$ has non-zero degree $\textnormal{deg}(\textit{\textbf{e}}) \neq 0$.  As a consequence, horizontal steps $(i,j) \rightarrow (i+1,j)$ with $\textit{\textbf{e}}_i = (1,0)$ are forbidden.  To keep track of positive and negative degree jumps, write $\mathbf{E}^{\pm}_J [ \gamma ] = \{ \textit{\textbf{e}} \ \textnormal{ in }  \gamma  :  \deg (\textbf{\textit{e}}) \in \mathbb{Z}_{\pm} \}$ and \begin{equation} \label{SzegoJumpsDegreePlusMinus} \mathbf{E}[\gamma] = \mathbf{E}_J^-[\gamma] \cup \mathbf{E}_J^+[\gamma]. \end{equation}

\noindent If jump degrees are restricted to $k = \pm 1$, Szeg\H{o} paths are Catalan paths \cite{StanleyCatalanBook}; see \textsection [\ref{SUBSECCatalanVKLS}].

\subsection{Weighted enumeration of Szeg\H{o} paths} \label{SUBSECWeightedEnumerationSzegoPaths} We now define a $\C$-weight $\mathfrak{W}$ for Szeg\H{o} paths. 

\begin{definition} \label{DefinitionSzegoPathsWeight} Let $v^{\textnormal{out}}=\{V_k^{\textnormal{out}}\}_{k=1}^{\infty}$, $v^{\textnormal{in}}=\{V_k^{\textnormal{in}}\}_{k=1}^{\infty}$ be two specializations.  Write $\overline{V}$ for the complex conjugate of $V \in \mathbb{C}$.  The $(v^{\textnormal{out}}, v^{\textnormal{in}})$-weight $\mathfrak{W}( \cdot | v^{\textnormal{out}}, v^{\textnormal{in}})$ on Szeg\H{o} paths $\gamma$ is defined by \begin{equation} \label{SzegoWeightFormula} \mathfrak{W}( {\gamma} |v^{\textnormal{out}}, v^{\textnormal{in}})= \prod_{\textbf{\textit{e}} \in \mathbf{E}_J^-[\gamma]} \overline{V_{-\deg ( \textbf{\textit{e}})}^{\textnormal{out}}} \prod_{\textbf{\textit{e}} \in \mathbf{E}_J^+[\gamma]} {V_{\deg ( \textbf{\textit{e}})}^{\textnormal{in}}}
 .\end{equation} \end{definition}

\noindent For example, $\mathfrak{W}(\gamma | v^{\textnormal{out}}, v^{\textnormal{in}}) = V_5^{\textnormal{in}} \overline{V_{5}^{\textnormal{out}}} V_{4}^{\textnormal{in}}  \overline{V_{4}^{\textnormal{out}}}  V_7^{\textnormal{in}}  \overline{V_{5}^{\textnormal{out}} } V_1^{\textnormal{in}}  V_3^{\textnormal{in}}  \overline{V_{2}^{\textnormal{out}}}  V_1^{\textnormal{in}}  \overline{V_{4}^{\textnormal{out}}}  V_4^{\textnormal{in}}  V_1^{\textnormal{in}}  \overline{V_{6}^{\textnormal{out}}}  V_5^{\textnormal{in}}  \overline{V_{4}^{\textnormal{out}}} V_2^{\textnormal{in}}  \overline{V_{3}^{\textnormal{out}}} $ for $\gamma$ in Figure [\ref{SzegoPath2020}].  (\ref{SzegoWeightFormula}) determines a $(v^{\textnormal{out}}, v^{\textnormal{in}})$-weighted enumeration problem for Szeg\H{o} paths.

\begin{problem} \label{ProblemSzegoPathEnumeration} Fix two specializations $v^{\textnormal{out}}=\{V_k^{\textnormal{out}}\}_{k=1}^{\infty}$, $v^{\textnormal{in}}=\{V_k^{\textnormal{in}}\}_{k=1}^{\infty}$.  Let $\Gamma_{1,0,0}(\ell)$ be the infinite set of Szeg\H{o} paths $\gamma$ of length $\ell$ and consider the weight $\mathfrak{W}(\gamma |v^{\textnormal{out}}, v^{\textnormal{in}})$ in \textnormal{(\ref{SzegoWeightFormula})}.  Determine \begin{equation} \label{SzegoPathWFunction} W_{1,0,0}(\ell | v^{\textnormal{out}}, v^{\textnormal{in}}):= \sum_{\gamma \in \Gamma_{1,0,0}(\ell)} \mathfrak{W}(\gamma |v^{\textnormal{out}}, v^{\textnormal{in}}). \end{equation} 
\end{problem}

\subsection{Convex action profiles from Szeg\H{o} paths} \label{SUBSECConvexActionProfiles} 
\noindent In the case $v^{\textnormal{out}} = v^{\textnormal{in}} = v$, Problem [\ref{ProblemSzegoPathEnumeration}] for Szeg\H{o} paths has a solution that leads to the definition of convex action profiles $\mathbf{f}( c | v; 0)$ \textcolor{black}{from \cite{Moll1}}.
\begin{proposition} \label{PropositionSzegoPaths} Choose a specialization $v=\{V_k\}_{k=1}^{\infty}$ so that \textcolor{black}{the symbol} \begin{equation} \label{Symbol} v(x) = \sum_{k=1}^{\infty} \big ( V_k e^{- \textnormal{\textbf{i}}k x} + \overline{V_k} e^{+ \textnormal{\textbf{i}}k x} \big )   \end{equation} \noindent satisfies $||v||_{\infty} = \sup_{x} |v(x)| <\infty$.  Then \textnormal{(\ref{SzegoPathWFunction})} for Szeg\H{o} paths is determined for $u \in \mathbb{C} \setminus \mathbb{R}$ by \begin{equation} \label{SzegoPathTUpFORMULA} \lim_{\ell' \rightarrow \infty} \sum_{\ell=0}^{\ell'} u^{- \ell-1} W_{1,0,0}(\ell | v,v) =  \textnormal{exp} \Bigg ( \int_{0}^{2\pi} \log \Bigg [ \frac{1}{u - v(x) } \Bigg ] \frac{dx}{2 \pi} \Bigg ). \end{equation} \end{proposition}
\begin{itemize}
\item \noindent \textit{Proof of \textnormal{Proposition [\ref{PropositionSzegoPaths}]}.} We have to show (\ref{SzegoPathTUpFORMULA}).  This equality has been verified in \textsection 5.4 in \cite{Moll1}.  To see this, consider the infinite Toeplitz matrix \begin{equation}\label{ToeplitzMatrix} L_{\bullet}(v) = \begin{bmatrix} 
0 & \overline{V_{1}} & \overline{V_{2}} & \overline{V_{3}} & \cdots  \\ 
V_{1} &0 & \overline{V_{1}} & \overline{V_{2}} & \ddots  \\
 V_{2} & V_{1} & 0& \overline{V_{1}} & \ddots \\
  V_{3} & V_2 & V_1 & 0 & \ddots \\
  \vdots & \vdots & \ddots & \ddots & \ddots \end{bmatrix} \end{equation}
  
  \noindent which is the $V_0=0$ case of (5.5) in \cite{Moll1}.  For $j \in \mathbb{Z}_{\bullet}$, let $\psi_j = [ 0 \ 0 \ \cdots 0 \ 1 \ 0 \cdots ]^T$ denote the column vector which is $1$ in the $j$th entry and $0$ otherwise.  \textcolor{black}{Note that $\mathbb{Z}_{\bullet} = \{0,1,2,\ldots\}$ and so the first entry of $\psi_j$ is in fact the $0^{\textnormal{th}}$ entry}. Consider the inner product $\langle \cdot, \cdot \rangle$ on the span of $\psi_j$ for which $\psi_j$ are orthonormal.  Directly from the definitions,  \begin{equation} \label{SzegoThisThing} W_{1,0,0}(\ell | v,v)  = \langle \psi_0, L_{\bullet}(v)^{\ell} \psi_0 \rangle.\end{equation} Using (\ref{SzegoThisThing}) and essential self-adjointness of (\ref{ToeplitzMatrix}) from \cite{Moll1}, the limit in (\ref{SzegoPathTUpFORMULA}) exists and is \begin{equation} \label{SzegoTheLove}{\mathbf{T}}^{\uparrow}(u | v;0) := \langle \psi_0, \frac{1}{u - L_{\bullet}(v)} \psi_0 \rangle \end{equation} \textcolor{black}{if $|| v||_{\infty} <\infty$}. Finally, \textnormal{Theorem 5.4.4} of \textnormal{\cite{Moll1}} says (\ref{SzegoTheLove}) is the right-hand side of (\ref{SzegoPathTUpFORMULA}). $\square$
  \end{itemize} 
  \noindent We chose the term ``Szeg\H{o} paths'' in this paper since in \cite{Moll1}, the equality of (\ref{SzegoTheLove}) and (\ref{SzegoPathTUpFORMULA}) is shown to be equivalent to Szeg\H{o}'s First Theorem in light of Remark 2 after Theorem 1.6.1 in Simon \cite{SimonSzego}.  We can now define ``convex action profiles'' \textcolor{black}{whose existence and uniqueness is due to \cite{Moll1}.}
\begin{definition} \label{DefinitionConvexActionProfiles} Let $v=\{V_k\}_{k=1}^{\infty}$ be a specialization so $||v||_{\infty} = \sup_x |v(x)| <\infty$ holds for $v(x)$ in \textnormal{(\ref{Symbol})}.  The convex action profile $\mathbf{f}(c | v;0)$ is the \textcolor{black}{unique} function of $c \in \mathbb{R}$ characterized by the conditions that $\mathbf{f}( c | v; 0) \sim |c|$ as $c \rightarrow \pm \infty$ and, for all $u \in \mathbb{C} \setminus \mathbb{R}$, $\mathbf{T}^{\uparrow}(u | v; 0)$ in \textnormal{(\ref{SzegoTheLove})} is both \begin{equation} \label{DefinitionConvexActionProfilesFORMULA} \lim_{\ell' \rightarrow \infty} \sum_{\ell=0}^{\ell'} u^{- \ell-1} W_{1,0,0}(\ell | v,v) = \textnormal{exp} \Bigg ( \int_{-\infty}^{+\infty} \log \Bigg [ \frac{1}{u-c} \Bigg ] \cdot \tfrac{1}{2} \mathbf{f}''(c | v; 0) dc \Bigg ) \end{equation} \noindent where $W_{1,0,0}(\ell | v,v)$ are the $(v,v)$-weighted sums over Szeg\H{o} paths in \textnormal{(\ref{SzegoPathWFunction})}. \end{definition}
\noindent Convex action profiles were introduced in Definition 1.4.1 in \cite{Moll1} and appeared originally as limit shapes in \textsection 4.2.1 of \cite{Ok2}.  In \cite{Moll1}, many properties of convex action profiles $\mathbf{f}(c | v;0)$ are established.  By Proposition [\ref{PropositionSzegoPaths}], convex action profiles are equivalently characterized for $p=0,1,2,\ldots$ by \begin{equation} \label{PushForwardFormulas} \int_{-\infty}^{+\infty} c^p \cdot \tfrac{1}{2} \mathbf{f}''(c | v; 0) dc = \int_0^{2\pi} v(x)^p \  \tfrac{dx}{2\pi}, \end{equation} \noindent namely $\tfrac{1}{2} \mathbf{f}''(c | v; 0) dc$ is a probability measure that can be described as the push-forward of the uniform measure $\tfrac{dx}{2\pi}$ on $[0,2\pi]$ along $v: [0,2\pi] \rightarrow \mathbb{R}$.  Since this measure is non-negative, $\mathbf{f}(c | v; 0)$ is convex in the variable $c \in \mathbb{R}$.  Moreover, the support of $\tfrac{1}{2} \mathbf{f}''(c|v;0) dc$ is contained in $[\inf_x v, \sup_x v]$ and is connected if $v(x)$ is continuous.  As a consequence, $\mathbf{f}(c | v;0) = |c|$ if $c \not\in [\inf_x v, \sup_x v]$.  See Figure [\ref{LLNFIG2020}] for a depiction of a convex action profile.  For further detail, see Figure 8 in \cite{Moll1}.

\subsection{Example: Catalan paths and Vershik-Kerov Logan-Shepp profile} \label{SUBSECCatalanVKLS} Consider the case
\begin{equation} \label{PlancherelSpecialization} V_k^{\textnormal{PL}} = \begin{cases} 1 \ \ \ \ \ \textnormal{if \ } k=1 \\ 0 \ \ \ \ \ \textnormal{otherwise} \end{cases} \end{equation}
\noindent of a Plancherel specialization $v_{\textnormal{PL}} = \{1,0,0,\ldots\}$.  In this case, the weight (\ref{SzegoWeightFormula}) simplifies \begin{equation} \mathfrak{W}(\gamma | v_{\textnormal{PL}}, v_{\textnormal{PL}}) 
= \begin{cases} 1 \ \ \ \ \textnormal{if $\gamma$ is a Catalan path} \\ 0 \ \ \ \ \textnormal{otherwise}, \end{cases} 
\end{equation} \noindent $W_{1,0,0}(\ell |v_{\textnormal{PL}}, v_{\textnormal{PL}})$ in (\ref{SzegoPathWFunction}) are the Catalan numbers, and $\mathbf{T}^{\uparrow}(u | v_{\textnormal{PL}};0)$ in (\ref{SzegoTheLove}) and (\ref{SzegoPathTUpFORMULA}) solves \begin{equation} \label{CatalanQuadraticTing} \mathbf{T}^{\uparrow}(u | v_{\textnormal{PL}};0) + \frac{1}{\mathbf{T}^{\uparrow}(u | v_{\textnormal{PL}};0)} =u.\end{equation} \noindent As is well-known and reviewed in Example 7.2.9 in Kerov \cite{Ke1}, for the specialization $v_{\textnormal{PL}}$, the solution of the algebraic equation (\ref{CatalanQuadraticTing}) is the Stieltjes transform of Wigner's semi-circle law and the associated convex action profile is the profile Vershik-Kerov \cite{KeVe} and Logan-Shepp \cite{LoSh} \begin{equation} \label{FormulaVKLS} \mathbf{f}(c | v_{\textnormal{PL}};0) = \begin{cases} \tfrac{2}{\pi} (c \arcsin ( \tfrac{c}{2}) + \sqrt{4 - c^2}),  \ \ \ \ \ | c| \leq 2 \\ |c| ,    \ \ \ \ \ \ \ \  \ \ \ \ \  \ \ \ \ \  \ \ \ \ \ \ \ \ \ \  \ \ \ \ \  \ \ \ \ \  \ \ \ \ \ |c| \geq 2. \end{cases} \end{equation}
\noindent As is consistent with (\ref{PushForwardFormulas}), $v_{\textnormal{PL}}(x) = 2\cos x$ so $\mathbf{f}(c | v_{\textnormal{PL}}; 0)$ is characterized for $p=0,1,2,\ldots$ by \begin{equation} \int_{-\infty}^{+\infty} c^{p} \cdot \tfrac{1}{2} \mathbf{f}''(c|v;0) dc = \int_{0}^{2 \pi} \Big ( 2 \cos x \Big ) ^p \tfrac{dx}{2 \pi}. \end{equation}




\section{Sliding paths and dispersive action profiles} \label{SECSlidingDispersive}

 \noindent In this section we define dispersive action profiles $\mathbf{f}( c | v ;\ebar)$ for any specialization $v = \{V_k\}_{k=1}^{\infty}$.  Dispersive action profiles were defined in \textsection 3 of \cite{Moll1}, implicitly studied in \cite{GerardKappeler2019}, and are limit shapes of Jack measures in Regimes I and III of Theorem [\ref{Theorem1LLN}] as we prove in \textsection [\ref{SECProof1}].  Our presentation of dispersive action profiles below is based on the notion of ``sliding paths.''  In \textsection [\ref{SUBSECSlidingPaths}] we define sliding paths as a generalization of Szeg\H{o} paths with ``slides''.  In \textsection [\ref{SUBSECWeightedEnumerationSlidingPaths}], for any two specializations $v^{\textnormal{out}}, v^{\textnormal{in}}$, we pose a $(v^{\textnormal{out}}, v^{\textnormal{in}})$-weighted enumeration problem for sliding paths of length $\ell$ with $m$ slides in Problem [\ref{ProblemSlidingPathEnumeration}].  In \textsection [\ref{SUBSECDispersiveActionProfiles}], we describe the solution $W_{1,0,m}( \ell | v,v)$ of this problem in the case $v^{\textnormal{out}}= v^{\textnormal{in}}=v$ in Proposition [\ref{PropositionSlidingPaths}] and define dispersive action profiles $\mathbf{f}(c | v; \ebar)$ in terms of $W_{1,0,m}( \ell | v,v)$.  In \textsection [\ref{SUBSECMotzkinNPS}], for Plancherel specializations $v_{\textnormal{PL}}$, we show that the $(v_{\textnormal{PL}}, v_{\textnormal{PL}})$-weighted enumeration of sliding paths is a weighted enumeration of Motzkin paths \cite{StanleyVol2} with associated dispersive action profile $\mathbf{f}(c | v_{\textnormal{PL}}, \ebar)$ from Nekrasov-Pestun-Shatashvili \cite{NekPesSha} and Do{\l}{\k{e}}ga-{\'S}niady \cite{DoSni}.

\subsection{Sliding paths} \label{SUBSECSlidingPaths} 
\begin{figure}[htb]
\centering
\includegraphics[width=0.55 \textwidth]{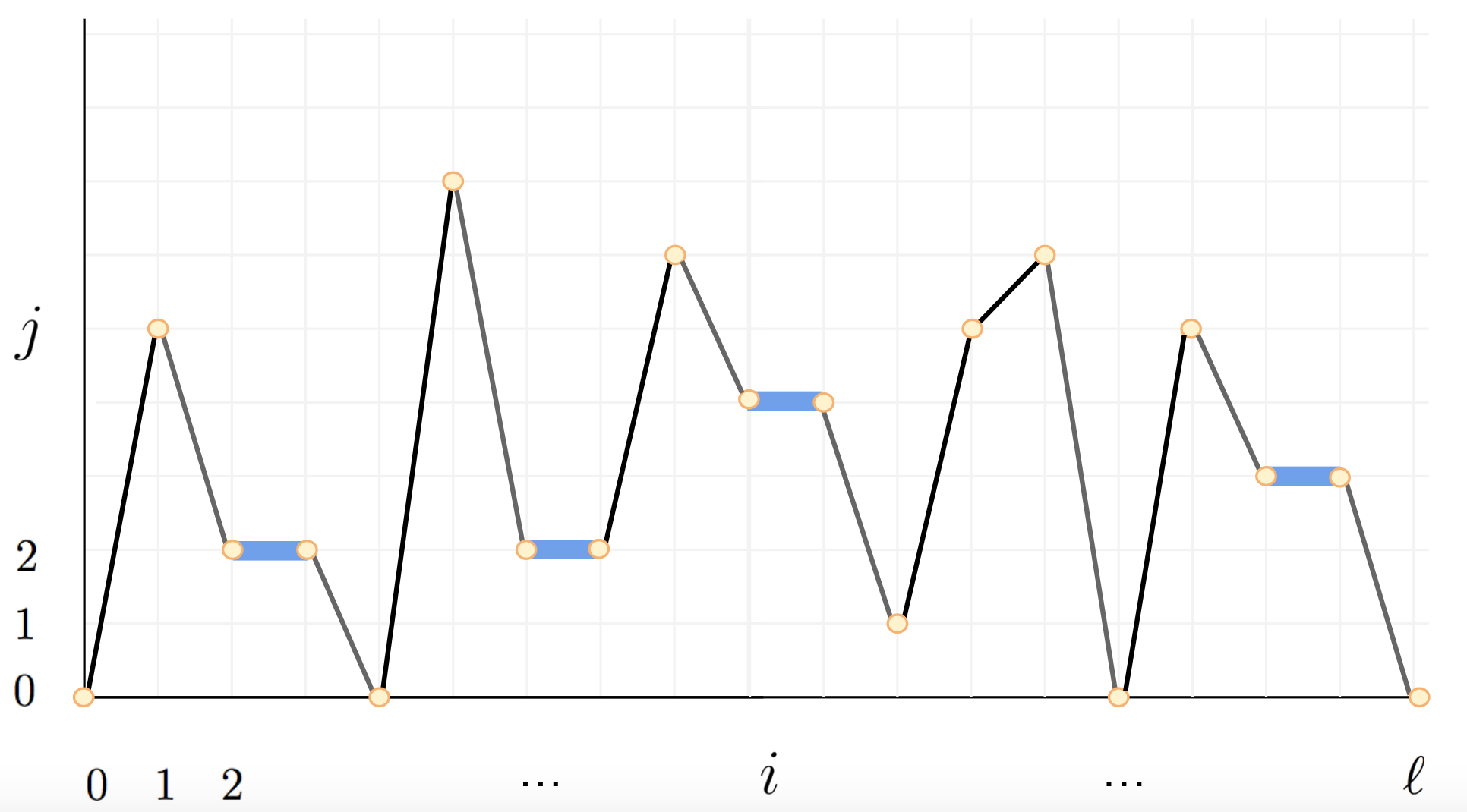}
\caption{A sliding path of length $\ell=18$ with $m=4$ slides.}
\label{SlidingPath2020}
\end{figure}

We now define a generalization of Szeg\H{o} paths in \textsection [\ref{SUBSECSzegoPaths}] we call \textit{sliding paths}.
\begin{definition} \label{DefinitionSlidingPaths} A sliding path $\mathbf{\gamma}$ of length $\ell$ is a lattice path $\textbf{y}_0, \textbf{y}_1, \ldots, \textbf{y}_{\ell}$ in $\mathbb{Z}_{\bullet}^2$ with
\begin{itemize}
\item Boundary conditions: $\textbf{y}_0 = (0, 0)$ and $\textbf{y}_{\ell} = (0,\ell)$ 
\item Step types: $\textbf{y}_{i+1} = \textbf{y}_{i} + \textbf{\textit{e}}_i$ with steps $\textit{\textbf{e}}_i = (1,k)$ for some integer $k \in \mathbb{Z}$ possibly zero.
\end{itemize}
\noindent If $k=0$ and $\textbf{y}_i = (i,j)$, we say that $\textit{\textbf{e}}_i = (1,0)$ is a slide at height $j$ and write $\textnormal{{height}}(\textit{\textbf{e}}_i) = j$.\\
\noindent If $k \neq 0$, we say that $\textit{\textbf{e}}_i = (1,k)$ is a jump of degree $k$ and write $\deg(\textit{\textbf{e}}_i) = k$.
\end{definition}
\noindent Since each step in a sliding path is by $1$ in the first component, the vertices in a sliding path must be of the form $\textit{\textbf{y}}_i = (i, j_i)$ for some $j_i \in \mathbb{Z}_{\bullet}$.  Unlike Szeg\H{o} paths, horizontal steps of the form $(i,j) \rightarrow (i+1,j)$ with $\textit{\textbf{e}}_i = (1,0)$ are allowed: we refer to such steps as ``slides''.  Otherwise, if $\textit{\textbf{e}}_i = (1,k)$ for $k \neq 0$, we call the step $\textit{\textbf{e}}_i$ a ``jump''.  In this way, the set of steps $\textit{\textbf{e}}_i$ in $\gamma$ is partitioned \begin{equation} \mathbf{E}[\gamma] =  \mathbf{E}_S[\gamma]\cup \mathbf{E}_J^-[\gamma] \cup \mathbf{E}_J^+[\gamma]   \end{equation}

\noindent into the set of slides $\mathbf{E}_S[\gamma]$ and jumps $\mathbf{E}_J^{\pm}[\gamma]$ as in (\ref{SzegoJumpsDegreePlusMinus}).  We depict a sliding path in \textcolor{black}{Figure [\ref{SlidingPath2020}].}

\noindent If the jump sizes are restricted to $k = \pm 1$, sliding paths are Motzkin paths \cite{StanleyVol2}.  We discuss the case of Motzkin paths in \textsection [\ref{SUBSECMotzkinNPS}].  First, we show that a weighted enumeration problem for sliding paths is solvable and is equivalent to a result of the author \cite{Moll1} and G\'{e}rard-Kappeler \cite{GerardKappeler2019}. 

\subsection{Weighted enumeration of sliding paths} \label{SUBSECWeightedEnumerationSlidingPaths} We now define a $\C$-weight $\mathfrak{W}$ for sliding paths.

\begin{definition} \label{DefinitionSlidingPathsWeight} Let $v^{\textnormal{out}}=\{V_k^{\textnormal{out}}\}_{k=1}^{\infty}$, $v^{\textnormal{in}}=\{V_k^{\textnormal{in}}\}_{k=1}^{\infty}$ be two specializations.  Write $\overline{V}$ for the complex conjugate of $V \in \mathbb{C}$.  The $(v^{\textnormal{out}}, v^{\textnormal{in}})$-weight $\mathfrak{W}( \cdot | v^{\textnormal{out}}, v^{\textnormal{in}})$ on sliding paths $\gamma$ is defined by \begin{equation} \label{SlidingWeightFormula} \mathfrak{W}( {\gamma} |v^{\textnormal{out}}, v^{\textnormal{in}} )= 
\prod_{\textit{\textbf{e}} \in \textnormal{\textbf{E}}_S[\gamma]} \textnormal{height}(\textit{\textbf{e}}) \prod_{\textbf{\textit{e}} \in \mathbf{E}_J^-[\gamma]} \overline{V_{-\deg ( \textbf{\textit{e}})}^{\textnormal{out}}} \prod_{\textbf{\textit{e}} \in \mathbf{E}_J^+[\gamma]} {V_{\deg ( \textbf{\textit{e}})}^{\textnormal{in}}}.\end{equation} \end{definition}

\noindent If a sliding path $\gamma$ has no slides, $\gamma$ is a Szeg\H{o} path and (\ref{SlidingWeightFormula}) reduces to (\ref{SzegoWeightFormula}).  The weight $\mathfrak{W}(\gamma |v^{\textnormal{out}}, v^{\textnormal{in}})$ determines the following $(v^{\textnormal{out}}, v^{\textnormal{in}})$-weighted enumeration problem for sliding paths.

\begin{problem} \label{ProblemSlidingPathEnumeration} Fix two specializations $v^{\textnormal{out}}=\{V_k^{\textnormal{out}}\}_{k=1}^{\infty}$, $v^{\textnormal{in}}=\{V_k^{\textnormal{in}}\}_{k=1}^{\infty}$.  Let $\Gamma_{1,0,m}(\ell)$ be the infinite set of sliding paths $\gamma$ of length $\ell$ with $m$ slides and consider the weight $\mathfrak{W}(\gamma |v^{\textnormal{out}}, v^{\textnormal{in}})$ in \textnormal{(\ref{SlidingWeightFormula})}.  Determine \begin{equation} \label{SlidingPathWFunction} W_{1,0,m}(\ell | v^{\textnormal{out}}, v^{\textnormal{in}}):= \sum_{\gamma \in \Gamma_{1,0,m}(\ell)} \mathfrak{W}(\gamma |v^{\textnormal{out}}, v^{\textnormal{in}}). \end{equation} 
\end{problem}

\subsection{Dispersive action profiles from sliding paths} \label{SUBSECDispersiveActionProfiles} 
\noindent In the case $v^{\textnormal{out}} = v^{\textnormal{in}} = v$, Problem [\ref{ProblemSlidingPathEnumeration}] has a solution that leads to the definition of dispersive action profiles $\mathbf{f}(c | v; \ebar)$.

\begin{proposition} \label{PropositionSlidingPaths} Choose $v=\{V_k\}_{k=1}^{\infty}$ so \textnormal{(\ref{Symbol})} satisfies $||v||_{\infty} = \sup_{x} |v(x)| <\infty$.  Fix $\textcolor{black}{\ebar<0}$.  There exists some $n_{*}(v) \in \{0,1,2,3,\ldots\} \cup \{ \infty\}$ and $2n_{*}(v)+1$ real numbers \begin{equation} \label{DispersiveActionProfileExtrema} \cdots <S_{i}^{\uparrow}(v; \ebar)< S_{i}^{\downarrow}(v; \ebar) < S_{i-1}^{\uparrow} (v; \ebar) \cdots < S_1^{\uparrow}(v; \ebar) < S_1^{\downarrow}(v; \ebar) < S_0^{\uparrow} (v; \ebar) \end{equation} \noindent depending on $v$ and $\ebar$ so that for any $u \in \mathbb{C} \setminus \mathbb{R}$, \textnormal{(\ref{SlidingPathWFunction})} for sliding paths \textcolor{black}{is determined by} \begin{equation} \label{SlidingPathTUpFORMULA}  \lim_{\ell' \rightarrow \infty} \sum_{\ell=0}^{\ell'} u^{- \ell -1} \sum_{m=0}^{\ell} W_{1,0,m}(\ell | v,v)  \textcolor{black}{\ebar}^m \ = \frac{1}{u - S_0^{\uparrow}(v; \ebar)} \prod_{i=1}^{n_{*}(v)} \frac{u - S_{i}^{\downarrow}( v; \ebar)}{u - S_{i}^{\uparrow}(v; \ebar)} \end{equation} \noindent a meromorphic function in $u$ \textcolor{black}{with} simple interlacing zeroes \textcolor{black}{$S_i^{\downarrow}(v; \ebar)$} and poles \textcolor{black}{$S_i^{\uparrow}(v; \ebar)$ in} \textnormal{(\ref{DispersiveActionProfileExtrema})}.
\end{proposition} \begin{itemize}
\item \textit{Proof of} \textnormal{Proposition [\ref{PropositionSlidingPaths}]}.  We have to show (\ref{SlidingPathTUpFORMULA}).  This equality has been verified in \textsection 3.3 in \cite{Moll1}. To see this, \textcolor{black}{set $V_0 = 0$ and replace $\ebar$ by $- \ebar$ in formula (1.3) in \cite{Moll1} to get} \begin{equation}\label{LaxMatrix} L_{\bullet}(v; \ebar) = \begin{bmatrix} 
0 & \overline{V_{1}} & \overline{V_{2}} & \overline{V_{3}} & \cdots  \\ 
V_{1} & \textcolor{black}{\ebar}& \overline{V_{1}} & \overline{V_{2}} & \ddots  \\
 V_{2} & V_{1} &  2 \textcolor{black}{\ebar}& \overline{V_{1}} & \ddots \\
  V_{3} & V_2 & V_1 & 3 \textcolor{black}{\ebar} & \ddots \\
  \vdots & \vdots & \ddots & \ddots & \ddots \end{bmatrix}. \end{equation}  \noindent Using $\psi_j$ and $\langle \cdot, \cdot \rangle$ in (\ref{SzegoThisThing}), the definitions imply \begin{equation} \label{SlidingThisThing} \sum_{m=0}^{\ell} W_{1,0,m}(\ell | v,v) \textcolor{black}{\ebar^m}  = \langle \psi_0, L_{\bullet}(v;\ebar)^{\ell} \psi_0 \rangle.\end{equation} Using (\ref{SlidingThisThing}) and essential self-adjointness of (\ref{LaxMatrix}) from \cite{Moll1}, the limit in (\ref{SlidingPathTUpFORMULA}) exists and is \begin{equation} \label{SlidingTheLove}\mathbf{T}^{\uparrow}(u | v;\ebar) := \langle \psi_0, \frac{1}{u - L_{\bullet}(v; \ebar)} \psi_0 \rangle. \end{equation} Finally, \textnormal{Proposition 1.2.1} of \textnormal{\cite{Moll1}} says (\ref{SlidingTheLove}) is the right-hand side of (\ref{SlidingPathTUpFORMULA}) if $||v||_{\infty}<\infty$.  $\square$
\end{itemize}
\noindent Proposition [\ref{PropositionSlidingPaths}] from \cite{Moll1} was independently found by G\'{e}rard-Kappeler \cite{GerardKappeler2019} as discussed in \cite{Moll2, Moll1}.  We now define ``dispersive action profiles'' \textcolor{black}{whose existence and uniqueness is due to \cite{Moll1}.}

\begin{definition} \label{DefinitionDispersiveActionProfiles} Let $v=\{V_k\}_{k=1}^{\infty}$ be a specialization so $||v||_{\infty} = \sup_x |v(x)| <\infty$ holds for $v(x)$ in \textnormal{(\ref{Symbol})}.  The dispersive action profile $\mathbf{f}(c | v;\ebar)$ is the \textcolor{black}{unique} function of $c \in \mathbb{R}$ characterized by the conditions that $\mathbf{f}( c | v; \ebar) \sim |c|$ as $c \rightarrow \pm \infty$ and, for all $u \in \mathbb{C} \setminus \mathbb{R}$, $\mathbf{T}^{\uparrow}(u | v; \ebar)$ in \textnormal{(\ref{SlidingTheLove})} is both\begin{equation} \label{DefinitionDispersiveActionProfilesFORMULA}  \lim_{\ell' \rightarrow \infty} \sum_{\ell=0}^{\ell'} u^{- \ell -1} \sum_{m=0}^{\ell} W_{1,0,m}(\ell | v,v)\textcolor{black}{\varepsilon^m}   = \textnormal{exp} \Bigg ( \int_{-\infty}^{+\infty} \log \Bigg [ \frac{1}{u-c} \Bigg ] \cdot \tfrac{1}{2} \mathbf{f}''(c | v; \ebar) dc \Bigg ) \end{equation} \noindent where $W_{1,0,m}(\ell | v,v) $ are the $(v,v)$-weighted sums over sliding paths in \textnormal{(\ref{SlidingPathWFunction})}. \end{definition}

\noindent Dispersive action profiles were introduced in Definition 1.2.2 in \cite{Moll1} and were implicitly studied in G\'{e}rard-Kappeler \cite{GerardKappeler2019} as shown in \cite{Moll2, Moll1}.  In \cite{Moll1}, many properties of dispersive action profiles $\mathbf{f}(c | v;\ebar)$ are established.  By Proposition [\ref{PropositionSlidingPaths}], dispersive action profiles are piecewise-linear functions with slopes $\pm 1$, $n_{*}(v)+1$ local minima $S_i^{\uparrow}(v; \ebar)$, and $n_{*}(v)$ local maxima $S_{i+1}^{\downarrow}(v; \ebar)$ from (\ref{DispersiveActionProfileExtrema}).  In light of further results in \cite{GerardKappeler2019, Moll2}, \textcolor{black}{for $\ebar < 0$} these local extrema have no accumulation points except possibly $- \infty$ if $n_{*}(v) = \infty$, are bounded above by $||v||_{\infty} = \sup_x |v(x)|$ for $v(x)$ in (\ref{Symbol}), and satisfy $S_{i-1}^{\uparrow}(v;\ebar) - S_{i}^{\downarrow}(v; \ebar) \in \textcolor{black}{| \ebar|} \mathbb{Z}_+$.  See Figure [\ref{LLNFIG2020}] for a depiction of a dispersive action profile.  For further detail, see Figure 1 in \cite{Moll1}.  Note that if $n_{*}(v) = \infty$ in Proposition [\ref{PropositionSlidingPaths}], the slope of the curve $\mathbf{f}(c | v; \ebar)$ in Figure [\ref{LLNFIG2020}] changes infinitely-many times in the direction $c \rightarrow - \infty$.  
\subsection{Example: Motzkin paths and Nekrasov-Pestun-Shatashvili profile} \label{SUBSECMotzkinNPS} 
\begin{figure}[htb]
\centering
\includegraphics[width=0.5 \textwidth]{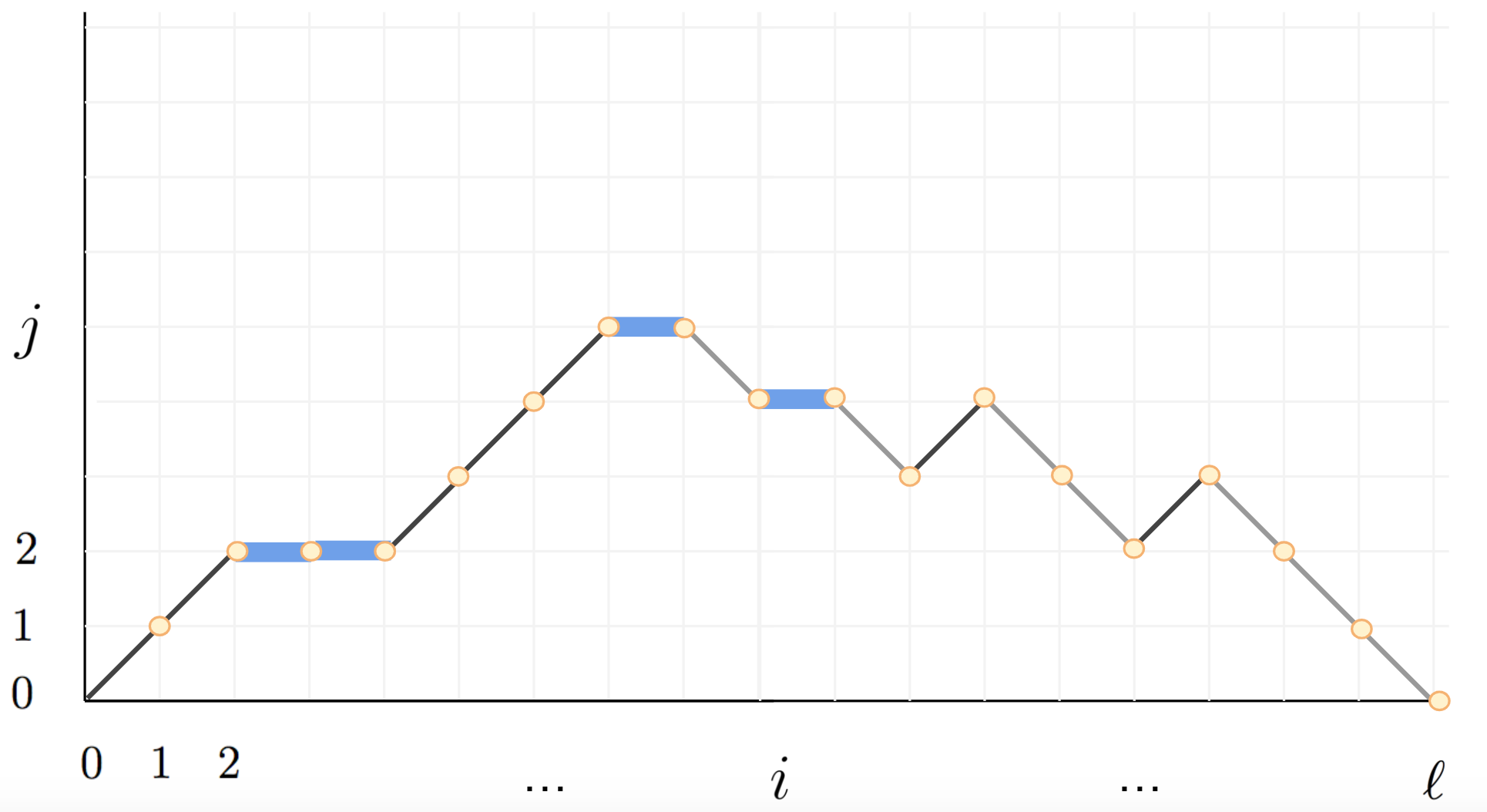}
\caption{A Motzkin path of length $\ell=18$ with $m=4$ slides.}
\label{MotzkinPath2020}
\end{figure}

Consider the case $v_{\textnormal{PL}}$ of a Plancherel specialization $v_{\textnormal{PL}} = \{1,0,0,\ldots\}$ in (\ref{PlancherelSpecialization}).  In this case, the weight (\ref{SlidingWeightFormula}) is supported on Motzkin paths but is not uniform.  For example, $\mathfrak{W}(\gamma |v_{\textnormal{PL}}, v_{\textnormal{PL}}) = 2 \cdot 2 \cdot 5 \cdot 4 = 80$ in Figure [\ref{MotzkinPath2020}].

\noindent By Proposition 7.1.1 of \cite{Moll1}, for $v_{\textnormal{PL}}(x) = 2 \cos x$ in (\ref{Symbol}), $\mathbf{T}^{\uparrow}(u | v_{\textnormal{PL}};\ebar)$ in (\ref{SlidingTheLove}) and (\ref{SlidingPathTUpFORMULA}) solves \begin{equation} \label{MotzkinQuadraticTing} \mathbf{T}^{\uparrow}(u + \ebar | v_{\textnormal{PL}};\ebar) + \frac{1}{\mathbf{T}^{\uparrow}(u | v_{\textnormal{PL}};\ebar)} =u\end{equation} \noindent a difference equation derived in \textsection 4 of Poghossian \cite{Poghossian} now understood as the Nekrasov-Shatashvili limit of Nekrasov's non-perturbative Dyson-Schwinger equations \cite{NekYI} in the case of a particular pure abelian gauge theory.  In this case, the Nekrasov-Okounkov measures on partitions \cite{NekOk} are the Poissonized Jack-Plancherel measures \cite{Ke4} as discussed in \textsection 3.2 of \cite{NekOk}.  As such, $\mathbf{f}(c | v_{\textnormal{PL}}; \ebar)$ is one of the limit shapes of Nekrasov-Okounkov measures by Nekrasov-Pestun-Shatashvili \cite{NekPesSha}.  This dispersive action profile $\mathbf{f}(c | v_{\textnormal{PL}}; \ebar)$ was independently discovered by Do{\l}{\k{e}}ga-{\'S}niady \cite{DoSni}.




\section{Ribbon paths and Jack measures on profiles} \label{SECRibbonAnisotropic}

\noindent In this section we prove Theorem [\ref{Theorem3AOE}], our main technical result for Jack measures $M(v^{\textnormal{out}}, v^{\textnormal{in}})$ \textcolor{black}{with arbitrary $v^{\textnormal{out}}, v^{\textnormal{in}}$ satisfying the assumptions in Definition [\ref{DefinitionJackMeasures}]}.  Our proof of Theorem [\ref{Theorem3AOE}] is based on the notion of ``connected ribbon paths.''  In \textsection [\ref{SUBSECRibbonPaths}] we define ribbon paths $\boldsymbol{\vec{\gamma}}$ on $n$ ``sites'' of lengths $\ell_1, \ldots, \ell_n$ as a generalization of ordered lists of $n$ sliding paths $\gamma_1, \ldots, \gamma_n$ that can have ``pairings'' between jumps of opposite degrees.  We then say that a ribbon path $\boldsymbol{\vec{\gamma}}$ is ``connected'' if its sites and pairings define a connected graph.  In \textsection [\ref{SUBSECWeightedEnumerationRibbonPaths}], we pose a $(v^{\textnormal{out}}, v^{\textnormal{in}})$-weighted enumeration problem for connected ribbon paths on $n$ sites of lengths $\ell_1, \ldots, \ell_n$ with $n-1+g$ pairings and $m$ slides in Problem [\ref{ProblemRibbonPathEnumeration}].  In \textsection [\ref{SUBSECJackMeasuresProfilesRibbonPaths}], we describe the solution $W_{n,g,m}(\ell_1, \ldots, \ell_n | v^{\textnormal{out}}, v^{\textnormal{in}})$ of this problem in Theorem [\ref{Theorem3AOE}] as the coefficient of $\hbar^{n-1+g} \ebar^{m}$ in the polynomial expansion of the joint cumulants $\kappa_n$ from (\ref{KeyRelation}).  In \textsection [\ref{SUBSECProofExpectedSizeProposition}], we illustrate our Theorem [\ref{Theorem3AOE}] relating Jack measures and ribbon paths in the case $n=1$, $\ell=2$: by computing $W_{1,g,m}(2|v,v)$, we prove Proposition [\ref{PropositionSizeExpectation}].\\
\\
\noindent In \textsection [\ref{SUBSECRibbonPathsComments}] we discuss our motivation for introducing ribbon paths and the larger context of our Theorem [\ref{Theorem3AOE}].  In \textsection [\ref{APPENDIXRibbonPaths}], we show that the joint cumulants in our Theorem [\ref{Theorem3AOE}] appear in Chapuy-Do{\l}{\k{e}}ga \cite{ChapuyDolega2020} and relate ribbon paths to the ribbon graphs on non-oriented surfaces in \cite{ChapuyDolega2020}.

\subsection{Ribbon paths} \label{SUBSECRibbonPaths} To define ribbon paths, we first need the notion of a ``pairing''.

\begin{definition} \label{DefinitionPairings} Let $\vec{\gamma} = (\gamma_1, \ldots, \gamma_n)$ denote an ordered sequence of $n$ lattice paths $\gamma_a$ in $\mathbb{Z}_{\bullet}^2$ indexed by $a=1,2,\ldots, n$ of lengths $\ell_1, \ldots, \ell_n$.  A pairing $\textit{\textbf{p}}$ is the data of two steps $(\textbf{\textit{e}}_{i;a}, \textbf{\textit{e}}_{i'; a'})$, where $\textbf{\textit{e}}_{i;a}$ is the $i$th step in $\gamma_a$ and $\textbf{\textit{e}}_{i';a'}$ is the $i'$th step in $\gamma_{a'}$, satisfying
\begin{itemize}
\item Degree condition: $\deg (\textbf{\textit{e}}_{i;a})= -k$ and $\deg (\textbf{\textit{e}}_{i';a'}) = +k$ for some $k=1,2,3,\ldots$
\item Ordering condition: either $a < a'$ or $a = a'$ and $i < i'$
\end{itemize}
\noindent We say that such a pairing $\textit{\textbf{p}} =(\textbf{\textit{e}}, \textbf{\textit{e}'}) $ is of size $k$ and write $\textnormal{size}(\textbf{p})  =\deg(\textbf{\textit{e}'}) = k$.
\end{definition}

\begin{figure}[htb]
\centering
\includegraphics[width=0.65 \textwidth]{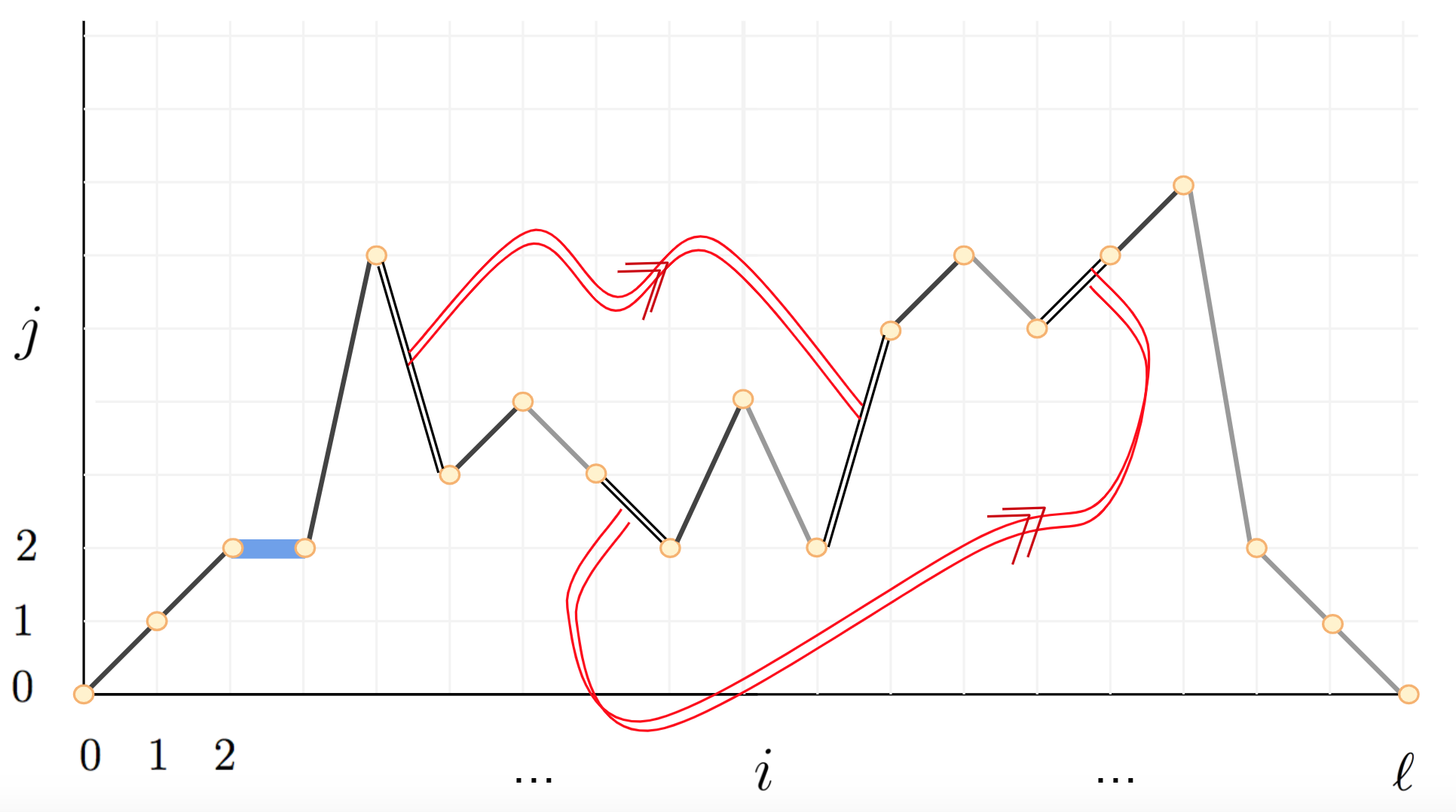}
\caption{A ribbon path $\boldsymbol{\vec{\gamma}} = (\gamma_1; \textit{\textbf{p}}_1, \textit{\textbf{p}}_2)$ on $n=1$ site of length $\ell_1 = 18$ with $m=1$ slide $\textbf{\textit{e}}$ of $\textnormal{height}(\textbf{\textit{e}})=2$, $g=2$ pairings $\textit{\textbf{p}}_1, \textbf{\textit{p}}_2$ of $\textnormal{size}(\textit{\textbf{p}}_1)=3$, $\textnormal{size}(\textbf{\textit{p}}_2)=1$.}
\label{RibbonPath2020}
\end{figure} 

\noindent A pairing can only consist of a down step by $-k$ and an up step by $+k$ in which the down step appears before the up step in the ambient ordering.  We may now define ribbon paths.

\begin{definition} \label{DefinitionRibbonPaths} A ribbon path $\boldsymbol{\vec{\gamma}}$ on $n$ sites of lengths $\ell_1, \ldots, \ell_n$ is an ordered sequence of $n$ sliding paths $\gamma_1, \ldots, \gamma_n$ from \textnormal{Definition [\ref{DefinitionSlidingPaths}]} together with a number $q \geq 0$ of pairings $\textit{\textbf{p}}_1, \ldots, \textit{\textbf{p}}_q$ as in \textnormal{Definition [\ref{DefinitionPairings}]} so that any $\textit{\textbf{e}} \in \textnormal{\textbf{E}} [\boldsymbol{\vec{\gamma}}]:=\textnormal{\textbf{E}}[\gamma_1] \cup \cdots \cup \textnormal{\textbf{E}}[\gamma_n]$ participates in at most one pairing.  Write $\boldsymbol{\vec{\gamma}} = (\gamma_1, \ldots, \gamma_n ; \textit{\textbf{p}}_1, \ldots, \textit{\textbf{p}}_q)$, $\textit{\textbf{p}}_i \in \boldsymbol{\vec{\gamma}}$, and $\mathbf{E}_{PJ}[\boldsymbol{\vec{\gamma}}]$ for the subset of steps in pairings of $\boldsymbol{\vec{\gamma}}$.
\end{definition}

\begin{figure}[htb]
\centering
\includegraphics[width=0.85 \textwidth]{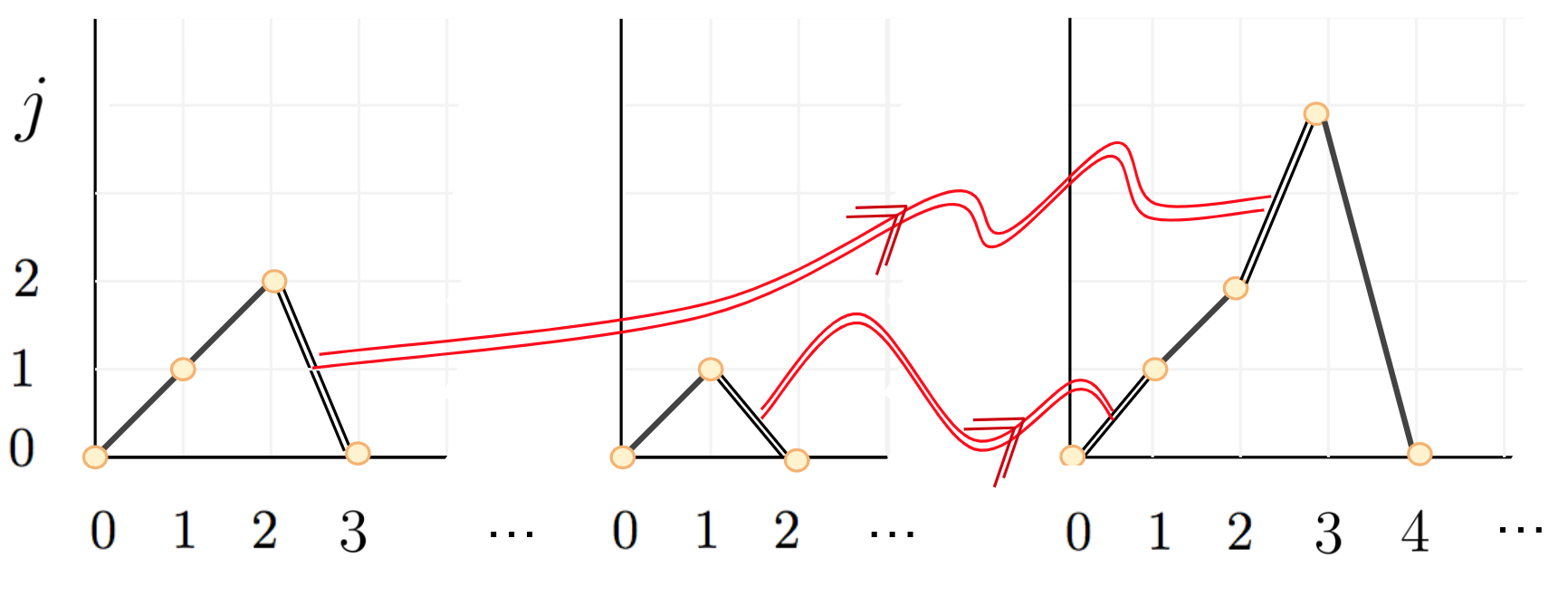}
\caption{A ribbon path on $n=3$ sites of lengths $\ell_1 = 3$, $\ell_2 = 2$, $\ell_3 = 4$ with $m=0$ slides and $q=2$ pairings $\textit{\textbf{p}}_1, \textit{\textbf{p}}_2$ of sizes $k(\textbf{\textit{p}}_1) = 2$ and $k(\textbf{\textit{p}}_2)=1$.}
\label{ConnectedRibbonPath2020}
\end{figure}

\noindent For ribbon paths $\boldsymbol{\vec{\gamma}}$ on $n=1$ and $n=3$ sites, see Figures [\ref{RibbonPath2020}] and [\ref{ConnectedRibbonPath2020}], respectively.  As is automatic from Definition [\ref{DefinitionRibbonPaths}], ribbon paths generalize both Szeg\H{o} paths and sliding paths: 

\begin{proposition} \label{DuhReduceDuhProposition} Let $\boldsymbol{\vec{\gamma}}$ be a ribbon path on $n$ sites.
\begin{enumerate}
\item If $\boldsymbol{\vec{\gamma}}$ has $0$ pairings, $\boldsymbol{\vec{\gamma}}$ is a disjoint union of $n$ sliding paths.
\item If $\boldsymbol{\vec{\gamma}}$ has $0$ pairings and $0$ slides, $\boldsymbol{\vec{\gamma}}$ is a disjoint union of $n$ Szeg\H{o} paths.
\end{enumerate}
\end{proposition}
\noindent By the degree condition in Definition [\ref{DefinitionPairings}], slides cannot participate in pairings, therefore \begin{equation} \label{RibbonPathStepSetDecomposition} \textnormal{\textbf{E}}[\boldsymbol{\vec{\gamma}}]= \textnormal{\textbf{E}}_{PJ}[ \boldsymbol{\vec{\gamma}}] \cup  \textnormal{\textbf{E}}_S[\boldsymbol{\vec{\gamma}}] \cup \textnormal{\textbf{E}}_{UJ}^+[\boldsymbol{\vec{\gamma}}] \cup  \textnormal{\textbf{E}}_{UJ}^-[\boldsymbol{\vec{\gamma}}]\end{equation} \noindent the set of steps in a ribbon path is the disjoint union of the set $\textnormal{\textbf{E}}_{PJ} [ \boldsymbol{\vec{\gamma}}]$ of paired jumps, the set of slides $\textnormal{\textbf{E}}_S[\boldsymbol{\vec{\gamma}}]$, and the sets $\textnormal{\textbf{E}}_{UJ}^{\pm}[ \boldsymbol{\vec{\gamma}}]$ of unpaired jumps of degrees in $\mathbb{Z}_{\pm}$.  Using (\ref{RibbonPathStepSetDecomposition}), we associate two partitions to every ribbon path $\boldsymbol{\vec{\gamma}}$ which we call the ``{unpaired jump profiles}'' of $\boldsymbol{\vec{\gamma}}$.
\begin{definition} \label{DefinitionUnpairedJumpProfiles} Let $\boldsymbol{\vec{\gamma}}$ be a ribbon path.  The unpaired jump profiles of $\boldsymbol{\vec{\gamma}}$ are the two partitions $\boldsymbol{\mu}^{\pm} [ \boldsymbol{\vec{\gamma}}]$ defined from \textnormal{(\ref{RibbonPathStepSetDecomposition})} by \begin{equation} \# \{i : \boldsymbol{\mu}^{{\pm}}_i [ \boldsymbol{\vec{\gamma}}] =k \} = \# \{ \mathbf{e} \in \mathbf{E}_{UJ}^{\pm} [ \boldsymbol{\vec{\gamma}}] \ : \ \deg (\mathbf{e}) = \pm k\} .\end{equation}  
\end{definition}

\noindent For example, in the ribbon path $\boldsymbol{\vec{\gamma}}$ in Figure [\ref{RibbonPath2020}], reading from left to right the unpaired jumps of positive degree have degrees $1, 1, 4, 1, 2, 1, 1$ while the unpaired jumps of negative degree have degrees $-1, -2, -1, -5, -1, -1$ and hence, sorting these in weakly decreasing order, the unpaired jump profiles of the ribbon path $\boldsymbol{\vec{\gamma}}$ in Figure [\ref{RibbonPath2020}] are \begin{eqnarray} \boldsymbol{\mu}^{+} [ \boldsymbol{\vec{\gamma}}] &=& (  1 \leq 1 \leq 1 \leq 1 \leq 1 \leq 2 \leq 4)  \\
\boldsymbol{\mu}^{-} [ \boldsymbol{\vec{\gamma}}] &=& (  1 \leq 1 \leq 1 \leq 1 \leq 2 \leq 5) .
\end{eqnarray}

\noindent These unpaired jump profiles are partitions of the same size $|  \boldsymbol{\mu}^{{+}}[ \boldsymbol{\vec{\gamma}}] | =  | \boldsymbol{\mu}^{{-}}[ \boldsymbol{\vec{\gamma}}] | = 11$. \begin{definition} \label{DefinitionSizeRibbonPath} The size of a ribbon path is the size of either unpaired jump profile \begin{equation} \mathbf{d} [ \boldsymbol{\vec{\gamma}}] = |  \boldsymbol{\mu}^{{+}}[ \boldsymbol{\vec{\gamma}}] | =  | \boldsymbol{\mu}^{{-}}[ \boldsymbol{\vec{\gamma}}] |. \end{equation} \end{definition} 

\noindent The size of a ribbon path is well-defined since (i) each of the sliding paths in $\boldsymbol{\vec{\gamma}}$ starts and ends at height $0$ and (ii) any pairing in $\boldsymbol{\vec{\gamma}}$ pairs steps of opposite degrees implies that $|\boldsymbol{\mu}^{{+}}[ \boldsymbol{\vec{\gamma}}]| =| \boldsymbol{\mu}^{{-}}[ \boldsymbol{\vec{\gamma}}]|$.  As an additional example, the unpaired jump profiles of the Szeg\H{o} path $\gamma$ in Figure [\ref{SzegoPath2020}] are  \begin{eqnarray} \boldsymbol{\mu}^{+} [ \boldsymbol{\vec{\gamma}}] &=& ( 1 \leq 1 \leq 1 \leq 2 \leq 3 \leq 4 \leq 4 \leq 5 \leq 5 \leq 7)  \\
\boldsymbol{\mu}^{-} [ \boldsymbol{\vec{\gamma}}] &=& (  2 \leq 3 \leq 4 \leq 4 \leq 4 \leq 5 \leq 5 \leq 6) .
\end{eqnarray}
\noindent and the size of this Szeg\H{o} path is $\mathbf{d} [\boldsymbol{\vec{\gamma}}]= 33$.  For any Szeg\H{o} path, its size is half the metric length of the piecewise-linear trajectory of this path if one chooses to connect the vertices of the path by a linear interpolation as in Figure [\ref{SzegoPath2020}] \textcolor{black}{since all of its jumps are unpaired}.  However, in this paper we reserve the term ``length'' for the metric length \textcolor{black}{$\ell$} of the ``site'' $[0, \ell]$, so that in Figure [\ref{SzegoPath2020}] we have a ribbon path of size $33$ with $0$ pairings and $0$ slides on a single site of length $18$.  \textcolor{black}{In \textsection [\ref{AppendixB3}], we relate the size $\mathbf{d} [ \boldsymbol{\vec{\gamma}}]$ of ribbon paths to the size $| \mathbf{M}|$ of ribbon graphs' constellations in Chapuy-Do{\l}{\k{e}}ga \cite{ChapuyDolega2020}.}  \\
\\
\noindent We now define connectivity for ribbon paths.

\begin{definition} \label{DefinitionReducedGraph} The reduced graph of a ribbon path $\boldsymbol{\vec{\gamma}} = (\gamma_1, \ldots, \gamma_n; \textit{\textbf{p}}_1, \ldots, \textit{\textbf{p}}_q)$ is the graph with $n$ vertices $a=1, 2, \ldots, n$ and edges $a \rightarrow a'$ for every pairing $\textbf{\textit{p}} = ( \textbf{\textit{e}}, \textbf{\textit{e'}})$ of $\textbf{\textit{e}} \in \textnormal{\textbf{E}}[\gamma_a]$ and $\textnormal{\textbf{E}}[\gamma_{a'}]$. \end{definition}

\begin{definition} \label{DefinitionConnectedRibbonPaths} A ribbon path $\vec{\gamma}$ is connected if its reduced graph is connected. \end{definition}

\begin{figure}[htb]
\centering
\includegraphics[width=0.85 \textwidth]{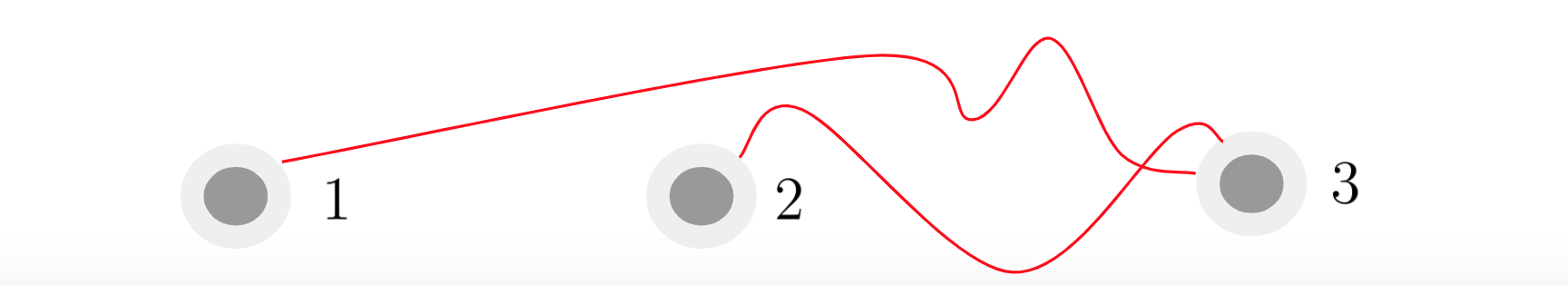}
\caption{The reduced graph of the connected ribbon path in Figure [\ref{ConnectedRibbonPath2020}] .}
\label{ReducedGraphConnectedRibbonPath2020}
\end{figure} 

\noindent The following special case is automatic from Definition [\ref{DefinitionConnectedRibbonPaths}].

\begin{proposition} 
\noindent \label{RibbonPathConnectedPropositionDuh}All ribbon paths $\boldsymbol{\vec{\gamma}}$ on $n=1$ site are connected even if $\boldsymbol{\vec{\gamma}}$ has no pairings.

\end{proposition}

\subsection{Weighted enumeration of ribbon paths} \label{SUBSECWeightedEnumerationRibbonPaths} We now define a $\mathbb{C}$-weight $\mathfrak{W}$ for ribbon paths.
 \begin{definition} \label{DefinitionRibbonPathsWeight} Let $v^{\textnormal{out}}=\{V_k^{\textnormal{out}}\}_{k=1}^{\infty}$, $v^{\textnormal{in}}=\{V_k^{\textnormal{in}}\}_{k=1}^{\infty}$ be two specializations.  Write $\overline{V}$ for the complex conjugate of $V \in \mathbb{C}$.  The $(v^{\textnormal{out}}, v^{\textnormal{in}})$-weight $\mathfrak{W}( \cdot | v^{\textnormal{out}}, v^{\textnormal{in}})$ on ribbon paths is defined by \begin{equation} \label{RibbonWeightFormula}  \mathfrak{W}( \boldsymbol{\vec{\gamma}} |v^{\textnormal{out}}, v^{\textnormal{in}} )= \prod_{ \textbf{p} \in \boldsymbol{\vec{\gamma}}}  \textnormal{size}(\textbf{p}) \prod_{\textit{\textbf{e}} \in \textnormal{\textbf{E}}_S[\boldsymbol{\vec{\gamma}}]} \textnormal{height}(\textit{\textbf{e}}) \prod_{\textbf{\textit{e}} \in \mathbf{E}_{UJ}^-[\boldsymbol{\vec{\gamma}}]} \overline{V_{-\deg ( \textbf{\textit{e}})}^{\textnormal{out}}} \prod_{\textbf{\textit{e}} \in \mathbf{E}_{UJ}^+[\boldsymbol{\vec{\gamma}}]} {V_{\deg ( \textbf{\textit{e}})}^{\textnormal{in}}}.\end{equation} \end{definition}
  
 \noindent For example, $\mathfrak{W}(\boldsymbol{\vec{\gamma}}| v^{\textnormal{out}}, v^{\textnormal{in}}) = V_1^{\textnormal{in}} V_1^{\textnormal{in}} \cdot {2} \cdot V_1^{\textnormal{in}} \cdot  {1} \cdot V_{1}^{\textnormal{in}} \overline{V^{\textnormal{out}}_{4}} =  2 (V_1^{\textnormal{in}})^4 \overline{V_{4}^{\textnormal{out}}}$ for $\boldsymbol{\vec{\gamma}}$ in Figure [\ref{ConnectedRibbonPath2020}].  In general, (\ref{RibbonWeightFormula}) are multinomials in $\overline{V_{k}^{\textnormal{out}}}, V_k^{\textnormal{in}}$ with positive integer coefficients.  If a ribbon path $\boldsymbol{\vec{\gamma}}$ has no pairings, $\boldsymbol{\vec{\gamma}}$ is a disjoint union of $n$ sliding paths and (\ref{SlidingWeightFormula}) reduces to the product of weights (\ref{SzegoWeightFormula}).  $\mathfrak{W}(\boldsymbol{\vec{\gamma}} |v^{\textnormal{out}}, v^{\textnormal{in}})$ determines the following $(v^{\textnormal{out}}, v^{\textnormal{in}})$-weighted enumeration problem for ribbon paths.

\begin{problem} \label{ProblemRibbonPathEnumeration}  Fix two specializations $v^{\textnormal{out}}=\{V_k^{\textnormal{out}}\}_{k=1}^{\infty}$, $v^{\textnormal{in}}=\{V_k^{\textnormal{in}}\}_{k=1}^{\infty}$.  Consider the complex weight $\mathfrak{W}(\boldsymbol{\vec{\gamma}} |v^{\textnormal{out}}, v^{\textnormal{in}})$ on ribbon paths in \textnormal{(\ref{SlidingWeightFormula})}.
\begin{enumerate}
\item Let $\Delta_{n,q,m}(\ell_1, \ldots, \ell_n)$ be the infinite set of ribbon paths $\boldsymbol{\vec{\gamma}}$ on $n$ sites of lengths $\ell_1, \ldots, \ell_n$ with $q$ pairings and $m$ slides.  Determine \begin{equation} \label{RibbonPathYFunction} Y_{n,q,m}(\ell_1, \ldots, \ell_n | v^{\textnormal{out}}, v^{\textnormal{in}}):= \sum_{\boldsymbol{\vec{\gamma}} \in \Delta_{n,q,m}(\ell_1, \ldots, \ell_n)} \mathfrak{W}(\boldsymbol{\vec{\gamma}} |v^{\textnormal{out}}, v^{\textnormal{in}}). \end{equation} 
\item Let $\Gamma_{n,g,m}(\ell_1, \ldots, \ell_n)$ be the infinite set of connected ribbon paths $\boldsymbol{\vec{\gamma}}$ on $n$ sites of lengths $\ell_1, \ldots, \ell_n$ with $n-1+g$ pairings and $m$ slides.  Determine \begin{equation} \label{RibbonPathWFunction} W_{n,g,m}(\ell_1, \ldots, \ell_n | v^{\textnormal{out}}, v^{\textnormal{in}}):= \sum_{\boldsymbol{\vec{\gamma}} \in \Gamma_{n,g,m}(\ell_1, \ldots, \ell_n)} \mathfrak{W}(\boldsymbol{\vec{\gamma}} |v^{\textnormal{out}}, v^{\textnormal{in}}). \end{equation} 
\end{enumerate}
\end{problem}

\subsection{Jack measures on profiles from ribbon paths} \label{SUBSECJackMeasuresProfilesRibbonPaths} 
\noindent Under suitable assumptions on $v^{\textnormal{out}}, v^{\textnormal{in}}$, we can describe the solution of Problem [\ref{ProblemRibbonPathEnumeration}].  This is our main technical result for Jack measures.
\begin{theorem} \label{Theorem3AOE} Fix $\ebar \in \mathbb{R}$, $\hbar>0$, and two specializations $v^{\textnormal{out}} = \{V_k^{\textnormal{out}}\}_{k=1}^{\infty}$, $v^{\textnormal{in}} = \{V_k^{\textnormal{in}}\}_{k=1}^{\infty}$ satisfying $\max (|V_k^{\textnormal{out}}|, | V_k^{\textnormal{in}}|) \leq A r^k$ for some $A>0$ and $0<r<1$.  Let $f_{\lambda}( c | \ebar, \hbar)$ be the random anisotropic partition profile sampled from the Jack measure $M(v^{\textnormal{out}}, v^{\textnormal{in}})$.  Consider any $n=1,2,3,\ldots$ and any $\vec{\ell}= (\ell_1, \ldots, \ell_n) \in \mathbb{Z}_+^n$.  Write $|| \vec{\ell}||_1 = \ell_1 + \cdots + \ell_n$.  \begin{enumerate}
\item The joint moments of the random variables $T_{\ell}^{\uparrow}(\ebar, \hbar)$ in \textnormal{(\ref{TransitionStatistics})} are polynomials in $\ebar$ and $\hbar$
\begin{equation} \label{JointMomentsAOE} \mathbb{E} \big [ T_{\ell_1}^{\uparrow} (\ebar,\hbar) \cdots T_{\ell_n}^{\uparrow} (\ebar, \hbar) \big ] = \sum_{q=0}^{\tfrac{1}{2} || \vec{\ell}||_1 } \sum_{m=0}^{|| \vec{\ell}||_1} Y_{n,q,m}(\ell_1, \ldots, \ell_n | v^{\textnormal{out}}, v^{\textnormal{in}})  \hbar^q  \textcolor{black}{\ebar^m}  \end{equation} \noindent with finite coefficients $Y_{n,q,m}(\ell_1, \ldots, \ell_n | v^{\textnormal{out}}, v^{\textnormal{in}})$ counting ribbon paths in \textnormal{(\ref{RibbonPathYFunction})}.\\
\item The joint cumulants of the random variables $T_{\ell}^{\uparrow}(\ebar, \hbar)$ in \textnormal{(\ref{TransitionStatistics})} are polynomials in $\ebar$ and $\hbar$
\begin{equation} \label{JointCumulantsAOE} \kappa_n \big ( T_{\ell_1}^{\uparrow} (\ebar, \hbar), \ldots, T_{\ell_n}^{\uparrow} (\ebar, \hbar) \big ) = \sum_{g=0}^{\tfrac{1}{2} || \vec{\ell}||_1 -n+1} \sum_{m=0}^{|| \vec{\ell}||_1}  W_{n,g,m}(\ell_1, \ldots, \ell_n | v^{\textnormal{out}}, v^{\textnormal{in}})  \hbar^{n-1+g} \textcolor{black}{\ebar^m}  \end{equation}
\noindent with finite coefficients $W_{n,\textcolor{black}{g},m}(\ell_1, \ldots, \ell_n | v^{\textnormal{out}}, v^{\textnormal{in}})$ counting connected ribbon paths in \textnormal{(\ref{RibbonPathWFunction})}. \end{enumerate}  \end{theorem}

\noindent \textcolor{black}{The scaling $\kappa_n \sim \hbar^{n-1}$ of the joint cumulants presented earlier in (\ref{KeyRelation}) follows from (\ref{JointCumulantsAOE}) in Regimes I and III for fixed $\ebar$ since as $\hbar \rightarrow 0$ the dominant terms in (\ref{JointCumulantsAOE}) are those with $W_{n,g,m}$ at $g=0$ and arbitrary $m$.  On the other hand, (\ref{KeyRelation}) also follows from (\ref{JointCumulantsAOE}) in Regime II with $\ebar, \hbar \rightarrow 0$ as $\ebar \sim \hbar^{1/2}$ since in this case the dominant terms in (\ref{JointCumulantsAOE}) are those with $W_{n,g,m}$ at $g=0$ and $m=0$, namely $\kappa_n \sim W_{n,0,0} \hbar^{n-1}$.  Similarly,} the expansion (\ref{LocatingRibbonPathsIntro}) is automatic from (\ref{JointCumulantsAOE}) and the definition (\ref{RibbonWeightFormula}) of the weight $\mathfrak{W}$ with
\begin{equation} \label{ConnectedTwoMuCountFormula} C_{n,g,m}(\ell_1, \ldots, \ell_n | \mu^{\textnormal{out}}, \mu^{\textnormal{in}} ) = \sum_{\boldsymbol{\vec{\gamma}} \in \Gamma_{n,g,m}(\ell_1, \ldots, \ell_n | \mu^{\textnormal{out}}, \mu^{\textnormal{in}})} \prod_{ \textbf{p} \in \boldsymbol{\vec{\gamma}}}  \textnormal{size}(\textbf{p}) \prod_{\textit{\textbf{e}} \in \textnormal{\textbf{E}}_S[\boldsymbol{\vec{\gamma}}]} \textnormal{height}(\textit{\textbf{e}})  \end{equation} where $\Gamma_{n,g,m}(\ell_1, \ldots, \ell_n | \mu^{\textnormal{out}}, \mu^{\textnormal{in}})$ is the finite set of connected ribbon paths in $\Gamma_{n,g,m}(\ell_1, \ldots, \ell_n)$ with unpaired jump profiles $\boldsymbol{\mu}^{-}[ \boldsymbol{\vec{\gamma}}] = \mu^{\textnormal{out}}$ and $\boldsymbol{\mu}^{+}[ \boldsymbol{\vec{\gamma}}] = \mu^{\textnormal{in}}$ as in Definition [\ref{DefinitionUnpairedJumpProfiles}].

\noindent \textit{Proof of \textnormal{Theorem [\ref{Theorem3AOE}]}.} We prove (\ref{JointMomentsAOE}) and (\ref{JointCumulantsAOE}) in Steps 6 and 7 of Steps 1-7 below. 
\subsubsection{\textbf{Step 1: Stanley's Cauchy kernel is a reproducing kernel}} \label{SUBSUBSECStep1} First, we recall the fact that Stanley's Cauchy kernel (\ref{CauchyIdentity}) defining Jack measures $M(v^{\textnormal{out}}, v^{\textnormal{in}})$ is a reproducing kernel for the Hilbert space $\mathcal{F}$ defined in \textsection [\ref{SUBSECJackPolynomials}] as the completion of $\C[\rho_1, \rho_2, \ldots]$ with respect to the $\hbar$-dependent inner product $\langle \cdot , \cdot \rangle_{\hbar}$.  For background on reproducing kernel Hilbert spaces, see Aronszajn \cite{AronszajnTheoryReproducingKernels}.

\begin{lemma} \label{CoherentStateLemma} Let $v = \{V_k\}_{k=1}^{\infty}$ be a specialization of $\rho_1, \rho_2, \ldots$ satisfying $\sum_{k=1}^{\infty} \frac{|V_k|^2}{k} < \infty$ and let $\widehat{\rho}_{-k} = \hbar k \frac{\partial}{\partial \rho_k}$ be the annihilation operators defined in \textnormal{(\ref{AnnihilationOperator})}.  Consider the exponential \begin{equation} \label{CoherentState} \Upsilon_v(\rho_1, \rho_2, \ldots | \hbar) = \textnormal{exp} \Bigg ( \frac{1}{\hbar} \sum_{k=1}^{\infty}\frac{ \overline{V_k} \rho_k }{k} \Bigg ). \end{equation}
\begin{enumerate}
\item $\Upsilon_v( \textcolor{black}{\rho} | \hbar)$ is an element of $\mathcal{F}$.
\item $\Upsilon_v( \textcolor{black}{\rho} | \hbar)$ is a joint eigenfunction of the annihilation operators \begin{equation} \label{HaveToLoveThisEigenvalueRelation} \widehat{\rho}_{-k} \Upsilon_v(\rho_1, \rho_2, \ldots | \hbar) = \overline{V}_k  \Upsilon_v(\rho_1, \rho_2, \ldots | \hbar) \end{equation} with eigenvalue $\overline{V_k}$ given by the complex conjugate of $V_k$.
\item $\Upsilon_v( \textcolor{black}{\rho}  | \hbar)$ is a reproducing kernel for $(\mathcal{F}, \langle \cdot, \cdot \rangle_{\hbar})$.
\end{enumerate}  \end{lemma} 

\noindent \textit{Proof of \textnormal{Lemma [\ref{CoherentStateLemma}]}}. Using the Taylor expansion for the exponential and formula (\ref{NormIntro}), \begin{equation} \label{CoherentStateNorm} || \Upsilon_v( \textcolor{black}{\rho} | \hbar) ||_{\hbar}^2 = \textnormal{exp} \Bigg ( \textcolor{black}{\frac{1}{\hbar}} \sum_{k=1}^{\infty} \frac{ |V_k|^2}{k} \Bigg ) \end{equation} 
\noindent \textcolor{black}{so $||\Upsilon_v ( \rho | \hbar) ||_{\hbar} < \infty$} \textcolor{black}{if} $\sum_{k=1}^{\infty} \frac{|V_k|^2}{k} <\infty$, proving (1).  \textcolor{black}{(2) is a direct calculation.}  \textcolor{black}{To prove (3)}, \textcolor{black}{use (2) and} $\widehat{\rho}_{\pm k}^{\dagger} = \widehat{\rho}_{\mp k}$ \textcolor{black}{for} $\widehat{\rho}_{+k}$ in (\ref{CreationOperator}) and $\widehat{\rho}_{-k}$ in (\ref{AnnihilationOperator}), for any $P( \rho_1, \rho_2, \ldots) \in \mathcal{F}$ we have \begin{eqnarray} \langle \Upsilon_v( \textcolor{black}{\rho} | \hbar) , P(\rho_1, \rho_2, \ldots ) \rangle_{\hbar} &=& \langle \Upsilon_v( \textcolor{black}{\rho} | \hbar) , P(\widehat{\rho}_1, \widehat{\rho}_2, \ldots ) \cdot 1 \rangle_{\hbar}  \\ &=&  \langle P(\widehat{\rho}_{-1}, \widehat{\rho}_{\textcolor{black}{-2}}, \ldots ) \cdot \Upsilon_v( \textcolor{black}{\rho} | \hbar), 1 \rangle_{\hbar}  \\ &=& \langle P(\overline{V}_1, \overline{V_2}, \ldots )\cdot  \Upsilon_v( \textcolor{black}{\rho} | \hbar), 1 \rangle_{\hbar}  \\ &=& P(V_1, V_2, \ldots ) \langle \Upsilon_v( \textcolor{black}{\rho} | \hbar), 1 \rangle_{\hbar} .\end{eqnarray} \textcolor{black}{Note that} $\langle \cdot, \cdot \rangle_{\hbar}$ conjugates the left entry.  Since $\langle \Upsilon_v(\textcolor{black}{\rho} | \hbar) , 1\rangle_{\hbar} =1$, \textcolor{black}{(3) follows:} \begin{equation} \label{Lizard} \textcolor{black}{\langle \Upsilon_v ( \rho | \hbar) , P(\rho_1, \rho_2, \ldots )  \rangle_{\hbar}  = P(V_1, V_2, \ldots ).} \ \ \square \end{equation}

\noindent \textcolor{black}{By (\ref{CoherentState}) and (\ref{Lizard}), for any two $v^{\textnormal{out}}, v^{\textnormal{in}}$ satisfying $\sum_{k=1}^{\infty} \frac{ |V_k^{\textnormal{out}} |^2 }{k} < \infty$ and $\sum_{k=1}^{\infty} \frac{ |V_k^{\textnormal{in}} |^2 }{k} < \infty$,  \begin{equation} \label{NewOutInNormalizationFactor} \langle \Upsilon_{v^{\textnormal{in} } } (\rho | \hbar ) , \Upsilon_{v^{\textnormal{out}}}( \rho | \hbar ) \rangle_{\hbar} = \Upsilon_{v^{\textnormal{out}} } ( v^{\textnormal{in}} | \hbar) = \textnormal{exp} \Bigg ( \frac{1}{\hbar} \sum_{k=1}^{\infty} \frac{\overline{V_k^{\textnormal{out}} } V_k^{\textnormal{in}} }{k} \Bigg ) \end{equation}\noindent is the normalization factor in the Definition [\ref{DefinitionJackMeasures}] of Jack measures $M(v^{\textnormal{out}}, v^{\textnormal{in}})$.  That being said,} \textcolor{black}{Jack polynomials play no role in the proof of Lemma [\ref{CoherentStateLemma}]}.   In fact, using Definition [\ref{DefinitionJacks}] of Jack polynomials as eigenfunctions of a self-adjoint operator in $\mathcal{F}$, the Cauchy identity (\ref{CauchyIdentity}) follows from Lemma [\ref{CoherentStateLemma}]\textcolor{black}{:} for any $v$ \textcolor{black}{satisfying $\sum_{k=1}^{\infty} \frac{ |V_k |^2 }{k} < \infty$} and any \textcolor{black}{orthonormal} basis $B_{\lambda}$ of $\mathcal{F}$, (\ref{Lizard}) implies \textcolor{black}{$\Upsilon_v (\rho |\hbar) = \sum_{\lambda} \overline{B_{\lambda}(V_1, V_2, \ldots)} B_{\lambda}(\rho_1, \rho_2, \ldots)$} \textcolor{black}{and so (\ref{NewOutInNormalizationFactor}) is also} \begin{equation} \sum_{\lambda} { \overline{B_{\lambda} (V_1^{\textnormal{out}}, V_2^{\textnormal{out}}, \ldots )} B_{\lambda} (V_1^{\textnormal{in}}, V_2^{\textnormal{in}}, \ldots )} = \textcolor{black}{\textnormal{exp} \Bigg ( \frac{1}{\hbar} \sum_{k=1}^{\infty} \frac{\overline{V_k^{\textnormal{out}}}  V_k^{\textnormal{in} } }{k} \Bigg )} .\end{equation}

\subsubsection{\textbf{Step 2: Joint moments from commuting operators}} \label{SUBSUBSECStep2} In the previous Step 1, we saw that the Jack measure $M(v^{\textnormal{out}}, v^{\textnormal{in}})$ itself is defined in terms of the decomposition of the reproducing kernel $\Upsilon_{v^{\textnormal{out}}} ( v^{\textnormal{in}} | \hbar)$ over the basis of Jack polynomials.  We now use this relationship to show that any operator acting diagonally on Jack polynomials provides a method of moments for Jack measures.

\begin{lemma} \label{BornRuleLemma} Suppose $\widehat{O}$ is a self-adjoint operator in $\mathcal{F}$ which acts diagonally on Jack polynomials
\begin{equation}\label{EigenvalueRelation} \widehat{O} P_{\lambda} (\rho_1, \rho_2, \ldots | \ebar, \hbar) = O_{\lambda} P_{\lambda} (\rho_1, \rho_2, \ldots | \ebar, \hbar) \end{equation}

\noindent with eigenvalue $O_{\lambda} \in \mathbb{R}$.  Then for any two specializations $v^{\textnormal{out}}, v^{\textnormal{in}}$ of $\rho_1, \rho_2, \ldots$ which define a Jack measure $M(v^{\textnormal{out}}, v^{\textnormal{in}})$ as in \textnormal{Definition [\ref{DefinitionJackMeasures}]}, the expected value of the random variable $O_{\lambda}$ is \begin{equation} \ \mathbb{E} [ O_{\lambda} ] = \frac{\big \langle \textcolor{black}{\Upsilon}_{v^{\textnormal{in}}} ( \textcolor{black}{\rho} | \hbar), \widehat{O} \ \textcolor{black}{\Upsilon}_{v^{\textnormal{out}}} ( \textcolor{black}{\rho} | \hbar) \big \rangle_{\hbar} }{ \textcolor{black}{\big \langle \textcolor{black}{\Upsilon}_{v^{\textnormal{in}}} ( \textcolor{black}{\rho} | \hbar), \textcolor{black}{\Upsilon}_{v^{\textnormal{out}}} ( \textcolor{black}{\rho} | \hbar) \big \rangle_{\hbar}}} \end{equation} the matrix element of $\widehat{O}$ between $\textcolor{black}{\Upsilon}_{v^{\textnormal{out}}} (\textcolor{black}{\rho} | \hbar)$ and $\textcolor{black}{\Upsilon}_{v^{\textnormal{in}}}(\textcolor{black}{\rho} | \hbar)$ from \textcolor{black}{\textnormal{(\ref{CoherentState})} normalized by \textnormal{(\ref{NewOutInNormalizationFactor})}}.\end{lemma}

\noindent \textit{Proof:} The eigenvalues $O_{\lambda}$ are real-valued functions $O: \mathbb{Y} \rightarrow \mathbb{R}$ on the discrete sample space $\mathbb{Y}$ of all partitions $\lambda$ and therefore define a random variable associated to the Jack measure $M(v^{\textnormal{out}}, v^{\textnormal{in}})$.  Using (\ref{JackMeasureLaw}), \textcolor{black}{formula (\ref{Lizard})} \textcolor{black}{and its complex conjugate}, \textcolor{black}{(\ref{NewOutInNormalizationFactor})}, (\ref{EigenvalueRelation}), \textcolor{black}{$\overline{O_{\lambda}} = O_{\lambda}$,} and $\widehat{O}^{\dagger} = \widehat{O}$, \begin{eqnarray} \mathbb{E} [ O_{\lambda} ] &=& \sum_{\lambda} O_{\lambda} \textnormal{Prob}(\lambda) \nonumber \\ 
&=& \sum_{\lambda} O_{\lambda}  P_{\lambda}^{\textnormal{norm}} ( V^{\textnormal{in}}_1, V_2^{\textnormal{in}}, \ldots | \ebar, \hbar) \overline{P^{\textnormal{norm}}_{\lambda}( V^{\textnormal{out}}_1, V_2^{\textnormal{out}},... | \ebar, \hbar)}  \cdot \textnormal{exp} \Big ( -  \frac{1}{\hbar} \sum_{k=1}^{\infty} \frac{\overline{V_k^{\textnormal{out}}} V_k^{\textnormal{in}}}{k} \Big ) \nonumber \\ 
&=& \sum_{\lambda}  O_{\lambda}  \langle \Upsilon_{v^{\textnormal{in}}} ( \textcolor{black}{\rho} | \hbar), P_{\lambda}^{\textnormal{norm}} ( \textcolor{black}{\rho} | \ebar, \hbar) \rangle_{\hbar}  { \langle P_{\lambda}^{\textnormal{norm}} ( \textcolor{black}{\rho} | \ebar, \hbar) ,  \Upsilon_{v^{\textnormal{out}}}( \textcolor{black}{\rho} | \hbar)\rangle_{\hbar}}\cdot \textcolor{black}{\big \langle \textcolor{black}{\Upsilon}_{v^{\textnormal{in}}} ( \textcolor{black}{\rho} | \hbar), \textcolor{black}{\Upsilon}_{v^{\textnormal{out}}} ( \textcolor{black}{\rho} | \hbar) \big \rangle_{\hbar}^{-1}} \nonumber \\ &=& \sum_{\lambda}  \langle \Upsilon_{v^{\textnormal{in}}} ( \textcolor{black}{\rho} | \hbar), P_{\lambda}^{\textnormal{norm}} ( \textcolor{black}{\rho} | \ebar, \hbar) \rangle_{\hbar}  { \langle \widehat{O} P_{\lambda}^{\textnormal{norm}} ( \textcolor{black}{\rho} | \ebar, \hbar) ,  \Upsilon_{v^{\textnormal{out}}}( \textcolor{black}{\rho} | \hbar)\rangle_{\hbar}}\cdot \textcolor{black}{\big \langle \textcolor{black}{\Upsilon}_{v^{\textnormal{in}}} ( \textcolor{black}{\rho} | \hbar), \textcolor{black}{\Upsilon}_{v^{\textnormal{out}}} ( \textcolor{black}{\rho} | \hbar) \big \rangle_{\hbar}^{-1}} \nonumber \\ &=& \sum_{\lambda}   \langle \Upsilon_{v^{\textnormal{in}}} ( \textcolor{black}{\rho} | \hbar), P_{\lambda}^{\textnormal{norm}} ( \textcolor{black}{\rho} | \ebar, \hbar) \rangle_{\hbar}  { \langle P_{\lambda}^{\textnormal{norm}} ( \textcolor{black}{\rho} | \ebar, \hbar) ,  \widehat{O} \Upsilon_{v^{\textnormal{out}}}( \textcolor{black}{\rho} | \hbar)\rangle_{\hbar}}\cdot \textcolor{black}{\big \langle \textcolor{black}{\Upsilon}_{v^{\textnormal{in}}} ( \textcolor{black}{\rho} | \hbar), \textcolor{black}{\Upsilon}_{v^{\textnormal{out}}} ( \textcolor{black}{\rho} | \hbar) \big \rangle_{\hbar}^{-1}} \nonumber \\ 
 &=& \big \langle \textcolor{black}{\Upsilon}_{v^{\textnormal{in}}} ( \textcolor{black}{\rho} | \hbar), \widehat{O} \ \textcolor{black}{\Upsilon}_{v^{\textnormal{out}}} ( \textcolor{black}{\rho} | \hbar) \big \rangle_{\hbar} \cdot \textcolor{black}{\big \langle \textcolor{black}{\Upsilon}_{v^{\textnormal{in}}} ( \textcolor{black}{\rho} | \hbar), \textcolor{black}{\Upsilon}_{v^{\textnormal{out}}} ( \textcolor{black}{\rho} | \hbar) \big \rangle_{\hbar}^{-1}}\nonumber \end{eqnarray} where in the last equality we recognize $\sum_{\lambda}  | P_{\lambda}^{\textnormal{norm}} \rangle \langle P_{\lambda}^{\textnormal{norm}} |$ as the resolution of the identity. $\square$

\noindent Lemma [\ref{BornRuleLemma}] for the expected value of a single random variable $O_{\lambda}$ can be equivalently stated as a result for the joint moments of several random variables $Q_1 |_{\lambda}, Q_2 |_{\lambda}, \ldots, Q_{\ell}|_{\lambda}, \ldots$ as follows.
\begin{lemma} \label{MultipleBornRuleLemma} Suppose $\widehat{Q}_{\ell}$ is a family of self-adjoint operators in $\mathcal{F}$ indexed by $\ell=1,2,3,\ldots$ which commute and are simultaneously diagonalized on Jack polynomials with eigenvalues $Q_{\ell} |_{\lambda}$.  Then for any two specializations $v^{\textnormal{out}}, v^{\textnormal{in}}$ of the variables $\rho_1, \rho_2, \ldots$ which define a Jack measure $M(v^{\textnormal{out}}, v^{\textnormal{in}})$ as in \textnormal{Definition [\ref{DefinitionJackMeasures}]} and any $\ell_1, \ldots, \ell_n$, the joint moments of the random variables $Q_{\ell} |_{\lambda}$ may be computed as matrix elements \begin{equation} \ \mathbb{E} [ Q_{\ell_1} |_{\lambda} \cdots Q_{\ell_n} |_{\lambda} ] = \frac{\big \langle \textcolor{black}{\Upsilon}_{v^{\textnormal{in}}} ( \textcolor{black}{\rho} | \hbar), \widehat{Q}_{\ell_1} \cdots \widehat{Q}_{\ell_n} \cdot \textcolor{black}{\Upsilon}_{v^{\textnormal{out}}} ( \textcolor{black}{\rho} | \hbar) \big \rangle_{\hbar}}{ \textcolor{black}{\big \langle \textcolor{black}{\Upsilon}_{v^{\textnormal{in}}} ( \textcolor{black}{\rho} | \hbar), \textcolor{black}{\Upsilon}_{v^{\textnormal{out}}} ( \textcolor{black}{\rho} | \hbar) \big \rangle_{\hbar}}} .\end{equation} \end{lemma}

\noindent \textit{Proof:} Follows automatically from Lemma [\ref{BornRuleLemma}] after choosing $\widehat{O} = \widehat{Q}_{\ell_1} \cdots \widehat{Q}_{\ell_n}$. $\square$

\subsubsection{\textbf{Step 3: Nazarov-Sklyanin hierarchy of commuting operators}} \label{SUBSUBSECStep3} In the previous Step 2, we saw that joint moments of random variables $Q_{\ell} |_{\lambda}$ associated to Jack measures $M(v^{\textnormal{out}}, v^{\textnormal{in}})$ can be realized as matrix elements involving the reproducing kernel\textcolor{black}{s $\Upsilon_{v^{\textnormal{out}}}, \Upsilon_{v^{\textnormal{in}}}$} if one can find commuting self-adjoint operators $\widehat{Q}_{\ell}$ which are diagonalized on Jack polynomials with eigenvalues $Q_{\ell}|_{\lambda}$.  To apply this reasoning to prove Theorem [\ref{Theorem3AOE}], we need commuting operators $\widehat{T}_{\ell}^{\uparrow}(\ebar, \hbar)$ which satisfy\begin{equation} \label{DreamRelation} \widehat{T}^{\uparrow}_{\ell}( \ebar, \hbar) P_{\lambda}( \rho_1, \rho_2, \ldots | \ebar, \hbar) = T^{\uparrow}_{\ell}(\ebar, \hbar) |_{\lambda} P_{\lambda}( \rho_1, \rho_2, \ldots | \ebar, \hbar) ,\end{equation} i.e. which are diagonalized on Jack polynomials with eigenvalues $T_{\ell}^{\uparrow}(\ebar, \hbar)|_{\lambda}$ from (\ref{TransitionStatistics}).  We now recall a remarkable explicit construction of $\widehat{T}^{\uparrow}_{\ell}(\ebar, \hbar)$ satisfying (\ref{DreamRelation}) from Nazarov-Sklyanin \cite{NaSk2}.

\begin{theorem} \label{NSTheorem} \textnormal{[Nazarov-Sklyanin \cite{NaSk2}]} \textcolor{black}{For any $\ebar \in \mathbb{R}$ and $\hbar >0$,} consider the infinite matrix \begin{equation} \label{QuantumLaxMatrix} \widehat{L}_{\bullet}(\ebar, \hbar) = \begin{bmatrix} 
0 & \widehat{\rho}_{-1} &  \widehat{\rho}_{-2} &  \widehat{\rho}_{-3} & \cdots  \\ 
 \widehat{\rho}_{1} & \textcolor{black}{\ebar} &  \widehat{\rho}_{-1} &  \widehat{\rho}_{-2} & \ddots  \\
  \widehat{\rho}_{2} &  \widehat{\rho}_{1} & 2 \textcolor{black}{\ebar}& \widehat{\rho}_{-1} & \ddots \\
   \widehat{\rho}_{3} &  \widehat{\rho}_2 &  \widehat{\rho}_1 & 3 \textcolor{black}{\ebar} & \ddots \\
  \vdots & \vdots & \ddots & \ddots & \ddots \end{bmatrix} \end{equation}  where $\widehat{\rho}_{\pm k}$ are the creation and annihilation operators \textnormal{(\ref{CreationOperator}), (\ref{AnnihilationOperator})}.  For $j \in \mathbb{Z}_{\bullet} = \{0,1,2,3,\ldots\}$, let $\psi_j = [ 0 \ 0 \ \cdots 0 \ 1 \ 0 \cdots ]^T$ denote the column vector which is $1$ in the $j$th entry and $0$ otherwise.  For any $\ell=0,1,2,3,\ldots$, let $\widehat{T}_{\ell}^{\uparrow}(\ebar, \hbar)$ denote the top-left entry of the $\ell$th power of \textnormal{(\ref{QuantumLaxMatrix})}, namely \begin{equation} \label{NSHierarchy} \widehat{T}^{\uparrow}_{\ell} (\ebar, \hbar) = \langle \psi_0 , \widehat{L}_{\bullet}(\ebar, \hbar)^{\ell} \psi_0 \rangle \end{equation} if $\langle \cdot, \cdot \rangle$ is defined on the span of $\psi_j$ by declaring $\psi_j$ orthonormal.  The following holds for \textnormal{(\ref{NSHierarchy})}:
  \begin{enumerate}
  \item For any $\ell$, $\widehat{T}^{\uparrow}_{\ell} (\ebar, \hbar)$ are unbounded self-adjoint operators on $\mathcal{F}$ all defined on $\C[\rho_1, \rho_2, \ldots]$.
  \item For any $\ell_1 \neq \ell_2$, $\widehat{T}^{\uparrow}_{\ell} (\ebar, \hbar)$ commute $[ \widehat{T}_{\ell_1}^{\uparrow} (\ebar, \hbar) , \widehat{T}_{\ell_2}^{\uparrow} (\ebar, \hbar) ] = 0$ if restricted to $\C[\rho_1, \rho_2, \ldots]$.
  \item For any $\ell$, $\widehat{T}_{\ell}(\ebar, \hbar)$ satisfy \textnormal{(\ref{DreamRelation})}, acting diagonally on Jack polynomials $P_{\lambda}( \rho_1, \rho_2, \ldots | \ebar, \hbar)$ with eigenvalues $T_{\ell}^{\uparrow}(\ebar, \hbar)|_{\lambda}$ defined from the anisotropic partition profile $f_{\lambda}(c  | \ebar, \hbar)$ by \textnormal{(\ref{TransitionStatistics})}.
  \end{enumerate}
\end{theorem}

\noindent In \textsection 8 of \cite{Moll2}, we verified that Theorem [\ref{NSTheorem}] agrees with Theorem 2 in \cite{NaSk2}.  \textcolor{black}{Note that $L_{\bullet}$ in (6.2) in \cite{Moll2} should have the form (\ref{QuantumLaxMatrix}) since $\alpha - 1 = \ebar / (-\varepsilon_2)$ is on the diagonal of $L$ in (6.1) in \cite{NaSk2}.} We refer to the set of commuting operators (\ref{NSHierarchy}) as the \textit{Nazarov-Sklyanin hierarchy}.  Let $\delta$ denote the Kronecker delta and $\widehat{\rho}_0 = 0$.  For $j, j' =0,1,2,\ldots$, the $(j,j')$th entry of $\widehat{L}_{\bullet}(\ebar, \hbar)$ in (\ref{QuantumLaxMatrix}) is \begin{equation} \label{MatrixElementGoodFormula} \langle \psi_{j} , \widehat{L}_{\bullet}(\ebar, \hbar) \psi_{j'} \rangle = \widehat{\rho}_{j' - j} +  \textcolor{black}{\ebar} j \delta(j-j') .\end{equation} \noindent The term $\widehat{\rho}_{j'-j}$ is a creation operator if $j<j'$, is $0$ if $j=j'$, or an annihilation operator if $j>j'$.  The term $  \textcolor{black}{\ebar} j \delta(j-j')$ only appears if $j=j'$, the diagonal of (\ref{QuantumLaxMatrix}), and vanishes if $\ebar=0$.  Define {\small \begin{equation} \nonumber \widehat{S}(j_1, \ldots, j_{\ell-1} | \ebar, \hbar) = \widehat{\rho}_{j_1} \cdot \Big (  \widehat{\rho}_{j_2 - j_1} +  \textcolor{black}{\ebar}  j_1 \delta(j_1 - j_2)  \Big ) \cdots \Big ( \widehat{\rho}_{j_{\ell-1} - j_{\ell-2}}  +   \textcolor{black}{\ebar}  j_{\ell-2} \delta(j_{\ell-2} - j_{\ell-1} ) \Big ) \cdot \widehat{\rho}_{-j_{\ell-1}}. \end{equation}} \noindent By (\ref{MatrixElementGoodFormula}), the $\ell$th operator $\widehat{T}^{\uparrow}_{\ell} (\ebar, \hbar)$ in the Nazarov-Sklyanin hierarchy (\ref{NSHierarchy}) is therefore \begin{equation} \label{ExactFormula}  \widehat{T}^{\uparrow}_{\ell}(\ebar, \hbar) = \sum_{j_1, \ldots, j_{\ell-1} =0}^{\infty} \widehat{S}(j_1, \ldots, j_{\ell-1} | \ebar, \hbar) \end{equation} \noindent When $\ell=3$, formula (\ref{ExactFormula}) for $T_3^{\uparrow} (\ebar, \hbar)$ coincides with the formula (\ref{OperatorIntro}) for $\mathcal{H}(\ebar, \hbar)$ in \textsection [\ref{SUBSECJackPolynomials}].

\subsubsection{\textbf{Step 4: Weighted enumeration of ribbon paths from Nazarov-Sklyanin path operators}} \label{SUBSUBSECStep4} In this step we relate the $(v^{\textnormal{out}}, v^{\textnormal{in}})$-weighted enumeration of ribbon paths introduced in \textsection [\ref{SUBSECWeightedEnumerationRibbonPaths}] to the Nazarov-Sklyanin hierarchy introduced in Step 3.  Recall by Definition [\ref{DefinitionRibbonPaths}] that a ribbon path $\boldsymbol{\vec{\gamma}}=(\gamma; \textit{\textbf{p}}_1, \ldots, \textit{\textbf{p}}_q)$ on $n=1$ site of length $\ell$ is $1$ sliding path $\gamma$ with some number $q$ of pairings.

\begin{lemma} \label{PathPathPathLemma} Fix $(j_1, \ldots, j_{\ell-1}) \in \mathbb{Z}_{\bullet}^{\ell-1}$.  Let $\Delta_{1,q,m}(\ell)[j_1, \ldots, j_{\ell-1}]$ denote the finite set of ribbon paths $\boldsymbol{\vec{\gamma}}=(\gamma; \textit{\textbf{p}}_1, \ldots, \textit{\textbf{p}}_q)$ on $1$ site whose $\gamma$ has vertices $\{(0,0), (1, j_1), \ldots, (\ell-1, j_{\ell-1}), (\ell,0)\}$ with $m$ slides and which has $q$ pairings.  Consider the Nazarov-Sklyanin path operator {\small \begin{equation} \nonumber \widehat{S}(j_1, \ldots, j_{\ell-1} | \ebar, \hbar) = \widehat{\rho}_{j_1} \cdot \Big (  \widehat{\rho}_{j_2 - j_1} +  \textcolor{black}{\ebar}  j_1 \delta(j_1 - j_2)  \Big ) \cdots \Big ( \widehat{\rho}_{j_{\ell-1} - j_{\ell-2}}  +  \textcolor{black}{\ebar}  j_{\ell-2} \delta(j_{\ell-2} - j_{\ell-1} ) \Big ) \cdot \widehat{\rho}_{-j_{\ell-1}} \end{equation} } \noindent from the expansion \textnormal{(\ref{ExactFormula})} of the Nazarov-Sklyanin hierarchy \textnormal{(\ref{NSHierarchy})}.  For any two specializations $v^{\textnormal{out}}, v^{\textnormal{in}}$, the reproducing kernel matrix elements of the path operator are polynomials in $\hbar$ and $\ebar$ \begin{equation} \label{DesiredIdentity20402040} \frac{\langle \textcolor{black}{\Upsilon}_{v^{\textnormal{in}}}, \widehat{S}(j_1, \ldots, j_{\ell-1} | \ebar, \hbar)  \cdot \textcolor{black}{\Upsilon}_{v^{\textnormal{out}}} \rangle_{\hbar}}{ \textcolor{black}{\langle \Upsilon_{v^{\textnormal{in}}} , \Upsilon_{v^{\textnormal{out}} } \rangle_{\hbar}}} = \textcolor{black}{\ebar^m}\sum_{q=0}^{\ell /2} \sum_{\boldsymbol{\vec{\gamma}} \in \Delta_{1,q,m}(\ell)[j_1, \ldots, j_{\ell-1}] } \mathfrak{W}(\boldsymbol{\vec{\gamma}} |v^{\textnormal{out}},v^{\textnormal{in}}) \hbar^q  \end{equation} \noindent whose coefficients are $(v^{\textnormal{out}}, v^{\textnormal{in}})$-weighted sums of ribbon paths with $\mathfrak{W}$ from \textnormal{Definition [\ref{DefinitionRibbonPathsWeight}]}. \end{lemma}

\noindent \textit{Proof of Lemma} [\ref{PathPathPathLemma}]: First, to each path operator $\widehat{S}(j_1, \ldots, j_{\ell-1} | \ebar, \hbar)$ we associate a sliding path $\gamma$ of length $\ell$ in $\mathbb{Z}_{\bullet}^2$ as in Definition [\ref{DefinitionSlidingPaths}] with vertices $\{(0,0), (1, j_1), \ldots, (\ell-1, j_{\ell-1}), (\ell,0)\}$.  The locations of the slides are predetermined by $j_1, \ldots, j_{\ell-1}$: a slide at height $j$ occurs whenever $j_{a} = j_{a+1}$.  Conversely, $\widehat{S}(j_1, \ldots, j_{\ell-1}| \ebar, \hbar)$ is determined by a sliding path $\gamma$ with such vertices: 
\begin{itemize}
\item To each jump of degree $+k$, associate \textcolor{black}{the} creation operator $\widehat{\rho}_{+k}$ as in (\ref{CreationOperator})
\item To each jump of degree $-k$, associate \textcolor{black}{the} annihilation operator $\widehat{\rho}_{-k}$ as in (\ref{AnnihilationOperator})
\item To each slide at height $j$, associate the scalar multiple $ \textcolor{black}{\ebar}  j$ of the identity operator in $\mathcal{F}$.
\end{itemize}

\begin{figure}[htb]
\centering
\includegraphics[width=0.7 \textwidth]{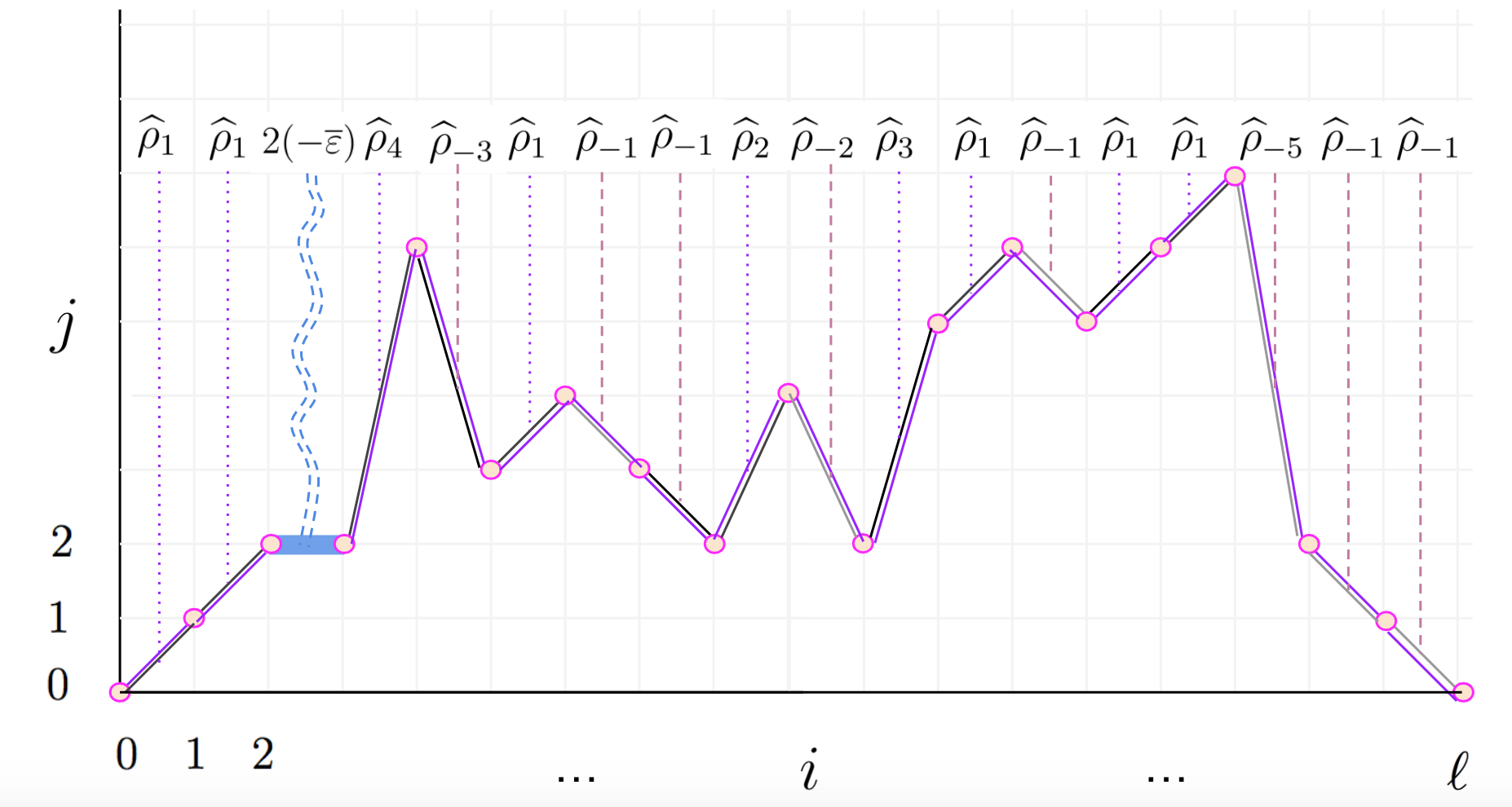}
 \vspace*{-5mm}
\caption{The path operator $\widehat{S}(1,2,\textcolor{black}{2},6,3,4,3,2,4,2,5,6,5,6,7,2,1)$ is the ordered product $\widehat{\rho}_{1} \widehat{\rho}_1 (2  \textcolor{black}{\ebar} ) \widehat{\rho}_4 \widehat{\rho}_{-3} \widehat{\rho}_1 \widehat{\rho}_{-1} \widehat{\rho}_{-1} \widehat{\rho}_2 \widehat{\rho}_{-2} \widehat{\rho}_3 \widehat{\rho}_1 \widehat{\rho}_{-1} \widehat{\rho}_1 \widehat{\rho}_1 \widehat{\rho}_{-5} \widehat{\rho}_{-1} \widehat{\rho}_{-1}$ parametrized by a sliding path $\gamma$ of length $\ell=18$ with $m=1$ slide at height $j=2$.}
\label{SoperatorFIG}
\end{figure}

\noindent By these rules, the path operator $\widehat{S}(j_1, \ldots, j_{\ell-1} | \ebar, \hbar)$ is recovered by taking the sliding path $\gamma$ and multiplying the associated operators in order from left to right.  For an illustration, see Figure [\ref{SoperatorFIG}].

\pagebreak

\noindent We have just exhibited a bijection between path operators and sliding paths $\gamma$.  To reflect this, write $\widehat{S}(\gamma | \ebar, \hbar) = \widehat{S}(j_1, \ldots, j_{\ell-1} | \ebar, \hbar)$ and $ \Delta_{1,q,m}(\ell)[\gamma] = \Delta_{1,q,m}(\ell)[j_1, \ldots, j_{\ell-1}]$.  Next, we show that each path operator itself has an explicit polynomial expansion in $\ebar$ and $\hbar$ indexed by ribbon paths: {\small  \begin{equation} \label{MegaDecomposition} { \widehat{S}(\gamma | \ebar, \hbar) =  \textcolor{black}{\ebar}^m  \prod_{\textit{\textbf{e}} \in \textnormal{\textbf{E}}_S({{\gamma}})}   \textnormal{height}(\textit{\textbf{e}})  \sum_{\boldsymbol{\vec{\gamma}} \in \Delta_{1,q,m}(\ell)[\gamma]} \hbar^q \prod_{ \textbf{p} \in \boldsymbol{\vec{\gamma}}}  \textnormal{size}(\textbf{p}) \prod_{\textbf{\textit{e}} \in \mathbf{E}_{UJ}^+[\boldsymbol{\vec{\gamma}}]} \widehat{\rho}_{\deg ( \textbf{\textit{e}})} \prod_{\textbf{\textit{e}} \in \mathbf{E}_{UJ}^-[\boldsymbol{\vec{\gamma}}]} \widehat{\rho}_{\deg ( \textbf{\textit{e}})} }\end{equation}}

\noindent To derive (\ref{MegaDecomposition}), consider the definition of path operators $\widehat{S}$ in the statement of Lemma [\ref{PathPathPathLemma}]:
\begin{itemize}
\item $\widehat{S}(j_1, \ldots, j_{\ell-1}| \ebar, \hbar)$ may have $j_a = j_{a+1}$.  In this case, the term $\widehat{\rho}_{j_{a+1} - j_a} +  \textcolor{black}{\ebar} j_a \delta(j_a - j_{a+1})$ in the path operator $\widehat{S}(\gamma | \ebar, \hbar)$ becomes $\widehat{\rho}_{0} +  \textcolor{black}{\ebar} j_a =  \textcolor{black}{\ebar} j_a$ since we have defined $\widehat{\rho}_0=0$.  On the other hand, in the sliding path $\gamma$, such $j_a = j_{a+1}$ corresponds to a slide $\textbf{\textit{e}} \in \mathbf{E}_S(\gamma)$ at height $j_a$, so the contribution of slides to the path operator is $ \textcolor{black}{\ebar} \textnormal{height}(\textbf{\textit{e}})$ for each of the $m$ slides, in agreement with the multiplicative factor on the right-hand side of (\ref{MegaDecomposition}).
\item $\widehat{S}(j_1, \ldots, j_{\ell-1} | \ebar, \hbar)$ is not normally-ordered: annihilation operators $\widehat{\rho}_{-k}$ can be to the left of creation operators $\widehat{\rho}_{+k'}$ in $\widehat{S}$.  In this case, replace $\widehat{\rho}_{-k} \widehat{\rho}_{+k'} =\widehat{\rho}_{+k'} \widehat{\rho}_{-k} + \hbar k \delta(k-k')$ using \begin{equation} \label{LeibnizLove} [\widehat{\rho}_{-k}, \widehat{\rho}_{+k'}] =  \hbar k \delta(k-k') . \end{equation} On the other hand, in the case $k=k'$, such $\widehat{\rho}_{-k}$ and $\widehat{\rho}_{+k}$ in the $\widehat{S}(\gamma | \ebar, \hbar)$ labeled by $\gamma$ are indexed by two edges $\textbf{\textit{e}}$ and $\textbf{\textit{e}'}$ of degrees $-k$ and $k$ in $\gamma$ which satisfy Definition [\ref{DefinitionPairings}] of a pairing $\textit{\textbf{p}}$ of size $k$.  By repeated applications of (\ref{LeibnizLove}), bring all $\widehat{\rho}_{-k}$ to the right of all $\widehat{\rho}_{k'}$ at the cost of an additive terms whenever $k=k'$ and a contribution $\hbar k$.  This results in (\ref{MegaDecomposition}) indexed by ribbon paths $\boldsymbol{\vec{\gamma}} = (\gamma | \textit{\textbf{p}}_1, \ldots, \textit{\textbf{p}}_q)$ where the number of pairings $q$ is the number of applications of (\ref{LeibnizLove}), each pairing has weight $\textnormal{size}(\textbf{p})$, and the normally-ordered product of creation and annihilation operators on the right-hand side of (\ref{MegaDecomposition}) are indexed by unpaired jumps $\mathbf{E}_{UJ}^{\pm}[\boldsymbol{\vec{\gamma}}]$ of positive and negative degrees, respectively.
\end{itemize}
\noindent Finally, the normally-ordered terms in (\ref{MegaDecomposition}) have reproducing kernel matrix elements
\begin{equation} \label{DoubleReproduction} {\small  {\langle \textcolor{black}{\Upsilon}_{v^{\textnormal{in}}}, \prod_{\textbf{\textit{e}} \in \mathbf{E}_{UJ}^+[\boldsymbol{\vec{\gamma}}]} \widehat{\rho}_{\deg ( \textbf{\textit{e}})} \prod_{\textbf{\textit{e}} \in \mathbf{E}_{UJ}^-[\boldsymbol{\vec{\gamma}}]} \widehat{\rho}_{-\deg ( \textbf{\textit{e}})} \cdot \textcolor{black}{\Upsilon}_{v^{\textnormal{out}}} \rangle_{\hbar}}= {\textcolor{black}{\langle \Upsilon_{v^{\textnormal{in}}}, \Upsilon_{v^{\textnormal{out}}} \rangle_{\hbar} } } \prod_{\textbf{\textit{e}} \in \mathbf{E}_{UJ}^-[\boldsymbol{\vec{\gamma}}]} \overline{V^{\textnormal{out}}_{-\deg ( \textbf{\textit{e}})}} \prod_{\textbf{\textit{e}} \in \mathbf{E}_{UJ}^+[\boldsymbol{\vec{\gamma}}]} {V}_{\deg ( \textbf{\textit{e}})}^{\textnormal{in}} } .\end{equation}

\noindent The identity (\ref{DoubleReproduction}) follows by applying $\widehat{\rho}_{-k} \textcolor{black}{\Upsilon}_{v} = \overline{V_k} \textcolor{black}{\Upsilon}_{v} $ \textcolor{black}{from (2)} in Lemma [\ref{CoherentStateLemma}] twice using $\widehat{\rho}_{\pm k}^{\dagger} = \widehat{\rho}_{\mp k}$.  Substituting (\ref{DoubleReproduction}) into (\ref{MegaDecomposition}) yields (\ref{DesiredIdentity20402040}). $\square$\\
\\
 \noindent The identical argument above implies a generalization of Lemma [\ref{PathPathPathLemma}] to the case of $n$ sites.
\begin{lemma} \label{MultiplePathPathPathLemma} Fix $\{\ell_a\}_{a=1}^n \subset \mathbb{Z}_{+}$ and $(j^{(a)}_1, \ldots, j^{(a)}_{\ell_a-1}) \in \mathbb{Z}_{\bullet}^{\ell_a-1}$.  Let $\Delta_{n,q,m}(\ell)[\{j^{(a)}_1, \ldots, j^{(a)}_{\ell-1}\}_{a=1}^n]$ denote the finite set of ribbon paths $\boldsymbol{\vec{\gamma}}=(\gamma_1, \ldots, \gamma_n; \textit{\textbf{p}}_1, \ldots, \textit{\textbf{p}}_q)$ on $n$ sites whose $\gamma_a$ has vertices $\{(0,0), (1, j_1^{(a)}), \ldots, (\ell_a-1, j^{(a)}_{\ell_a-1}), (\ell_a,0)\}$ with $m_a$ slides so that $m = \sum_{a=1}^m m_a$. The reproducing kernel matrix element of the ordered product of $n$ Nazarov-Sklyanin path operators is proportional to $\textcolor{black}{\ebar}^m$ and is a polynomial in $\hbar$ whose coefficients are $(v^{\textnormal{out}}, v^{\textnormal{in}})$-weighted sums of ribbon paths  \begin{equation} \label{DesiredIdentity20402040MULTIPLE} \frac{\langle \textcolor{black}{\Upsilon}_{v^{\textnormal{in}}} , \prod_{a=1}^n \widehat{S}(j^{(a)}_1, \ldots, j^{(a)}_{\ell_a-1} | \ebar, \hbar)  \cdot \textcolor{black}{\Upsilon}_{v^{\textnormal{out}}} \rangle_{\hbar} } { \textcolor{black}{\langle \Upsilon_{v^{\textnormal{in}}}, \Upsilon_{v^{\textnormal{out}} } \rangle_{\hbar}}} = \sum_{q=0}^{\ell /2} \sum_{\boldsymbol{\vec{\gamma}} \in \Delta_{n,q,m}(\ell)[\{j^{(a)}_1, \ldots, j^{(a)}_{\ell-1}\}_{a=1}^n] } \mathfrak{W}(\boldsymbol{\vec{\gamma}} |v^{\textnormal{out}}, v^{\textnormal{in}}) \hbar^q  \textcolor{black}{\ebar}^m \nonumber \end{equation} \noindent in $\Delta_{1,q,m}(\ell)[\{j^{(a)}_1, \ldots, j^{(a)}_{\ell-1}\}_{a=1}^n]$ where $\mathfrak{W}(\gamma|v^{\textnormal{out}}, v^{\textnormal{in}})$ is the weight in \textnormal{Definition [\ref{DefinitionRibbonPathsWeight}]}. \end{lemma} 

\subsubsection{\textbf{Step 5: Estimates for weighted enumeration of ribbon paths}} \label{SUBSUBSECStep5}  In the previous Step 4, we fixed $\ell_1, \ldots, \ell_n$ and $n$ sequences $(j_1^{(a)}, \ldots, j^{(a)}_{\ell_a-1}) \in \mathbb{Z}_{\bullet}^{\ell_a-1}$.  We now fix $\ell_1, \ldots, \ell_n$ but instead sum over all $(j_1^{(a)}, \ldots, j^{(a)}_{\ell_a-1}) \in \mathbb{Z}_{\bullet}^{\ell_a-1}$ and show that the resulting $\textcolor{black}{(v^{\textnormal{out}}, v^{\textnormal{in}})}$-weighted sum over ribbon paths in (\ref{JointMomentsAOE}) is finite.  To do so, we use the \textcolor{black}{regularity} assumption on $v^{\textnormal{out}}, v^{\textnormal{in}}$ \textcolor{black}{in Theorem [\ref{Theorem3AOE}]}.

\begin{lemma} \label{EstimateLemma} Fix $\ebar \in \mathbb{R}$, $\hbar>0$, and $v^{\textnormal{out}} = \{V^{\textnormal{out}}_k\}_{k=1}^{\infty}$, $v^{\textnormal{in}} = \{V^{\textnormal{in}}_k\}_{k=1}^{\infty}$ two specializations satisfying $\max (|V_k^{\textnormal{out}}|, |V_k^{\textnormal{in}}|) \leq A r^k$ for some $A>0$ and $0<r<1$.  For $\vec{\ell}= (\ell_1, \ldots, \ell_n) \in \mathbb{Z}_+^n$, let $|| \vec{\ell}||_1 = \ell_1 + \cdots + \ell_n$.  Let $\Delta_{n,q,m}(\ell_1, \ldots, \ell_n)$ be the infinite set of ribbon paths $\boldsymbol{\vec{\gamma}}$ on $n$ sites of lengths $\ell_1, \ldots, \ell_n$ with $q$ pairings and $m$ slides.   Consider the weight $\mathfrak{W}(\boldsymbol{\vec{\gamma}} | v^{\textnormal{out}}, v^{\textnormal{in}} )$ on ribbon paths in \textnormal{Definition [\ref{DefinitionRibbonPathsWeight}]}.  Then the following infinite series converges
\begin{equation} \label{OMGFINITE} \sum_{q=0}^{\tfrac{1}{2} || \vec{\ell} ||_1} \sum_{m=0}^{|| \vec{\ell}||_1} \sum_{\boldsymbol{\vec{\gamma}} \in \Delta_{n,q,m}(\ell_1, \ldots, \ell_n )} \mathfrak{W}(\boldsymbol{\vec{\gamma}} | v^{\textnormal{out}}, v^{\textnormal{in}}) \ \hbar^q \textcolor{black}{\ebar}^m < \infty. \end{equation}\end{lemma}

\noindent \textit{Proof of Lemma} [\ref{EstimateLemma}]: A ribbon path on $n$ sites $\boldsymbol{\vec{\gamma}} = (\gamma_1, \ldots, \gamma_n ; \textit{\textbf{p}}_1, \ldots, \textit{\textbf{p}}_q)$ is the data of $n$ sliding paths $\gamma_a$ together with some number $q$ of pairings.  Let $\{(i, j_i^{(a)})\}_{i=0}^{\ell_a} \subset \mathbb{Z}_{\bullet}^2$ denote the $\ell_a$ vertices of $\gamma_a$. Let $\max \gamma_a = \max_{i} j_i^{(a)}$ be the maximum height of $\gamma_a$.  For $M \in \mathbb{Z}_{\bullet}$, define \begin{equation} \Delta_{n,q,m}^{M} (\ell_1, \ldots, \ell_n)  = \{ \boldsymbol{\vec{\gamma}} \in \Delta_{n,q,m}(\ell_1, \ldots, \ell_n ) \ : \ \max_{a=1, \ldots, n} \max \gamma_a = M\} .\end{equation} We can now rewrite the sum in (\ref{OMGFINITE}) by keeping track of the maximum heights: \textcolor{black}{we need to show}
\begin{equation}\sum_{q=0}^{\tfrac{1}{2} || \vec{\ell} ||_1} \sum_{m=0}^{|| \vec{\ell}||_1} \sum_{M= 0}^{\infty} \sum_{\boldsymbol{\vec{\gamma}} \in \Delta^{M}_{n,q,m}(\ell_1, \ldots, \ell_n ) } \mathfrak{W}(\boldsymbol{\vec{\gamma}} | v^{\textnormal{out}}, v^{\textnormal{in}}) \ \hbar^q \textcolor{black}{\ebar}^m \textcolor{black}{<\infty.} \end{equation}

\noindent For \textcolor{black}{any fixed} $\boldsymbol{\vec{\gamma}} \in \Delta^{M}_{n,q,m}(\ell_1, \ldots, \ell_n )$, using the decomposition (\ref{RibbonPathStepSetDecomposition}), we derive three bounds:
\begin{enumerate}
\item $\boldsymbol{\vec{\gamma}}$ has a total of $m$ slides.  By Definition [\ref{DefinitionRibbonPathsWeight}], slides in $\mathbf{E}_S[\boldsymbol{\vec{\gamma}}]$ contribute weight $j \leq M$.
\item $\boldsymbol{\vec{\gamma}}$ has a total of $q$ pairings.   By Definition [\ref{DefinitionRibbonPathsWeight}], pairings $\textit{\textbf{p}} \in \boldsymbol{\vec{\gamma}}$ contribute weight $k \leq M$.
\item $\boldsymbol{\vec{\gamma}}$ has a total of $n$ sliding paths which overall achieve maximum height $M$.  In order to achieve these maximum heights, there must exist some number $U \textcolor{black}{ \leq || \vec{\ell}||_1}$ of {unpaired up jumps} of positive degrees $k^+_1, \ldots, k^+_U \geq 1$ and some number $D \textcolor{black}{ \leq || \vec{\ell}||_1}$ of {unpaired down jumps} of negative degrees $-k^-_1, \ldots, -k^-_D$ so $k^{+}_1 + \cdots + k^+_U \geq M$ and $k^{-}_1 + \cdots + k^-_D  \geq M$.  By Definition [\ref{DefinitionRibbonPathsWeight}], unpaired jumps $\textbf{\textit{e}} \in \mathbf{E}_{UJ}^{\pm} [\boldsymbol{\vec{\gamma}}]$ contribute weights $\overline{V_k^{\textnormal{out}}}$ and $V_k^{\textnormal{in}}$, respectively, hence overall these unpaired jumps together contribute weight bounded by 
\begin{equation} \Big | \prod_{i=1}^U \overline{V_{k^{+}_i}^{\textnormal{out}}} \prod_{j=1}^{D} {V}^{\textnormal{in}}_{k^-_j} \Big | \leq \textcolor{black}{A^{U+D}} \prod_{i=1}^U r^{k_i^+} \prod_{j=1}^D r^{k_j^-} \leq \textcolor{black}{\max (A^{2 || \vec{\ell}||_1} ,1 )} \ r^{2M} .\end{equation} Here we made use of our regularity assumption on the specializations $v^{\textnormal{out}}, v^{\textnormal{in}}$.
\end{enumerate}

\noindent By (1), (2), (3), \textcolor{black}{since $\max ( A^{2 || \vec{\ell} ||_1}, 1)$ is independent of $q$, $m$, and $M$} for fixed $\vec{\ell}$, it is enough to show 
\begin{equation} \label{WhatIsLeft}  \sum_{q=0}^{\tfrac{1}{2} || \vec{\ell} ||_1} \sum_{m=0}^{|| \vec{\ell}||_1} \sum_{M= 0}^{\infty}  \ \ r^{2M} M^{q+m} \hbar^q \textcolor{black}{\ebar}^m\sum_{\boldsymbol{\vec{\gamma}} \in \Delta^{M}_{n,q,m}(\ell_1, \ldots, \ell_n ) }  1<\infty. \end{equation}
\noindent To finish, estimate the number of ribbon paths in $\Delta^{M}_{n,q,m}(\ell_1, \ldots, \ell_n )$:
\begin{itemize}
\item There are $(M+1)^{|| \vec{\ell}||_1}$ choices of sliding paths $\gamma_1, \ldots, \gamma_n$ since the $j$ coordinate of vertices $(i, j_i^{(a)})$ in $\gamma_a$ take values in $0, 1, \ldots, M$. 
\item There are $\binom{|| \vec{\ell}||_1}{2}$ choices of pairs of edges in $\gamma_1, \ldots, \gamma_n$ which may take place in a pairings.  
\end{itemize} \noindent It now suffices to show for fixed $\ell_1, \ldots, \ell_n$ that \begin{equation} \label{WhatIsReallyLeft} \sum_{q=0}^{\tfrac{1}{2} || \vec{\ell} ||_1} \sum_{m=0}^{|| \vec{\ell}||_1}  \hbar^q \textcolor{black}{\ebar}^m \ \  \Bigg ( \sum_{M= 0}^{\infty} r^{2M} M^{m+q} (M+1)^{ || \vec{\ell}||_1}  \binom{|| \vec{\ell}||_1}{2}^q \Bigg )    <\infty. \end{equation}
\noindent For each $q,m$, the series in parentheses converges due to the $r^{2M}$ term, completing the proof. $\square$

\subsubsection{\textbf{Step 6: Joint moments are weighted enumerations of ribbon paths}} \label{SUBSUBSECStep6} We now prove (\ref{JointMomentsAOE}).\\
\\
\noindent \textit{Proof of Part I of Theorem \textnormal{[\ref{Theorem3AOE}]}}: Consider the Nazarov-Sklyanin hierarchy (\ref{NSHierarchy}) of self-adjoint operators $\widehat{{T}}^{\uparrow}_{\ell}(\ebar, \hbar)$.    By part (3) of Theorem [\ref{NSTheorem}] in Step 3, $\widehat{{T}}^{\uparrow}_{\ell}(\ebar, \hbar)$ act diagonally on Jack polynomials with eigenvalue $T_{\ell}^{\uparrow}(\ebar, \hbar)|_{\lambda}$ defined in (\ref{TransitionStatistics}).  Applying Lemma [\ref{MultipleBornRuleLemma}] from Step 2 in the case $\widehat{Q}_{\ell} = \widehat{T}^{\uparrow}_{\ell}(\ebar, \hbar)$, we can realize the joint moments of the random variables $T_{\ell}^{\uparrow}(\ebar, \hbar)|_{\lambda}$ sampled from the Jack measure $M(v^{\textnormal{out}}, v^{\textnormal{in}})$ as the matrix elements
\begin{equation} \label{WhereDeyAtDo} \mathbb{E} \big [ T_{\ell_1}^{\uparrow} (\ebar,\hbar) \cdots T_{\ell_n}^{\uparrow} (\ebar, \hbar) \big ]  = \frac{\langle \textcolor{black}{\Upsilon}_{v^{\textnormal{in}}} , \widehat{T}^{\uparrow}_{\ell_1} (\ebar, \hbar) \cdots \widehat{T}^{\uparrow}_{\ell_n} (\ebar, \hbar) \cdot \textcolor{black}{\Upsilon}_{v^{\textnormal{out}}} \rangle_{\hbar} }{\textcolor{black}{\langle \Upsilon_{v^{\textnormal{in}}}, \Upsilon_{v^{\textnormal{out}}} \rangle_{\hbar}}}.  \end{equation} 
\noindent \textcolor{black}{Using} (\ref{ExactFormula}) \textcolor{black}{to write} each $\widehat{T}^{\uparrow}_{\ell_a}(\ebar, \hbar)$ as an infinite sum of path operators, the right side of (\ref{WhereDeyAtDo}) is 
\begin{equation} \label{CiddleCorn0} \sum_{j^{(1)}_1, \ldots, j^{(1)}_{\ell_1-1} =0}^{\infty} \cdots \sum_{j^{(n)}_1, \ldots, j^{(n)}_{\ell_n-1} =0}^{\infty} \frac{\langle \textcolor{black}{\Upsilon}_{v^{\textnormal{in}}} , \textcolor{black}{\prod_{a=1}^n}\widehat{S}(j^{(a)}_1, \ldots, j^{(a)}_{\ell_{a-1}-1} | \ebar, \hbar) \cdot \textcolor{black}{\Upsilon}_{v^{\textnormal{out}}} \rangle_{\hbar}}{\textcolor{black}{\langle \Upsilon_{v^{\textnormal{in}}}, \Upsilon_{v^{\textnormal{out}}} \rangle_{\hbar}}} .\end{equation} \noindent For each fixed $\{j_1^{(a)}, \ldots, j_{\ell_a-1}^{(a)}\}_{a=1}^n$, apply Lemma [\ref{MultiplePathPathPathLemma}] from Step 4 to expand each of these matrix element of path operators over ribbon paths: (\ref{CiddleCorn0}) becomes
\begin{equation} \label{CiddleCorn} \sum_{j^{(1)}_1, \ldots, j^{(1)}_{\ell_1-1} =0}^{\infty} \cdots \sum_{j^{(n)}_1, \ldots, j^{(n)}_{\ell_n-1} =0}^{\infty} \sum_{q=0}^{\tfrac{1}{2} || \vec{\ell}||_1} \sum_{m=0}^{|| \vec{\ell}||_1} \sum_{\boldsymbol{\vec{\gamma}} \in \Delta_{n,q,m} (\ell_1, \ldots, \ell_n)[\{j_1^{(a)}, \ldots, j_{\ell_a -1}^{(a)}\}_{a=1}^n] } \mathfrak{W}(\boldsymbol{\vec{\gamma}} | v^{\textnormal{out}}, v^{\textnormal{in}})  \hbar^q \textcolor{black}{\ebar}^m \nonumber \end{equation}
\noindent where $\Delta_{n,q,m}(\ell_1, \ldots, \ell_n) [ \{j_1^{(a)}, \ldots, j_{\ell_a-1}^{(a)}\}_{a=1}^n]$ is the finite set of ribbon paths with fixed underlying $n$ sliding paths $\gamma_1, \ldots, \gamma_n$.  Since we sum over all possible $j_1^{(a)}, \ldots, j_{\ell_a - 1}^{(a)}$, we have proven \begin{equation}   \mathbb{E} \big [ T_{\ell_1}^{\uparrow} (\ebar,\hbar) \cdots T_{\ell_n}^{\uparrow} (\ebar, \hbar) \big ]  = \sum_{q=0}^{\tfrac{1}{2} || \vec{\ell} ||_1} \sum_{m=0}^{|| \vec{\ell}||_1} \sum_{\boldsymbol{\vec{\gamma}} \in \Delta_{n,q,m}(\ell_1, \ldots, \ell_n )} \mathfrak{W}(\boldsymbol{\vec{\gamma}} | v^{\textnormal{out}}, v^{\textnormal{in}}) \ \hbar^q \textcolor{black}{\ebar}^m \end{equation} \noindent which is (\ref{JointMomentsAOE}) by the definition of $Y_{n,q,m}(\ell_1, \ldots, \ell_n |v^{\textnormal{out}}, v^{\textnormal{in}})$ in (\ref{RibbonPathYFunction}).  Finally, by Lemma [\ref{EstimateLemma}] in Step 5, these joint moments are finite and the proof of Part I of Theorem [\ref{Theorem3AOE}] is complete. $\square$

\subsubsection{\textbf{Step 7: Joint cumulants are weighted enumerations of connected ribbon paths}} \label{SUBSUBSECStep7} In the previous step, we showed that joint moments (\ref{JointMomentsAOE}) for Jack measures $M(v^{\textnormal{out}}, v^{\textnormal{in}})$ are polynomials in $\hbar$ and $\textcolor{black}{\ebar}$ whose coefficients are a $(v^{\textnormal{out}}, v^{\textnormal{in}})$-weighted count of ribbon paths.  In this final step we prove (\ref{JointCumulantsAOE}) relating joint cumulants and connected ribbon paths.  The following argument is standard in combinatorics and featured in Chapter 5 of Stanley \cite{StanleyVol2}.  For further background on the relationship between joint cumulants and the enumeration of connected structures, see Novak \cite{NovakFreeLectures}.  \\
\\
\noindent \textit{Proof of Part II of Theorem \textnormal{[\ref{Theorem3AOE}]}}: First, we decompose an arbitrary ribbon path $\boldsymbol{\vec{\gamma}}$ as a disjoint union of connected ribbon paths. Let $[n]=\{1,2,\ldots, n\}$ and let $\mathsf{P}(n)$ denote the set of set partitions $\pi$ of $[n]$.  A set partition $\pi$ of $[n]$ is a collection of disjoint subsets $B \subset [n]$ so that $\cup_{B \in \pi} B = [n]$.  For every ribbon path $\boldsymbol{\vec{\gamma}}$ \textcolor{black}{on $n$ sites}, let $\pi( \boldsymbol{\vec{\gamma}})$ be the set partition of $[n]$ defined by $B \in \pi(\boldsymbol{\vec{\gamma}})$ if and only if $B$ is a connected component of the reduced graph of $\boldsymbol{\vec{\gamma}}$ in Definition [\ref{DefinitionReducedGraph}].  For any fixed ribbon path $\boldsymbol{\vec{\gamma}} = (\gamma_1, \ldots, \gamma_n; \textit{\textbf{p}}_1, \ldots, \textit{\textbf{p}}_q)$ and any fixed $B \in \pi(\boldsymbol{\vec{\gamma}})$ , define $\boldsymbol{\vec{\gamma}}^{(B)}$ to be the ribbon path with sliding paths $\{\gamma_a\}_{a \in B}$ and those pairings $\textit{\textbf{p}}$ among $\textit{\textbf{p}}_1, \ldots, \textit{\textbf{p}}_q$ which pair jumps in $\gamma_a$ and $\gamma_{a'}$ for $a,a' \in B$.  By Definition [\ref{DefinitionConnectedRibbonPaths}], $\boldsymbol{\vec{\gamma}}^{(B)}$ is a connected ribbon path supported on the $|B|$ sites $B \subset [n]$.  \textcolor{black}{As a consequence, any ribbon path $\boldsymbol{\vec{\gamma}}$ decomposes into connected components} \begin{equation} \label{Decomposition} \boldsymbol{\vec{\gamma}} = \coprod_{B \in \pi(\boldsymbol{\vec{\gamma}})}  \boldsymbol{\vec{\gamma}}^{(B)}. \end{equation}
\noindent Next, observe that the weight $\mathfrak{W}(\boldsymbol{\vec{\gamma}} | v^{\textnormal{out}}, v^{\textnormal{in}})$ in Definition [\ref{DefinitionRibbonPathsWeight}] is multiplicative, hence \begin{equation} \label{WeightedDecomposition} \mathfrak{W}(\boldsymbol{\vec{\gamma}}| v^{\textnormal{out}}, v^{\textnormal{in}}) = \prod_{B \in \pi(\boldsymbol{\vec{\gamma}})}  \mathfrak{W}(\boldsymbol{\vec{\gamma}}^{(B)}  |v^{\textnormal{out}}, v^{\textnormal{in}}) .\end{equation} \noindent Substituting (\ref{WeightedDecomposition}) into the $(v^{\textnormal{out}}, v^{\textnormal{in}})$-weighted sum (\ref{JointMomentsAOE}) gives \begin{equation} \label{HereWeAreGraduation}  \mathbb{E} \big [ T_{\ell_1}^{\uparrow} (\ebar,\hbar) \cdots T_{\ell_n}^{\uparrow} (\ebar, \hbar) \big ]  = \sum_{q=0}^{\tfrac{1}{2} || \vec{\ell} ||_1} \sum_{m=0}^{|| \vec{\ell}||_1} \sum_{\boldsymbol{\vec{\gamma}} \in \Delta_{n,q,m}(\ell_1, \ldots, \ell_n )} \prod_{B \in \pi(\boldsymbol{\vec{\gamma}})}  \mathfrak{W}(\boldsymbol{\vec{\gamma}}^{(B)}  |v^{\textnormal{out}}, v^{\textnormal{in}}) \ \hbar^q \textcolor{black}{\ebar}^m . \end{equation}
\noindent We now relabel the sums in (\ref{HereWeAreGraduation}) as follows.  Let $\boldsymbol{\vec{\gamma}}$ be a ribbon path in $\Delta_{n,q,m}(\ell_1, \ldots, \ell_n )$.
\begin{itemize}
\item $\boldsymbol{\vec{\gamma}}$ is a ribbon path on $n$ sites of lengths $\ell_1, \ldots, \ell_n$.  For each $B \in \pi(\boldsymbol{\vec{\gamma}})$, the connected component $\boldsymbol{\vec{\gamma}}^{(B)}$ is a connected ribbon path on $|B|$ sites of lengths $\vec{\ell}^{(B)} =\{\ell_a \  : \ a \in B\}$. 
\item $\boldsymbol{\vec{\gamma}}$ has $m$ slides.  For each $B \in \pi(\boldsymbol{\vec{\gamma}})$, there is some $m_B \geq 0$ so that $\boldsymbol{\vec{\gamma}}^{(B)}$ has $m_B$ slides.
\item $\boldsymbol{\vec{\gamma}}$ has $q$ pairings.  For each $B \in \pi(\boldsymbol{\vec{\gamma}})$, there is some $q_B \geq 0$ so that $\boldsymbol{\vec{\gamma}}^{(B)}$ has $q_B$ pairings.  However, by assumption this $\boldsymbol{\vec{\gamma}}^{(B)}$ is itself a connected ribbon path, and since it is on $|B|$ sites, $q_B$ cannot be arbitrary but must satisfy $q_B \geq |B|-1$.  This is a fundamental fact: to connect $|B|$ vertices it requires at least $|B|-1$ edges.  Define $g_B = q_B - |B|+1$.
\end{itemize}
\noindent With these \textcolor{black}{relabelings, we can pass from $\Delta$ to $\Gamma$ in Problem [\ref{ProblemRibbonPathEnumeration}]:} the right side of (\ref{HereWeAreGraduation}) is
 \begin{equation} \sum_{\pi \in \mathsf{P}(n)} \prod_{B \in \pi} \Bigg ( \sum_{g_B =0}^{\tfrac{1}{2} || \vec{\ell}^{(B)} ||_1 \textcolor{black}{+ 1 -|B|} } \sum_{m_B=0}^{|| \vec{\ell}^{(B)} ||_1} \sum_{\boldsymbol{\vec{\gamma}}^{(B)} \in \textcolor{black}{\Gamma_{|B|, g_B, m_B}} (\vec{\ell}^{(B)}) } \mathfrak{W}(\boldsymbol{\vec{\gamma}}^{(B)} | v^{\textnormal{out}}, v^{\textnormal{in}}) \ \hbar^{|B|-1+g_B} \textcolor{black}{\ebar}^{m_B} \Bigg )  \end{equation}
  \noindent so by definition of the sum $W$ over connected ribbon paths in (\ref{RibbonPathWFunction}), for ${T}^{\uparrow}_{\ell} = {T}^{\uparrow}_{\ell}(\ebar, \hbar)$ we obtain
 { \begin{equation} \label{WhenIsThisOver} \mathbb{E} \big [ T_{\ell_1}^{\uparrow} \cdots T_{\ell_n}^{\uparrow} \big ]  = \sum_{\pi \in \mathsf{P}(n)} \prod_{B \in \pi} \Bigg ( \sum_{g_B =0}^{\tfrac{1}{2} || \vec{\ell}^{(B)} ||_1 } \sum_{m_B=0}^{|| \vec{\ell}^{(B)} ||_1} W_{|B|, g_B, m_B} (\vec{\ell}^{(B)} | v^{\textnormal{out}}, v^{\textnormal{in}}) \ \hbar^{|B|-1 + g_B} \textcolor{black}{\ebar}^{m_B} \Bigg ).  \end{equation}} \noindent We can now prove the desired result (\ref{JointCumulantsAOE}) by induction in $n$.  First, the $n=1$ case of (\ref{JointCumulantsAOE}) \begin{equation} \kappa_1(T^{\uparrow}_{\ell} (\ebar, \hbar)) = \sum_{g=0}^{\ell/2} \sum_{m=0}^{\ell} W_{1,g,m}(\ell | v^{\textnormal{out}}, v^{\textnormal{in}}) \hbar^{g} \textcolor{black}{\ebar}^m \end{equation} follows immediately from the $n=1$ case of (\ref{JointMomentsAOE}) for two simple reasons:
 \begin{itemize}
 \item $\kappa_1(T) = \mathbb{E}[T]$ the first cumulant is the expected value of any random variable $T$, and
 \item $W_{1,g,m}(\ell |v^{\textnormal{out}}, v^{\textnormal{in}}) = Y_{1,g,m}(\ell | v^{\textnormal{out}}, v^{\textnormal{in}})$ since all ribbon paths on $1$ site are connected.
 \end{itemize}
 
 \noindent Next, assume by induction that (\ref{JointCumulantsAOE}) holds for any number of sites strictly less than $n$.  In particular, given any $B \subset [n]$ with $|B| \leq n-1$, let 
$\kappa_{|B|} (\{T^{\uparrow}_{\ell_a}(\ebar, \hbar)\}_{a \in B} )$ denote the $|B|$th cumulant of the random variables $T^{\uparrow}_{\ell_a}(\ebar, \hbar)$ indexed by $a \in B \subset [n]$.  The inductive hypothesis then reads \begin{equation} \label{GoodInductiveHypothesis} \kappa_{|B|} (\{T^{\uparrow}_{\ell_a}(\ebar, \hbar)\}_{a \in B} ) = \sum_{g_B=0}^{\frac{1}{2} ||\vec{\ell}^{(B)} ||_1 } \sum_{m_B=0}^{ ||\vec{\ell}^{(B)} ||_1} W_{|B|, g_B, m_B} (\vec{\ell}^{(B)} | v^{\textnormal{out}}, v^{\textnormal{in}} ) \ \hbar^{|B|-1+g_B} \textcolor{black}{\ebar}^{m_B} \end{equation} if we use the same notation $\vec{\ell}^{(B)}= \{\ell_a\}_{a \in B}$ as above.  We can now prove (\ref{JointCumulantsAOE}) for $\kappa_n$.  Recall the recursive formula which expresses $\kappa_n$ in terms of those $\kappa_{n'}$ with $n' < n$: for $T_{\ell}^{\uparrow} = T_{\ell}^{\uparrow}(\ebar, \hbar)$, \begin{equation} \label{DestinationX} \kappa_n\big ( T_{\ell_1}^{\uparrow}, \ldots, \widehat{T}^{\uparrow}_{\ell_n}  \big ) = \mathbb{E} \big [T_{\ell_1}^{\uparrow}\cdots \widehat{T}^{\uparrow}_{\ell_n} \big ] - \sum_{\pi \in \mathsf{P}^*(n)} \prod_{B \in \pi} \kappa_{|B|}  (\{T^{\uparrow}_{\ell_a}\}_{a \in B} ). \end{equation} \noindent The crucial feature of (\ref{DestinationX}) is that the sum is not over all set partitions but only over \begin{equation} \mathsf{P}^*(n) := \{ \pi \in \mathsf{P}(n) \ : \ \pi \neq \{[n]\} \} \end{equation} which excludes $\pi =\{ [n] \}$ with only one $B = [n]$.  Substituting formula (\ref{WhenIsThisOver}) for the joint moment and the inductive hypothesis (\ref{GoodInductiveHypothesis}) into (\ref{DestinationX}) yields (\ref{JointCumulantsAOE}) since (\ref{DestinationX}) is then a difference of a sums over $\mathsf{P}(n)$ and a sum over $\mathsf{P}^*(n)$.  This completes the proof of Theorem [\ref{Theorem3AOE}]. $\square$

\subsection{Example: ribbon paths on $1$ site of length $2$ and expected area of $\Omega_{\lambda}(\ebar, \hbar)$} \label{SUBSECProofExpectedSizeProposition} 
$ \ $\\
$ \ $\\
\noindent \textit{Proof 1 of Proposition} [\ref{PropositionSizeExpectation}]: The operator $\widehat{T}^{\uparrow}_2 = \sum_{j_1=0}^{\infty} \widehat{\rho}_{j_1} \widehat{\rho}_{-j_1}$ known as the degree operator acts diagonally on Jack polynomials with eigenvalue $\hbar | \lambda|$ since $\widehat{\rho}_{-j_1} = \hbar j_1 \frac{\partial}{\partial \rho_{j_1}}$.  By Lemma [\ref{BornRuleLemma}], $\mathbb{E}[| \lambda| ] = \frac{1}{\hbar} { \langle \textcolor{black}{\Upsilon_v} , \widehat{T}^{\uparrow}_2  \textcolor{black}{\Upsilon_v} \rangle_{\hbar} }{ \langle \textcolor{black}{\Upsilon_v} ,\textcolor{black}{\Upsilon_v} \rangle_{\hbar}^{-1}}$.  Using $\widehat{\rho}_{-j_1} \textcolor{black}{\Upsilon_v }= \overline{V}_{j_1} \textcolor{black}{\Upsilon_v}$ which is \textcolor{black}{(2)} in Lemma [\ref{CoherentStateLemma}] and $\widehat{\rho_{j_1}}^{\dagger} = \widehat{\rho}_{-j_1}$, one computes $\langle \textcolor{black}{\Upsilon_v}, \widehat{T}^{\uparrow}_2  \cdot \textcolor{black}{\Upsilon_v} \rangle_{\hbar} = \textcolor{black}{\langle \Upsilon_v, \Upsilon_v \rangle } \sum_{j_1=0}^{\infty} |V_{j_1}|^2 $ which completes the proof. $\square$\\
\\
\noindent \textit{Proof 2 of Proposition} [\ref{PropositionSizeExpectation}]: By (\ref{TransitionStatistics}), $T_2^{\uparrow} |_{\lambda}= \hbar | \lambda|$.  The case $n=1$, $\ell=2$ of formula (\ref{JointMomentsAOE}) in Theorem [\ref{Theorem3AOE}] with $v^{\textnormal{out}} = v^{\textnormal{in}} = v$ implies $\mathbb{E}[T_2^{\uparrow}|_{\lambda} ] = \sum_{q=0}^{1} \sum_{m=0}^2 Y_{1,q,m}(2|v,v) \hbar^q \textcolor{black}{\ebar}^m$.  But $Y_{1,q,m}(2|v,v) = W_{1,q,m}(2|v,v) = 0$ if $(q,m) \neq (0,0)$ as can be directly checked, while $Y_{1,0,0}(2|v,v) =W_{1,0,0}(2|v,v)= \sum_{j_1=0}^{\infty} |V_{j_1}|^2$ comes from Szeg\H{o} paths $\gamma = \{(0,0), (1, j_1), (2,0)\}$ of length $2$.  Note $\widehat{S}(j_1 | \ebar, \hbar) = \widehat{\rho}_{j_1} \widehat{\rho}_{-j_1}$ in Proof 1 is the path operator associated to this $\gamma$. $\square$





\section{Limit shapes of Jack measures on profiles} \label{SECProof1} 

\noindent In \textsection [\ref{SUBSECLLNProof}] we prove Theorem [\ref{Theorem1LLN}], our first main result for Jack measures: as $\hbar \rightarrow 0$, the random profile $f_{\lambda}(c | \ebar, \hbar) \rightarrow \mathbf{f}( c | v; \ebar)$ forms a limit shape which is the convex action profile $\mathbf{f}( c| v;0)$ in \textsection [\ref{SECSzegoConvex}] or dispersive action profile $\mathbf{f}( c | v;\ebar)$ in \textsection [\ref{SECSlidingDispersive}] depending on whether or not $\ebar \rightarrow 0$ as $\hbar \rightarrow 0$.  In \textsection [\ref{SUBSECJackPlancherelLimitShapes}], we illustrate our Theorem [\ref{Theorem1LLN}] for Poissonized Jack-Plancherel measures $M(v_{\textnormal{PL}}, v_{\textnormal{PL}})$.  In \textsection [\ref{APPENDIXdePoissonization}] we show that the same limit shapes form for the Jack-Plancherel measures $M_d(v_{\textnormal{PL}}, v_{\textnormal{PL}})$ as $d \rightarrow \infty$.
\subsection{Proof of first main result} \label{SUBSECLLNProof} We now prove Theorem [\ref{Theorem1LLN}] in Steps 1-6 below.

\subsubsection{\textbf{Step 1: Linear statistics $O_p (\ebar, \hbar)|_{\lambda}$ are polynomials in $T_{\ell}^{\uparrow} (\ebar, \hbar)|_{\lambda}$}} For $p=1,2,\ldots$, define \begin{equation} \label{TheOh} O_p( \ebar, \hbar) |_{\lambda} := \int_{-\infty}^{+\infty} c^p \cdot \tfrac{1}{2} f_{\lambda}''(c | \ebar, \hbar) dc. \end{equation} \noindent For $\lambda$ sampled from a Jack measure, $O_p(\ebar, \hbar)|_{\lambda}$ is a random variable we call the  \textit{$p$th linear statistic} of $f_{\lambda}(c | \ebar, \hbar)$.  Our Theorem [\ref{Theorem1LLN}] under consideration concerns the asymptotic behavior of these linear statistics.  To apply Theorem [\ref{Theorem3AOE}], we need to relate $O_p(\ebar, \hbar)|_{\lambda}$ to $T_{\ell}^{\uparrow}(\ebar, \hbar)|_{\lambda}$.  These $T_{\ell}^{\uparrow}(\ebar, \hbar)$ are defined by (\ref{TransitionStatistics}), and we can rewrite (\ref{TransitionStatistics}) using (\ref{TheOh}) as \begin{equation} \label{KMKRelationOT} \sum_{\ell=0}^{\infty} T_{\ell}^{\uparrow}(\ebar, \hbar)|_{\lambda} u^{- \ell-1} = \frac{1}{u}  \cdot \textnormal{exp} \Bigg ( \sum_{p=1}^{\infty} O_p (\ebar, \hbar) \frac{ u^{-p}}{p} \Bigg ) .\end{equation} 

\noindent As a consequence of (\ref{KMKRelationOT}), we have a polynomial relation:

\begin{lemma} \label{KMKPolynomialsLemma} \textnormal{[\textsection 3.3 in Kerov [\cite{Ke1}]} For $p \in \mathbb{Z}_+$, there are polynomials $\textnormal{KMK}_p$ independent of $\lambda$ so \begin{equation} \label{GreatPolynomialRelation} O_p(\ebar, \hbar)|_{\lambda} = \textnormal{KMK}_p \Big (T_1^{\uparrow} (\ebar, \hbar) |_{\lambda}, \ldots, T_p^{\uparrow}(\ebar, \hbar)|_{\lambda} \Big ).  \end{equation} Moreover, these $\textnormal{KMK}_p$ have no constant term. \end{lemma}

\noindent $\textnormal{KMK}_p$ originate in Kerov's Markov-Kre\u{\i}n Correspondence in \cite{Ke1}.  \textcolor{black}{Formula} (3.3.7) in \cite{Ke1} is:
\begin{eqnarray} O_1 &=& T_1^{\uparrow} \\ O_2 &=& 2 T_2^{\uparrow} - (T_1^{\uparrow})^2 \\ O_3 &=& 3 T_3^{\uparrow} - 3 T_1^{\uparrow} T_2^{\uparrow} + (T_1^{\uparrow})^3 \\ O_4 &=& 4 T_4^{\uparrow} - 4 T_1^{\uparrow} T_3^{\uparrow} - 2 (T_2^{\uparrow})^2 + 4 (T_1^{\uparrow})^2 T_2^{\uparrow} - (T_1^{\uparrow})^4. \end{eqnarray}

\subsubsection{\textbf{Step 2: \textcolor{black}{Convergence in probability} of linear statistics to constants as $\hbar \rightarrow 0$}} We now show that as $\hbar \rightarrow 0$, the random variables  $O_p(\ebar, \hbar) |_{\lambda}$ in (\ref{TheOh}) sampled from $M(\textcolor{black}{v^{\textnormal{out}},v^{\textnormal{in}}})$ \textcolor{black}{with arbitrary $v^{\textnormal{out}},v^{\textnormal{in}}$ satisfying the assumptions of Theorem [\ref{Theorem3AOE}] converge in probability} \begin{equation} O_p(\ebar, \hbar) |_{\lambda} \rightarrow \mathbf{O}_p(\textcolor{black}{v^{\textnormal{out}},v^{\textnormal{in}}}; \ebar) \end{equation} \noindent to deterministic constants $\mathbf{O}_p(\textcolor{black}{v^{\textnormal{out}},v^{\textnormal{in}}}; \ebar)$ \textcolor{black}{in Regimes I, II, and III}.  \textcolor{black}{Since $\mathbf{O}_p (v^{\textnormal{out}}, v^{\textnormal{in}}; \ebar)$ is a constant}, it is enough to show \textcolor{black}{the weak law of large numbers} \begin{eqnarray} \label{ExpectationLLNWeWant} \mathbb{E} [ O_p(\ebar, \hbar) |_{\lambda} ] & \rightarrow & \mathbf{O}_p(\textcolor{black}{v^{\textnormal{out}},v^{\textnormal{in}}}; \ebar) \\ \label{VarianceLLNWeWant} \textnormal{Var} [ O_p(\ebar, \hbar) |_{\lambda} ] & \rightarrow & 0 \end{eqnarray} 

\noindent where $\textnormal{Var}[O] = \mathbb{E}[O^2] - \mathbb{E}[O]^2$.  To establish (\ref{ExpectationLLNWeWant}) and (\ref{VarianceLLNWeWant}), we first prove asymptotic factorization of the joint moments $\mathbb{E} \big [ T_{\ell_1}^{\uparrow}(\ebar, \hbar) \cdots T_{\ell_n}^{\uparrow}(\ebar, \hbar) \big ]$ of the $T_{\ell}^{\uparrow}(\ebar, \hbar)$ using Theorem [\ref{Theorem3AOE}] from \textsection [\ref{SECRibbonAnisotropic}]. \begin{lemma} \label{AsymptoticFactorizationLemmaGood} For any $n \in \mathbb{Z}_+$ and $\textcolor{black}{\vec{\ell}} = (\ell_1, \ldots, \ell_n)\in \mathbb{Z}_{\geq 0}^n$, there are constants $\mathbf{T}_{\ell}^{\uparrow}(\textcolor{black}{v^{\textnormal{out}},v^{\textnormal{in}}}; \ebar)$ so as $\hbar \rightarrow 0$ \textcolor{black}{with $\ebar$ fixed}, \begin{equation} \label{AsymptoticFactorizationLemmaFormula} \mathbb{E} \big [ T_{\ell_1}^{\uparrow}(\ebar, \hbar) \cdots T_{\ell_n}^{\uparrow}(\ebar, \hbar) \big ] \rightarrow \mathbf{T}_{\ell_1}^{\uparrow}(\textcolor{black}{v^{\textnormal{out}},v^{\textnormal{in}}}; \ebar) \cdots \mathbf{T}_{\ell_n}^{\uparrow}(\textcolor{black}{v^{\textnormal{out}},v^{\textnormal{in}}}; \ebar). \end{equation} \noindent \textcolor{black}{Also, if $\ebar , \hbar \rightarrow 0$ so $\ebar \sim \hbar^{1/2}$ as in \textnormal{Regime II}, \textnormal{(\ref{AsymptoticFactorizationLemmaFormula})} holds with $\mathbf{T}^{\uparrow}_{\ell}(\textcolor{black}{v^{\textnormal{out}},v^{\textnormal{in}}}; 0) := \lim\limits_{\ebar \rightarrow 0} \mathbf{T}_{\ell}^{\uparrow}(\textcolor{black}{v^{\textnormal{out}},v^{\textnormal{in}}}; \ebar)$.}  \end{lemma} \begin{itemize}
\item \textit{Proof of} Lemma [\ref{AsymptoticFactorizationLemmaGood}].  By Theorem [\ref{Theorem3AOE}], the $\hbar \rightarrow 0$ limit is \textcolor{black}{finite and given by} \begin{equation}  \label{TransparentThing} \mathbb{E} \big [ T_{\ell_1}^{\uparrow}(\ebar, \hbar) \cdots T_{\ell_n}^{\uparrow}(\ebar, \hbar) \big ]  \rightarrow \sum_{m=0}^{|| \vec{\ell}||_1} Y_{n,0,m}(\ell_1, \ldots, \ell_n | \textcolor{black}{v^{\textnormal{out}},v^{\textnormal{in}}}) \textcolor{black}{\ebar}^m \end{equation} \noindent \textcolor{black}{if $\ebar$ is fixed.} By the definition of $Y$ in (\ref{RibbonPathYFunction}), in the case $q=0$, we have \begin{equation} \label{RibbonPathYFunctionNoPairings} Y_{n,0,m}(\ell_1, \ldots, \ell_n | v^{\textnormal{out}}, v^{\textnormal{in}}):= \sum_{\boldsymbol{\vec{\gamma}} \in \Delta_{n,0,m}(\ell_1, \ldots, \ell_n)} \mathfrak{W}(\boldsymbol{\vec{\gamma}} |v^{\textnormal{out}}, v^{\textnormal{in}}) \end{equation} 
\noindent a $(\textcolor{black}{v^{\textnormal{out}},v^{\textnormal{in}}})$-weighted sum over ribbon paths on $n$ sites of lengths $\ell_1, \ldots, \ell_n$ with no pairings.  By Proposition [\ref{DuhReduceDuhProposition}], for ribbon paths $\boldsymbol{\vec{\gamma}} = (\gamma_1, \ldots, \gamma_n)$ with no pairings, $\boldsymbol{\vec{\gamma}}$ is a disjoint union of $n$ sliding paths and \textcolor{black}{(\ref{RibbonWeightFormula})} is a product of weights \textcolor{black}{(\ref{SlidingWeightFormula})}.  As a consequence,  \begin{equation} \sum_{m=0}^{|| \vec{\ell}||_1} Y_{n,0,m} (\ell_1, \ldots, \ell_n |\textcolor{black}{v^{\textnormal{out}},v^{\textnormal{in}}} ) \textcolor{black}{\ebar}^m = \prod_{a=1}^n \Bigg ( \sum_{m_a=0}^{\ell_a} Y_{1,0,m_a}(\ell_a |\textcolor{black}{v^{\textnormal{out}},v^{\textnormal{in}}}) \textcolor{black}{\ebar}^{m_a} \Bigg ) \end{equation} \noindent The limit (\ref{TransparentThing}) now implies (\ref{AsymptoticFactorizationLemmaFormula}) for the constant \begin{equation} \label{TheConstantWeFoundAsymptoticFactorizationLemma} \mathbf{T}^{\uparrow}_{\ell} ( \textcolor{black}{v^{\textnormal{out}},v^{\textnormal{in}}}; \ebar) = \sum_{m=0}^{\ell} Y_{1,0,m}(\ell |\textcolor{black}{v^{\textnormal{out}},v^{\textnormal{in}}}) \textcolor{black}{\ebar}^{m}. \end{equation} \textcolor{black}{In the case $\ebar \sim \hbar^{1/2}$, unlike \textnormal{(\ref{TransparentThing})}, Theorem [\ref{Theorem3AOE}] now implies that as $\ebar, \hbar \rightarrow 0$, \begin{equation} \mathbb{E} \big [ T_{\ell_1}^{\uparrow}(\ebar, \hbar) \cdots T_{\ell_n}^{\uparrow}(\ebar, \hbar) \big ]  \rightarrow Y_{n,0,0}(\ell_1, \ldots, \ell_n | \textcolor{black}{v^{\textnormal{out}},v^{\textnormal{in}}}) \end{equation} the joint moment asymptotics are governed by the weighted enumeration of ribbon paths on $n$ sites with no pairings ($q=0$) nor slides ($m=0$).  Applying the argument above to Szeg\H{o} paths instead of sliding paths proves the result with $\mathbf{T}_{\ell}^{\uparrow}(\textcolor{black}{v^{\textnormal{out}},v^{\textnormal{in}}}; 0) = Y_{1,0,0}(\ell | \textcolor{black}{v^{\textnormal{out}},v^{\textnormal{in}}})$. $\square$} \end{itemize}

\noindent \textcolor{black}{In \textsection [\ref{DispersiveActionProfileLocation}], under the additional hypothesis $v^{\textnormal{out}} = v^{\textnormal{in}} = v$, we will relate $\mathbf{T}_{\ell}^{\uparrow}(v,v; \ebar)$ and its $\ebar \rightarrow 0$ limit to the dispersive and convex action profiles introduced earlier in \textsection [\ref{SUBSECDispersiveActionProfiles}] and \textsection [\ref{SUBSECConvexActionProfiles}], respectively.  First,} we prove (\ref{ExpectationLLNWeWant}) and (\ref{VarianceLLNWeWant}).  Consider the constants $\mathbf{T}_{\ell}^{\uparrow}(\textcolor{black}{v^{\textnormal{out}},v^{\textnormal{in}}}; \ebar)$ in Lemma [\ref{AsymptoticFactorizationLemmaGood}] and define \begin{equation} \label{NewConstantPolynomialRelation} \mathbf{O}_p (\textcolor{black}{v^{\textnormal{out}},v^{\textnormal{in}}}; \ebar) := \textnormal{KMK}_p \Big ( \mathbf{T}_{1}^{\uparrow}(\textcolor{black}{v^{\textnormal{out}},v^{\textnormal{in}}}; \ebar), \ldots, \mathbf{T}_{p}^{\uparrow}(\textcolor{black}{v^{\textnormal{out}},v^{\textnormal{in}}}; \ebar) \Big ) \end{equation} using the polynomials $\textnormal{KMK}_p$ from Lemma [\ref{KMKPolynomialsLemma}].  Using the polynomial relations (\ref{GreatPolynomialRelation}), (\ref{NewConstantPolynomialRelation}) and repeated application of Lemma [\ref{AsymptoticFactorizationLemmaGood}], for any $p_1, \ldots, p_N$ we have \textcolor{black}{in Regimes I, II, and III} \begin{equation} \label{CoolerAsymptoticFactorization} \mathbb{E} \Big [ O_{p_1}(\ebar, \hbar) \cdots O_{p_N}(\ebar, \hbar) \Big ] \rightarrow \mathbf{O}_{p_1}(\textcolor{black}{v^{\textnormal{out}},v^{\textnormal{in}}}; \ebar) \cdots \mathbf{O}_{p_N} (\textcolor{black}{v^{\textnormal{out}},v^{\textnormal{in}}}; \ebar) .\end{equation}

\noindent This asymptotic factorization (\ref{CoolerAsymptoticFactorization}) implies both (\ref{ExpectationLLNWeWant}) and (\ref{VarianceLLNWeWant}) in the cases $N=1,2$, respectively.  

\subsubsection{\textbf{Step 3: From linear statistics to limit shapes}} \label{SUBSUBSECWhereLimitShapesAt} At this stage, \textcolor{black}{although we have shown in (\ref{ExpectationLLNWeWant}) and (\ref{VarianceLLNWeWant}) the weak law of large numbers for Jack measure linear statistics with arbitrary $v^{\textnormal{out}}, v^{\textnormal{in}}$}, we have not established weak convergence $f_{\lambda}( c | \ebar, \hbar) \rightarrow \mathbf{f}( c | \textcolor{black}{v^{\textnormal{out}},v^{\textnormal{in}}}; \ebar)$ on \textcolor{black}{an explicit} limit shape, since a priori we do not know a function $\mathbf{f}  (c | \textcolor{black}{v^{\textnormal{out}},v^{\textnormal{in}}}; \ebar)$ of $c$ exists so for all $p=1,2,3,\ldots$, \begin{equation} \label{WhatOhReallyIs} \mathbf{O}_p(\textcolor{black}{v^{\textnormal{out}},v^{\textnormal{in}}}; \ebar) = \int_{-\infty}^{+\infty} c^p \tfrac{1}{2} \mathbf{f}'' (c | \textcolor{black}{v^{\textnormal{out}},v^{\textnormal{in}}}; \ebar) dc \end{equation} \noindent where $\mathbf{O}_p(\textcolor{black}{v^{\textnormal{out}},v^{\textnormal{in}}}; \ebar)$ are the constants defined in (\ref{NewConstantPolynomialRelation}).  By (\ref{KMKRelationOT}), it is enough to find $\mathbf{f}( c| \textcolor{black}{v^{\textnormal{out}},v^{\textnormal{in}}}; \ebar)$ so that for all $u \in \mathbb{C} \setminus \mathbb{R}$ one has \begin{equation} \label{ImplicitLimitShapeByKMK} \sum_{\ell=0}^{\infty} \mathbf{T}_{\ell}^{\uparrow}(\textcolor{black}{v^{\textnormal{out}},v^{\textnormal{in}}}; \ebar) u^{-\ell-1} = \textnormal{exp} \Bigg ( \int_{-\infty}^{+\infty} \log \Bigg [ \frac{1}{u-c} \Bigg ] \tfrac{1}{2} \mathbf{f}''(c | \textcolor{black}{v^{\textnormal{out}},v^{\textnormal{in}}}; \ebar) dc \Bigg ). \end{equation}  In the next steps, we find such $\mathbf{f}(c| \textcolor{black}{v^{\textnormal{out}},v^{\textnormal{in}}}; \ebar)$ and describe them explicitly \textcolor{black}{assuming $v^{\textnormal{out}} = v^{\textnormal{in}} = v$.  For the remainder of the proof, $v^{\textnormal{out}} = v^{\textnormal{in}} = v$ and we may abbreviate ``$(v^{\textnormal{out}}, v^{\textnormal{in}}) = (v,v)$'' by $v$.}

\subsubsection{\textbf{Step 4: The limit shape in Regime I is the dispersive action profile}} \label{DispersiveActionProfileLocation} In Regime I, $\hbar \rightarrow 0$ with $\textcolor{black}{\ebar <0}$ fixed.  If it exists, the limit shape in this case is implicitly characterized by (\ref{ImplicitLimitShapeByKMK}) where $\mathbf{T}^{\uparrow}_{\ell}(v; \ebar)$ are the constants in Lemma [\ref{AsymptoticFactorizationLemmaGood}].  In the proof of Lemma [\ref{AsymptoticFactorizationLemmaGood}], we gave an explicit formula for $\mathbf{T}^{\uparrow}_{\ell}(v;\ebar)$ in (\ref{TheConstantWeFoundAsymptoticFactorizationLemma}).  By Proposition [\ref{RibbonPathConnectedPropositionDuh}], ribbon paths on $n=1$ site are all connected, so it follows that $Y_{1,0,m}(\ell | v, v) = W_{1,0,m}(\ell |v,v)$ and so we may rewrite (\ref{TheConstantWeFoundAsymptoticFactorizationLemma}) \textcolor{black}{via (\ref{SlidingThisThing})} as \begin{equation} \label{TheConstantWeFoundAsymptoticFactorizationLemmaREWRITTEN} \mathbf{T}^{\uparrow}_{\ell}(v;\ebar) = \sum_{m=0}^{\ell} W_{1,0,m}(\ell | v,v) \textcolor{black}{\ebar}^m. \end{equation} \noindent Substituting (\ref{TheConstantWeFoundAsymptoticFactorizationLemmaREWRITTEN}) into (\ref{ImplicitLimitShapeByKMK}) yields formula (\ref{DefinitionDispersiveActionProfilesFORMULA}) in the definition of dispersive action profiles, hence the limit shape in Regime I is the dispersive action profile $\mathbf{f}(c | v; \ebar)$ from \textsection [\ref{SECSlidingDispersive}].

\subsubsection{\textbf{Step 5: The limit shape in Regime II is the convex action profile}} In Regime II, $\hbar, \ebar \rightarrow 0$ at a comparable rate $\frac{\ebar^2}{\hbar} = \gamma = (\sqrt{\alpha} - \frac{1}{\sqrt{\alpha}})^2$.  This includes the case $\hbar \rightarrow 0$ with $\ebar=0$ fixed at $\alpha=1$.  By Lemma [\ref{AsymptoticFactorizationLemmaGood}], in Regime II we now have the asymptotic factorization of joint moments \begin{equation} \mathbb{E}[ T_{\ell_1}^{\uparrow} (\ebar, \hbar) \cdots T_{\ell_n}^{\uparrow} (\ebar, \hbar) ] \rightarrow \mathbf{T}_{\ell_1}^{\uparrow} ( v;0) \cdots \mathbf{T}^{\uparrow}_{\ell_n}( v;0) \end{equation} \noindent where $\mathbf{T}^{\uparrow}_{\ell}(v;0)$ are defined by the $\ebar \rightarrow 0$ limit of (\ref{TheConstantWeFoundAsymptoticFactorizationLemmaREWRITTEN}) \textcolor{black}{since $Y_{1,0,0} (\ell | v, v) = W_{1,0,0}(\ell | v,v)$}: \begin{equation} \label{TheConstantWeFoundAsymptoticFactorizationLemmaRegime2} \mathbf{T}^{\uparrow}_{\ell}(v;0) = W_{1,0,0}(\ell |v,v). \end{equation} \noindent Substituting (\ref{TheConstantWeFoundAsymptoticFactorizationLemmaRegime2}) into (\ref{ImplicitLimitShapeByKMK}) yields formula (\ref{DefinitionConvexActionProfilesFORMULA}) in the definition of convex action profiles, hence the limit shape in Regime II is the convex action profile $\mathbf{f} ( c | v;0)$ from \textsection [\ref{SECSzegoConvex}] independent of $\gamma$.  
\subsubsection{\textbf{Step 6: The limit shape in Regime III is determined by Regime I}} In Regime III, $\hbar \rightarrow 0$ with $\textcolor{black}{\ebar >0}$ fixed.  Just as in Step 4, if it exists, the limit shape in this case is implicitly characterized by (\ref{ImplicitLimitShapeByKMK}) where $\mathbf{T}^{\uparrow}_{\ell}(v; \ebar)$ are the constants in (\ref{TheConstantWeFoundAsymptoticFactorizationLemmaREWRITTEN}).  These constants obey the symmetry relation \begin{equation} \label{SymmetryRelation} \mathbf{T}^{\uparrow}_{\ell} ( - v; - \ebar) = (-1)^{\ell} \mathbf{T}^{\uparrow}_{\ell}(v; \ebar) \end{equation}

\noindent as a consequence of (\ref{TheConstantWeFoundAsymptoticFactorizationLemmaREWRITTEN}) and the definition of $W_{1,0,m}(\ell | v, v)$ as a $(v,v)$-weighted sum of sliding paths with $\ell$ steps \textcolor{black}{and $m$ slides}.  Indeed, the $\ell$ steps in a sliding path are either jumps of degree $k \neq 0$, in which case the weight $V_k$ acquires a minus sign upon scaling the specialization $v \mapsto -v$, or slides of height $j$, corresponding to a sign in $\textcolor{black}{\ebar}^m$ upon $\ebar \mapsto -\ebar$.  Using the formula (\ref{ImplicitLimitShapeByKMK}) which determines $\mathbf{f} ( c| v; \ebar)$ from $\mathbf{T}^{\uparrow}_{\ell}(v; \ebar)$, the symmetry relation (\ref{SymmetryRelation}) implies the non-local relation {\small \begin{equation} \label{NonLocalRelationOneThree} \textnormal{exp} \Bigg (\int_{-\infty}^{+\infty}  \log \Bigg [ \frac{1}{(-u)-c} \Bigg ] \tfrac{1}{2} \mathbf{f} '' ( c | -v;- \ebar) dc \Bigg ) = - \textnormal{exp} \Bigg ( \int_{-\infty}^{+\infty} \log \Bigg [ \frac{1}{u-c} \Bigg ] \tfrac{1}{2} \mathbf{f}''(c | v; \ebar) dc \Bigg )  \end{equation}} for all $u \in \mathbb{C} \setminus \mathbb{R}$.  This completes the proof of Theorem [\ref{Theorem1LLN}]. $\square$

\subsection{Example: limit shapes for Poissonized Jack-Plancherel measures} \label{SUBSECJackPlancherelLimitShapes}
In \textsection [\ref{SUBSECJackMeasures}], we saw that Jack measures $M(v_{\textnormal{PL}}, v_{\textnormal{PL}})$ with two identical Plancherel specializations $v_{\textnormal{PL}}$ are the Poissonized Jack-Plancherel measures -- a mixture of the Jack-Plancherel measures $M_d^{\textnormal{PL}}$ of Kerov \cite{Ke4} \textcolor{black}{with law} \textcolor{black}{\textnormal{(\ref{JackPlancherelDoubleHookLaw})}} by a Poisson distribution with intensity $\frac{1}{\hbar}$.  By our discussion in \textsection [\ref{SUBSECCatalanVKLS}] and \textsection [\ref{SUBSECMotzkinNPS}], as $\hbar \rightarrow 0$, our Theorem [\ref{Theorem1LLN}] implies that the limit shapes for such random profiles is the Vershik-Kerov Logan-Shepp profile $\mathbf{f}( c | v_{\textnormal{PL}};0)$ in Regime II \cite{KeVe, LoSh} or the Nekrasov-Pestun-Shatashvili profile $\mathbf{f}(c | v; \ebar)$ in Regimes I and III \cite{NekPesSha}.  We establish the analogous result for $M_d^{\textnormal{PL}}$ in \textsection [\ref{APPENDIXdePoissonization}].




\section{Gaussian fluctuations of Jack measures on profiles} \label{SECProof2}

\noindent In \textsection [\ref{SUBSECCLTProof}] we prove Theorem [\ref{Theorem2CLT}], our second main result for Jack measures $M(v, v)$: in the limit as $\hbar \rightarrow 0$ \textcolor{black}{with $\ebar$ fixed}, the fluctuations of the random profile $f_{\lambda}(c | \ebar, \hbar)$ around the limit shape $\mathbf{f} (c | v; \ebar)$ occur at scale $\hbar^{1/2}$ and converge to a mean zero Gaussian process $\mathbf{G}( c | v; \ebar)$ with explicit covariance given in formula (\ref{SecondCovarianceFormula}) below.  \textcolor{black}{Moreover, if $\ebar, \hbar \rightarrow 0$ so that $\ebar^2 \sim { \gamma \hbar}$, the profile fluctuations occur at the same scale $\hbar^{1/2}$ and are given by $\mathbf{G} (c | v; 0) +  {\gamma}^{1/2} \mathbf{X}( c| v;0)$, a Gaussian process with the same covariance as in the $\ebar = 0$ fixed case except now with a deterministic mean shift $\mathbf{X}(c | v; 0)$.}  In \textsection [\ref{SUBSECJackPlancherelGaussianFluctuations}], we illustrate Theorem [\ref{Theorem2CLT}] for the Poissonized Jack-Plancherel measures $M(v_{\textnormal{PL}}, v_{\textnormal{PL}})$.  In \textsection [\ref{APPENDIXdePoissonization}] we show that a truncated version of the Gaussian process $\mathbf{G}(c | v_{\textnormal{PL}}; \ebar)$ describes the profile fluctuations for the original Jack-Plancherel measures $M_d(v_{\textnormal{PL}}, v_{\textnormal{PL}})$ as $d \rightarrow \infty$.

\subsection{Proof of second main result} \label{SUBSECCLTProof} We now prove Theorem [\ref{Theorem2CLT}] in Steps 1-10 below.

\subsubsection{\textbf{Step 1: Define $N$-decorated ribbon paths and connected $N$-decorated ribbon paths}} Recall $O_p( \ebar, \hbar) |_{\lambda}$ in (\ref{TheOh}) are polynomials in $T^{\uparrow}_{\ell}(\ebar, \hbar)|_{\lambda}$ by Lemma [\ref{KMKPolynomialsLemma}].  By Theorem [\ref{Theorem3AOE}], the $n$th joint cumulants of $T^{\uparrow}_{\ell}(\ebar, \hbar)|_{\lambda}$ can be expressed as $(v^{\textnormal{out}}, v^{\textnormal{in}})$-weighted sums of connected ribbon paths on $n$ sites.  \textcolor{black}{In Steps 1-3, we} establish a similar result for the $N$th joint cumulants of $O_p( \ebar, \hbar) |_{\lambda}$.

\begin{definition} \label{DefinitionNDecoration} An $N$-decoration of $n \in \mathbb{Z}_+$ is $\vec{\nu}=(n_1, \ldots, n_N) \in \mathbb{Z}_{ +}^N$ so $\sum_{j=1}^N n_j = n$. \end{definition}

\begin{definition} \label{DefinitionNDecoratedRibbonPaths} An $N$-decorated ribbon path $(\boldsymbol{\vec{\gamma}}, \vec{\nu})$ is a ribbon path $\boldsymbol{\vec{\gamma}}$ on $n$ sites of lengths $\ell_1, \ldots, \ell_n$ as in \textnormal{Definition [\ref{DefinitionRibbonPaths}]} with an $N$-decoration $\vec{\nu} = (n_1, n_2, \ldots, n_N)$ of $n$.  We relabel \begin{equation} \label{RelabeledSiteLengths} (\ell_1, \ldots, \ell_n) = (\ell_1^{(1)}, \ldots, \ell_{n_1}^{(1)} , \ell_{1}^{(2)}, \ldots,  \ell_{n_2}^{(2)} , \ldots, \ell_{1}^{(N)}, \ldots, \ell_{n_N}^{(N)}) \end{equation} \noindent as ordered sequences and say that the site of length $\ell_{i}^{(j)}$ in $\boldsymbol{\vec{\gamma}}$ is a site decorated by $j \in \{1,2, \ldots, N\}$.
\end{definition}

\noindent \begin{definition} \label{DefinitionReducedGraphAlongDecoration} The reduced graph of an $N$-decorated ribbon path $(\boldsymbol{\vec{\gamma}}, \vec{\nu})$ along $\vec{\nu}$ is the graph with $N$ vertices $j=1, 2, \ldots, N$ and edges $j \rightarrow j'$ for every pairing $\textbf{\textit{p}} = ( \textbf{\textit{e}}, \textbf{\textit{e'}})$ of an edge $\textbf{\textit{e}} \in \textnormal{\textbf{E}}[\gamma_i^{(j)}]$ in a site decorated by $j$ with an edge $\textcolor{black}{\textbf{\textit{e'}} \in} \textnormal{\textbf{E}}[\gamma_{i'}^{(j')}]$ in a site decorated by $j'$.  Here \begin{equation} \boldsymbol{\vec{\gamma}} = (\gamma_1^{(1)}, \cdots, \gamma_{n_1}^{(1)} , \ldots, \gamma_{1}^{(N)}, \ldots, \gamma_{n_N}^{(N)}; \mathbf{p}_1, \ldots, \mathbf{p}_q) \end{equation} we relabel the $n$ sliding paths in the ribbon path $\boldsymbol{\vec{\gamma}}$ using the $N$-decoration $\vec{\nu}$.\end{definition}

\begin{definition} \label{DefinitionConnectedAlongDecoration} An $N$-decorated ribbon path $(\boldsymbol{\vec{\gamma}}, \vec{\nu})$ is connected if its reduced graph along the $N$-decoration $\vec{\nu}$ in \textnormal{Definition [\ref{DefinitionReducedGraphAlongDecoration}]} is connected. \end{definition}

\begin{figure}[htb]
\centering
\includegraphics[width=0.6 \textwidth]{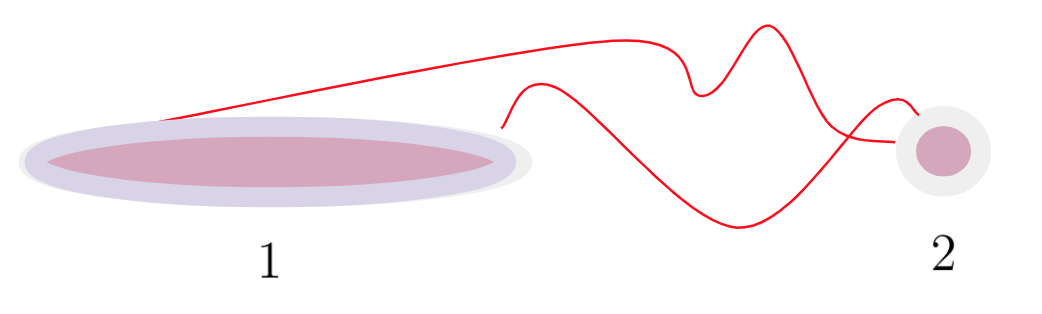}
\caption{The reduced graph of the ribbon path on $n=3$ sites in Figure [\ref{ConnectedRibbonPath2020}] along the $2$-decoration $\vec{\nu} = (n_1, n_2)=(2,1)$.  Since this graph with $N=2$ vertices is connected, $(\boldsymbol{\vec{\gamma}}, \vec{\nu})$ is connected.  Compare to the reduced graph in Figure [\ref{ReducedGraphConnectedRibbonPath2020}].}
\label{ReducedGraphAlongDecorationFIG2020}
\end{figure}

\noindent If an $N$-decorated ribbon path $(\boldsymbol{\vec{\gamma}}, \vec{\nu})$ has $N$-decoration $\vec{\nu} = (n_1, \ldots, n_N)$ with each $n_j =1$, we must have $N=n$ and Definition [\ref{DefinitionReducedGraphAlongDecoration}] reduces to Definition [\ref{DefinitionReducedGraph}].  In this case, the connectedness of $(\boldsymbol{\vec{\gamma}}, \vec{\nu})$ in Definition [\ref{DefinitionConnectedAlongDecoration}] reduces to the original Definition [\ref{DefinitionConnectedRibbonPaths}].

\subsubsection{\textbf{Step 2: Pose weighted enumeration problem for $N$-decorated ribbon paths}} As in \textsection [\ref{SUBSECWeightedEnumerationRibbonPaths}], 
\begin{problem} \label{ProblemDecoratedRibbonPathEnumeration} Fix two specializations $v^{\textnormal{out}}=\{V_k^{\textnormal{out}}\}_{k=1}^{\infty}$, $v^{\textnormal{in}}=\{V_k^{\textnormal{in}}\}_{k=1}^{\infty}$.  
Fix $\vec{p}=(p_1, \ldots, p_N)$, write $|| \vec{p}||_1 = p_1 + \cdots + p_N$, and let $\mathfrak{A}_{\vec{p}}(\vec{\ell}; \vec{\nu})$ be an arbitrary \textcolor{black}{$\mathbb{C}$-valued function of} $\vec{p} = (p_1, \ldots, p_N)$, an $N$-decoration $\vec{\nu}$ of $n$, and $\vec{\ell} = (\ell_1, \ldots, \ell_n)$.  Assume $\mathfrak{A}_{\vec{p}}(\vec{\ell}; \vec{\nu})$ vanishes if $n> || \vec{p}||_1$ or any $\ell_i > \textcolor{black}{|| \vec{p}||_1}$.  Consider the complex weight $\mathfrak{W}(\boldsymbol{\vec{\gamma}} |v^{\textnormal{out}}, v^{\textnormal{in}})$ on ribbon paths in \textcolor{black}{\textnormal{(\ref{RibbonWeightFormula})}}.  \begin{enumerate}
\item Let $\Delta^{\textnormal{N-dec}}_{n,q,m}(\ell_1, \ldots, \ell_n)$ be the infinite set of $N$-decorated ribbon paths $(\boldsymbol{\vec{\gamma}}, \vec{\nu})$ on $n$ sites of lengths $\ell_1, \ldots, \ell_n$ with $q$ pairings and $m$ slides.  Determine \begin{equation} \label{DecoratedRibbonPathYFunction} Y^{\textnormal{dec}}_{N,q,m}(p_1, \ldots, p_N | v^{\textnormal{out}}, v^{\textnormal{in}}):= \sum_{n=1}^{|| \vec{p}||_1} \sum_{\ell_1, \ldots, \ell_n = 1}^{||\vec{p}||_1}   \sum_{(\boldsymbol{\vec{\gamma}} , \vec{\nu}) \in \Delta^{\textnormal{N-dec}}_{n,q,m}(\ell_1, \ldots, \ell_n)}  \mathfrak{A}_{\vec{p}}(\vec{\ell}; \vec{\nu}) \cdot \mathfrak{W}(\boldsymbol{\vec{\gamma}} |v^{\textnormal{out}}, v^{\textnormal{in}}). \end{equation} 
\item Let $\Gamma^{\textnormal{N-dec}}_{n,g,m}(\ell_1, \ldots, \ell_n)$ be the infinite set of connected $N$-decorated ribbon paths $(\boldsymbol{\vec{\gamma}}, \vec{\nu})$ on $n$ sites of lengths $\ell_1, \ldots, \ell_n$ with $n-1+g$ pairings and $m$ slides .  Determine \begin{equation} \label{DecoratedRibbonPathWFunction} W^{\textnormal{dec}}_{N,g,m}(p_1, \ldots, p_N | v^{\textnormal{out}}, v^{\textnormal{in}}):= \sum_{n=1}^{|| \vec{p}||_1} \sum_{\ell_1, \ldots, \ell_n = 1}^{||\vec{p}||_1}  \sum_{(\boldsymbol{\vec{\gamma}} , \vec{\nu}) \in \Gamma_{n,g,m}^{\textnormal{N-dec}}(\ell_1, \ldots, \ell_n)}\mathfrak{A}_{\vec{p}}(\vec{\ell}; \vec{\nu}) \cdot  \mathfrak{W}(\boldsymbol{\vec{\gamma}} |v^{\textnormal{out}}, v^{\textnormal{in}}). \end{equation} 
\end{enumerate}

\end{problem}

\subsubsection{\textbf{Step 3: Polynomial expansion of joint moments and joint cumulants of linear statistics}}

\noindent We now solve a case of Problem [\ref{ProblemDecoratedRibbonPathEnumeration}] and prove that the $N$th joint moments (cumulants) of the linear statistics $O_p(\ebar, \hbar)|_{\lambda}$ are $(\textcolor{black}{v^{\textnormal{out}}, v^{\textnormal{in}}})$-weighted sums of (connected) $N$-decorated ribbon paths.

 \begin{lemma} \label{LinearStatisticsCumulantExpansionLemma} Fix $\vec{p}=(p_1, \ldots, p_N)$.  There is a specific \textcolor{black}{function} $\mathfrak{A}^{\textnormal{KMK}}_{\vec{p}}(\vec{\ell}; \vec{\nu})$ in \textnormal{Problem [\ref{ProblemDecoratedRibbonPathEnumeration}]} so:
 
   \begin{enumerate}
\item The joint moments of the random variables $O_{p}(\ebar, \hbar)$ in \textnormal{(\ref{TheOh})} are polynomials in $\ebar$ and $\hbar$ \begin{equation} \label{LinearStatisticJointMomentsAOE} \mathbb{E} \big [ O_{p_1} (\ebar,\hbar) \cdots O_{p_N} (\ebar, \hbar) \big ] = \sum_{q=0}^{\tfrac{1}{2} || \vec{p}||_1 } \sum_{m=0}^{|| \vec{p}||_1} Y^{\textnormal{dec}}_{N,q,m}(p_1, \ldots, p_N | v^{\textnormal{out}}, v^{\textnormal{in}})  \hbar^q  \textcolor{black}{\ebar}^m  \end{equation} \noindent with $Y^{\textnormal{dec}}_{N,q,m}(\vec{p} \  | v^{\textnormal{out}}, v^{\textnormal{in}})$ counting $N$-decorated ribbon paths $(\boldsymbol{\vec{\gamma}}, \vec{\nu})$ in \textnormal{(\ref{DecoratedRibbonPathYFunction})}.\\
\item The joint cumulants of the random variables $O_p(\ebar, \hbar)$ in \textnormal{(\ref{TransitionStatistics})} are polynomials in $\ebar$ and $\hbar$ \begin{equation} \label{LinearStatisticJointCumulantsAOE} \kappa_N \big ( O_{p_1} (\ebar, \hbar), \ldots, O_{p_N} (\ebar, \hbar) \big ) = \sum_{g=0}^{\tfrac{1}{2} || \vec{p}||_1 -N+1} \sum_{m=0}^{|| \vec{p}||_1}  W^{\textnormal{dec}}_{N,g,m}(p_1, \ldots, p_N | v^{\textnormal{out}}, v^{\textnormal{in}})  \hbar^{N-1+g} \textcolor{black}{\ebar}^m  \end{equation}
\noindent with $W^{\textnormal{dec}}_{N,q,m}(\vec{p} \ | v^{\textnormal{out}}, v^{\textnormal{in}})$ counting connected $N$-decorated ribbon paths $(\boldsymbol{\vec{\gamma}}, \vec{\nu})$ in \textnormal{(\ref{DecoratedRibbonPathWFunction})}. \end{enumerate}  
 \end{lemma}

\noindent Parts (1) and (2) of Lemma [\ref{LinearStatisticsCumulantExpansionLemma}] follow readily from our Theorem [\ref{Theorem3AOE}].
\begin{itemize}
\item \textit{Proof of} Lemma [\ref{LinearStatisticsCumulantExpansionLemma}]: Fix $j$.  By Lemma [\ref{KMKPolynomialsLemma}], there are constants $A_{p_j}(\ell^{(j)}_1, \ldots, \ell^{(j)}_{n_j})$ so \begin{equation} \label{KMKYesButWithLabels} O_{p_j}(\ebar, \hbar) = \sum_{n_j=1}^{p_j} \sum_{\ell_1^{(j)}, \ldots, \ell_{n_j}^{(j)} = 1 }^{p_j} A_{p_j}(\ell^{(j)}_1, \ldots, \ell^{(j)}_{n_j}) T_{\ell_1^{(j)}}^{\uparrow}(\ebar, \hbar) \cdots T_{\ell_{n_j}^{(j)}}^{\uparrow}(\ebar, \hbar)  .\end{equation} 
\noindent In (\ref{KMKYesButWithLabels}), we use indices $n_j$ and $\vec{\ell}^{(j)} = \ell_1^{(j)}, \ldots, \ell_{n_j}^{(j)}$ which reflect the label $j$ of the variable $p_j$ indexing our linear statistic.  Multiplying (\ref{KMKYesButWithLabels}) over $1 \leq j\leq N$ determines $\mathfrak{A}^{\textnormal{KMK}}_{\vec{p}}( \vec{\ell} | \vec{\nu})$: \begin{eqnarray} 
\label{RelabelStep} \prod_{j=1}^N O_{p_j}(\ebar, \hbar) & =&   \sum_{n=1}^{||\vec{p}||_1} \sum_{\ell_1, \ldots, \ell_n =1}^{||\vec{p}||_1} \sum_{\vec{\nu} \in \textnormal{N-dec}(n)}  \mathfrak{A}^{\textnormal{KMK}}_{\vec{p}}(\vec{\ell} \ | \ \vec{\nu}) T_{\ell_1}^{\uparrow}(\ebar, \hbar) \cdots T_{\ell_n}^{\uparrow}(\ebar, \hbar) .\end{eqnarray}

\noindent In (\ref{RelabelStep}), we use (\ref{RelabeledSiteLengths}). The data of the product $\ell_1^{(j)}, \ldots, \ell_{n_j}^{(j)}$ over $j=1,2,\ldots, N$ is the same as the data of $\ell_1, \ldots, \ell_n$ together with an $N$-decoration $\vec{\nu} = (n_1, \ldots, n_N)$ of $n$.  The desired weight $\mathfrak{A}^{\textnormal{KMK}}_{\vec{p}}(\vec{\ell} \ | \ \vec{\nu})$ is determined by this change of variables and $A_{p_j}(\ell^{(j)}_1, \ldots, \ell^{(j)}_{n_j})$ in (\ref{KMKYesButWithLabels}).
\item To prove (1), take the expected value of (\ref{RelabelStep}) and apply Part (1) of Theorem [\ref{Theorem3AOE}].
\item To prove (2), repeat the identical argument in Step 7 in \textsection [\ref{SUBSUBSECStep7}] but with Definition [\ref{DefinitionReducedGraphAlongDecoration}] in place of Definition [\ref{DefinitionReducedGraph}]. This completes the proof of Lemma [\ref{LinearStatisticsCumulantExpansionLemma}]. $\square$
\end{itemize}

\subsubsection{\textbf{Step 4: Weak convergence of rescaled centered linear statistics to Gaussians as $\hbar \rightarrow 0$}} Recall $O_p(\ebar, \hbar)|_{\lambda}$ defined in (\ref{TheOh}) and $\mathbf{O}_p(\textcolor{black}{v^{\textnormal{out}}, v^{\textnormal{in}}}; \ebar)$ in \textcolor{black}{(\ref{NewConstantPolynomialRelation})}.  Define \begin{equation} \label{RescaledLinearStatistics}\varphi_p(\ebar, \hbar)|_{\lambda} := \frac{ O_p(\ebar, \hbar)|_{\lambda} - \mathbf{O}_p(\textcolor{black}{v^{\textnormal{out}}, v^{\textnormal{in}}}; \ebar)}{ \hbar^{1/2}}. \end{equation} \noindent  
\begin{lemma}  \label{CLTOccursLemma} As $\hbar \rightarrow 0$, \textcolor{black}{if $\ebar$ is fixed}, the rescaled centered linear statistics \textnormal{(\ref{RescaledLinearStatistics})} converge weakly \begin{equation} \label{HereYaGoLilWodie} \varphi_p(\ebar, \hbar)|_{\lambda} \rightarrow \mathbf{G}_p(\textcolor{black}{v^{\textnormal{out}}, v^{\textnormal{in}}}; \ebar) \end{equation} \noindent to mean zero Gaussian random variables $\mathbf{G}_p(\textcolor{black}{v^{\textnormal{out}}, v^{\textnormal{in}}}; \ebar)$ whose covariance we denote by \begin{equation}\label{CovarianceFirstPresented} \boldsymbol{\Sigma}_{p_1, p_2} (\textcolor{black}{v^{\textnormal{out}}, v^{\textnormal{in}}}; \ebar) = \textnormal{Cov}[\mathbf{G}_{p_1}(\textcolor{black}{v^{\textnormal{out}}, v^{\textnormal{in}}}; \ebar), \mathbf{G}_{p_2}(\textcolor{black}{v^{\textnormal{out}}, v^{\textnormal{in}}}; \ebar)].\end{equation}  \textcolor{black}{Also, if $\ebar, \hbar \rightarrow 0$ so $\ebar^2 \sim \gamma \hbar$ for some $\gamma  \geq 0$, \textnormal{(\ref{RescaledLinearStatistics})} converge weakly to Gaussian random variables \begin{equation} \varphi_p(\ebar, \hbar)|_{\lambda} \rightarrow \mathbf{G}_p(\textcolor{black}{v^{\textnormal{out}}, v^{\textnormal{in}}}; 0) + \gamma^{1/2} \mathbf{X}_p(\textcolor{black}{v^{\textnormal{out}}, v^{\textnormal{in}}};0) \end{equation} given by the $\mathbf{G}_p(\textcolor{black}{v^{\textnormal{out}}, v^{\textnormal{in}}};0)$ in the $\ebar=0$ case of \textnormal{(\ref{HereYaGoLilWodie})} with deterministic mean shifts $\mathbf{X}_p(\textcolor{black}{v^{\textnormal{out}}, v^{\textnormal{in}}}; 0)$.}  \end{lemma} 
\noindent \textit{Proof of} Lemma [\ref{CLTOccursLemma}]:  \textcolor{black}{First, let $\hbar \rightarrow 0$ with $\ebar$ fixed.}  By definition, it is enough to show \begin{eqnarray} \label{CLTAsymptoticsK1} \kappa_1 \big ( \varphi_p(\ebar, \hbar)|_{\lambda} \big ) & \rightarrow& 0 \\ \label{CLTAsymptoticsK2} \kappa_2 \big ( \varphi_{p_1}(\ebar, \hbar)|_{\lambda} , \varphi_{p_2}(\ebar, \hbar)|_{\lambda} \big )  & \rightarrow & \boldsymbol{\Sigma}_{p_1, p_2} (\textcolor{black}{v^{\textnormal{out}}, v^{\textnormal{in}}}; \ebar) \\ \label{CLTAsymptoticsK3andAbove} \kappa_N(\varphi_{p_1}(\ebar, \hbar)|_{\lambda}, \ldots, \varphi_{p_N}(\ebar, \hbar)|_{\lambda}) & \rightarrow & 0  \end{eqnarray}

\noindent for some quantity $\boldsymbol{\Sigma}_{p_1, p_2} (\textcolor{black}{v^{\textnormal{out}}, v^{\textnormal{in}}}; \ebar)$ and for all $N \geq 3$ in (\ref{CLTAsymptoticsK3andAbove}).  

\begin{itemize}
\item To prove (\ref{CLTAsymptoticsK1}), polynomiality in $\hbar$ of the $N=1$ case of Part (1) of Lemma [\ref{LinearStatisticsCumulantExpansionLemma}] implies \begin{equation} \label{EbarFixedFactYerd} \mathbb{E} [O_p(\ebar, \hbar)|_{\lambda}] = \mathbf{O}_p(\textcolor{black}{v^{\textnormal{out}}, v^{\textnormal{in}}}; \ebar) +  {o}(\hbar^{1/2}) .\end{equation}
\item To establish the remaining limits (\ref{CLTAsymptoticsK2}), (\ref{CLTAsymptoticsK3andAbove}), for $N \geq 2$ we have \begin{equation} \label{ShiftAndScaleTheKappaN} \kappa_N\big (\varphi_{p_1} (\ebar, \hbar)|_{\lambda}, \ldots, \varphi_{p_N}(\ebar, \hbar)|_{\lambda} \big ) = \hbar^{-N/2} \kappa_N\big (O_{p_1} (\ebar, \hbar)|_{\lambda}, \ldots, O_{p_N}(\ebar, \hbar)|_{\lambda} \big )  \end{equation} \noindent using the invariance of $\kappa_N$ under shifts by constants $\mathbf{O}_{p_j}(\textcolor{black}{v^{\textnormal{out}}, v^{\textnormal{in}}}; \ebar)$ for $N \geq 2$ and multi-linearity.  By Part (2) of Lemma [\ref{LinearStatisticsCumulantExpansionLemma}], the $N$th joint cumulants of $O_{p}(\ebar, \hbar)|_{\lambda}$ in \textnormal{(\ref{TheOh})} obey
\begin{equation} \label{TheKeyRelationInLimit}  \kappa_N \big ( O_{p_1} (\ebar, \hbar) |_{\lambda}, \ldots, O_{p_N}( \ebar, \hbar)|_{\lambda} \big ) \sim \hbar^{N-1} \end{equation} 
as $\hbar \rightarrow 0$.  Formulas (\ref{ShiftAndScaleTheKappaN}) and (\ref{TheKeyRelationInLimit}) thus imply \begin{equation} \label{StillTheOne} \kappa_N\big (\varphi_{p_1} (\ebar, \hbar)|_{\lambda}, \ldots, \varphi_{p_N}(\ebar, \hbar)|_{\lambda} \big ) \sim \hbar^{N/2-1}  \end{equation} which proves (\ref{CLTAsymptoticsK2}) and (\ref{CLTAsymptoticsK3andAbove}).  The $\boldsymbol{\Sigma}_{p_1, p_2}(\textcolor{black}{v^{\textnormal{out}}, v^{\textnormal{in}}}; \ebar)$ in (\ref{CLTAsymptoticsK2}) independent of $\hbar$ is \begin{equation}  \label{DefCovAsLimit} \boldsymbol{\Sigma}_{p_1, p_2}(\textcolor{black}{v^{\textnormal{out}}, v^{\textnormal{in}}}; \ebar):= \lim_{\hbar \rightarrow 0} \frac{1}{\hbar} \kappa_2 \big ( O_{p_1}(\ebar, \hbar)|_{\lambda}, O_{p_2}(\ebar, \hbar)|_{\lambda} \big ). \end{equation}
\end{itemize}

\noindent \textcolor{black}{Next, let $\ebar, \hbar \rightarrow 0$ with $\ebar^2 \sim \gamma \hbar$.  By definition, it is enough to show \begin{eqnarray} \label{Regime2CLTAsymptoticsK1} \kappa_1 \big ( \varphi_p(\ebar, \hbar)|_{\lambda} \big ) & \rightarrow& \gamma^{1/2} \mathbf{X}_p (\textcolor{black}{v^{\textnormal{out}}, v^{\textnormal{in}}};0) \\ \label{Regime2CLTAsymptoticsK2} \kappa_2 \big ( \varphi_{p_1}(\ebar, \hbar)|_{\lambda} , \varphi_{p_2}(\ebar, \hbar)|_{\lambda} \big )  & \rightarrow & \boldsymbol{\Sigma}_{p_1, p_2} (\textcolor{black}{v^{\textnormal{out}}, v^{\textnormal{in}}}; 0) \\ \label{Regime2CLTAsymptoticsK3andAbove} \kappa_N(\varphi_{p_1}(\ebar, \hbar)|_{\lambda}, \ldots, \varphi_{p_N}(\ebar, \hbar)|_{\lambda}) & \rightarrow & 0 . \end{eqnarray} In Regime II, \textnormal{(\ref{EbarFixedFactYerd})} no longer holds: the correction to $\mathbb{E}[ O_p(\ebar, \hbar)|_{\lambda}] \rightarrow \mathbf{O}_p(\textcolor{black}{v^{\textnormal{out}}, v^{\textnormal{in}}}; \ebar)$ is $O(\hbar^{1/2})$, not $o(\hbar^{1/2})$.  Precisely, if $\ebar \sim \gamma^{1/2} \hbar^{1/2}$, the $N=1$ case of Part (1) of Lemma [\ref{LinearStatisticsCumulantExpansionLemma}] now implies \begin{equation} \mathbb{E} [ \varphi_p( \ebar, \hbar)|_{\lambda} ] \sim \gamma^{1/2} Y_{1,0,1}^{\textnormal{dec}} (p | \textcolor{black}{v^{\textnormal{out}}, v^{\textnormal{in}}}) \end{equation} \noindent since at $m=1$, $\ebar^m = \ebar^1 = \gamma^{1/2} \hbar^{1/2}$, so (\ref{Regime2CLTAsymptoticsK1}) holds with $\mathbf{X}_p (\textcolor{black}{v^{\textnormal{out}}, v^{\textnormal{in}}};0) = Y_{1,0,1}^{\textnormal{dec}}(p|\textcolor{black}{v^{\textnormal{out}}, v^{\textnormal{in}}})$.  On the other hand, for $N \geq 2$, substituting $\ebar \sim \gamma^{1/2} \hbar^{1/2}$ into Part (2) of Lemma [\ref{LinearStatisticsCumulantExpansionLemma}], by polynomiality, the dominant contributions to the $N$th joint cumulants only come from the terms with both $g=0$ and $m=0$.  For this reason, (\ref{StillTheOne}) still holds and implies (\ref{Regime2CLTAsymptoticsK2}) and (\ref{Regime2CLTAsymptoticsK3andAbove}). $\square$}

\subsubsection{\textbf{Step 5: Computation of covariance I: $2$-decorated ribbon paths}} Next, we have:
\begin{lemma} \label{FirstCovarianceLemma} The limiting covariance $\boldsymbol{\Sigma}_{p_1, p_2}(\textcolor{black}{v^{\textnormal{out}}, v^{\textnormal{in}}}; \ebar)$ in \textnormal{\textcolor{black}{(\ref{CovarianceFirstPresented}),} (\ref{DefCovAsLimit})} has the exact formula 
{ \begin{equation} \label{FirstCovarianceFormula} \boldsymbol{\Sigma}_{p_1, p_2}(\textcolor{black}{v^{\textnormal{out}}, v^{\textnormal{in}}}; \ebar) = \sum_{m=0}^{p_1 + p_2} W_{2,0,m}^{\textnormal{dec}} (p_1, p_2 | \textcolor{black}{v^{\textnormal{out}}, v^{\textnormal{in}}}) \textcolor{black}{\ebar}^m \end{equation}}\noindent where $W_{2,0,m}^{\textnormal{dec}}$ counts connected $2$-decorated ribbon paths with $1$ pairing in \textnormal{(\ref{DecoratedRibbonPathWFunction})} \textcolor{black}{with $\mathfrak{A} = \mathfrak{A}^{\textnormal{KMK}}$}.    \end{lemma}

\begin{itemize}
\item \textit{Proof of} Lemma [\ref{FirstCovarianceLemma}]: Apply Part (2) of Lemma [\ref{LinearStatisticsCumulantExpansionLemma}] at $N=2$ \textcolor{black}{and $\mathfrak{A} = \mathfrak{A}^{\textnormal{KMK}}$}. $\square$
\end{itemize}

\subsubsection{\textbf{Step 6: Computation of covariance II: welding operator}} \label{SUBSUBSECWelding}
\noindent We now give a second more transparent formula for the covariance $\boldsymbol{\Sigma}_{p_1, p_2}(v; \ebar)$ realized above as a combinatorial sum in (\ref{FirstCovarianceFormula}).  \textcolor{black}{To do so, for the remainder of the proof, we assume our Jack measure is defined with $v^{\textnormal{out}} = v^{\textnormal{in}} = v$.}

\begin{definition} Consider the infinite sets of complex variables $v^{{(1)}} = \{V_k^{{(1)}}\}_{k=1}^{\infty}$, $v^{{(2)}} = \{V_k^{{(2)}}\}_{k=1}^{\infty}$.  Write $\overline{V_k}$ for the complex conjugate of $V_k$ and write $\frac{\partial}{\partial V_k}$ and $\frac{\partial}{\partial \overline{V_k}}$ for the two Wirtinger derivatives.  The welding operator is the second-order differential operator 
\begin{equation}\label{WeldingOperatorFormula} \mathcal{K} := \sum_{k=1}^{\infty} k \frac{ \partial^2}{\partial \overline{V_k}^{{(1)}} \partial V_k^{{(2)}}}. \end{equation} \end{definition}

\noindent In the proof below, we use $\mathcal{K}$ to ``weld'' edges $\textit{\textbf{e}}_1, \textit{\textbf{e}}_{2}$ in $2$-decorated ribbon paths to form pairings.

\noindent \begin{lemma} \label{SecondCovarianceLemma} Fix $\ebar \in \mathbb{R}$ and two specializations $v^{{(1)}} = \{V_k^{{(1)}}\}_{k=1}^{\infty}$, $v^{{(2)}} = \{V_k^{{(2)}}\}_{k=1}^{\infty}$ satisfying $\max (|V_k^{{(1)}}|, | V_k^{{(2)}}|) \leq A r^k$ for some $A>0$ and $0<r<1$.  The covariance $\boldsymbol{\Sigma}_{p_1, p_2}(v; \ebar)$ of the Gaussian fluctuations of the rescaled centered linear statistics \textcolor{black}{for fixed $v$} can be computed \textcolor{black}{as} \begin{equation} \label{SecondCovarianceFormula} \boldsymbol{\Sigma}_{p_1, p_2}(v; \ebar) = \mathcal{K} \Bigg (  \int_{-\infty}^{+\infty} c_1^{p_1} \tfrac{1}{2} \mathbf{f}''(c_1| v^{{(1)}}; \ebar) dc_1 \Bigg ) \Bigg (  \int_{-\infty}^{+\infty} c_2^{p_2} \tfrac{1}{2} \mathbf{f}''(c_2 | v^{{(2)}}; \ebar) dc_2 \Bigg ) \Bigg |_{v^{{(1)}} = v^{{(2)}} = v}
\end{equation} \noindent directly from the limit shapes $\mathbf{f}( c | \widetilde{v}; \ebar)$ \textcolor{black}{for $\widetilde{v}$ near $v$} using the welding operator $\mathcal{K}$ in \textnormal{(\ref{WeldingOperatorFormula})}.  \end{lemma}  
 
 \begin{itemize}
 \item \textit{Proof of} Lemma [\ref{SecondCovarianceLemma}]: By (\ref{RescaledLinearStatistics}) and (\ref{WhatOhReallyIs}), we have to prove \begin{equation} \boldsymbol{\Sigma}_{p_1, p_2}(v; \ebar) =\mathcal{K} \mathbf{O}_{p_1}(v^{{(1)}}; \ebar) \mathbf{O}_{p_2} (v^{{(2)}}; \ebar) \Bigg |_{v^{{(1)}} = v^{{(2)}} = v} .\end{equation}
 \noindent By (\ref{ExpectationLLNWeWant}), $N=1$ in Part (2) of Lemma [\ref{LinearStatisticsCumulantExpansionLemma}], and Lemma [\ref{FirstCovarianceLemma}], it is enough to prove \begin{equation} W_{2,0,m}^{\textnormal{dec}}(p_1, p_2 | v, v) = \sum_{m_1 + m_2 = m}  \mathcal{K} \ W_{1,0,m_1}^{\textnormal{dec}}(p_1 |v^{(1)},v^{(1)}) W_{1,0,m_2}^{\textnormal{dec}}(p_2 |v^{(2)},v^{(2)}) \Bigg |_{v^{{(1)}} = v^{{(2)}} = v}.\end{equation}
 \noindent By our assumptions on $v^{{(1)}}, v^{{(2)}}$, we may exchange the derivatives in $\mathcal{K}$ with the infinite sums over ribbon paths with $0$ pairings defining $W^{\textnormal{dec}}_{1,0,m_j}(p_j | v^{(j)}, v^{(j)})$ for $j=1,2$ in (\ref{DecoratedRibbonPathWFunction}) since the summands are polynomials in $\overline{V_k}^{{(1)}}$ and $V_k^{{(2)}}$ and polynomials converge uniformly on compact disks $|V_k^{\textnormal{out}} |, |V_k^{\textnormal{in}}| \leq A r^k$.  \textcolor{black}{In $\mathfrak{W}$ in (\ref{RibbonWeightFormula})}, the welding operator $\mathcal{K}$ removes the contributions $\overline{V_k^{(1)}} V_k^{(2)}$ of an unpaired edges $\textit{\textbf{e}}_1, \textit{\textbf{e}}_2$ of degrees $-k$, $k$ on sites decorated by $1$, $2$ and inserts a factor of $k$.  Combinatorially, the same effect can be achieved by inserting a pairing $\textit{\textbf{p}} = (\textit{\textbf{e}}_1, \textit{\textbf{e}}_2)$ with $\textnormal{size}(\textit{\textbf{p}})=k$ \textcolor{black}{(which is an operation $\boldsymbol{\vec{\gamma}} \rightarrow \boldsymbol{\vec{\gamma}} \cup \mathbf{p}$)} then calculating $\mathfrak{W}(\boldsymbol{\vec{\gamma}} \cup \mathbf{p} | v,v)$ using the weight (\ref{RibbonWeightFormula}).  The result is by definition $W_{2,0,m}^{\textnormal{dec}}(p_1, p_2 | v,v)$. $\square$
 \end{itemize}

\subsubsection{\textbf{Step 7: Gaussian random specializations and welding operators}} \label{SUBSUBSECFGF} In the previous steps, we have established that fluctuations are asymptotically Gaussian and have derived exact formulas for the covariances of the Gaussian random variables $\mathbf{G}_p(v; \ebar)$.  However, at this stage we do not have an intrinsic description of the Gaussian process $\mathbf{G}(c | v; \ebar)$ on the line determined \textcolor{black}{implicitly} by \begin{equation} \label{WouldBeNiceToIBPThis} \mathbf{G}_p(v; \ebar) = \int_{-\infty}^{+\infty} c^p \tfrac{1}{2} \mathbf{G}''(c | v; \ebar) dc . \end{equation} We provide such a description in Steps 8-10.  In order to do so, in this step we define class of Gaussian random specializations $\Phi^{(\mathsf{g})} = \{\Phi^{(\mathsf{g})}_k \}_{k=1}^{\infty}$ for any sequence $\mathsf{g}:=(\mathsf{g}_{11}, \mathsf{g}_{22}, \ldots)$ of positive real numbers $\mathsf{g}_{kk}>0$ and use this to define $\boldsymbol{\varphi} = \{\boldsymbol{\varphi}_k\}_{k=1}^{\infty}$ which appears in the statement of Theorem [\ref{Theorem2CLT}].  For the remainder of this section, fix an arbitrary sequence $\mathsf{g} = (\mathsf{g}_{kk})_{k=1}^{\infty} \subset \mathbb{R}_+$. \begin{definition} Let $\mathbb{A}$ be the set of specializations $v=\{V_k\}_{k=1}^{\infty}$ with finitely-many non-zero $V_k$.  \end{definition} \noindent $\mathbb{A}$ is a complex affine space with global coordinates $(\overline{V_k}, V_k)_{k=1}^{\infty}$ where $\overline{V_k}$ is the complex conjugate.
\begin{definition} \label{GeneralGSpecializationSpace}Let $\mathcal{H}^{(\mathsf{g})}( \mathbb{Z}_+)$ be the Hilbert space completion of $\mathbb{A}$ with respect to
\begin{equation} \mathsf{g}( v^{(1)}, v^{(2)} ) = \sum_{k=1}^{\infty} \mathsf{g}_{kk} \overline{V_k^{(1)} } V_k^{(2)} .\end{equation} \end{definition} 
\noindent The notation $\mathbb{A} \subset \mathcal{H}^{(\mathsf{g})}( \mathbb{Z}_+)$ reflects the fact that $v \in \mathbb{A}$ defines $v: \mathbb{Z}_+ \rightarrow \mathbb{C}$ by $k \mapsto V_k$.

\begin{definition} Let ${\Phi}^{(\mathsf{g})}  = ({\Phi}^{(\mathsf{g})}_{1}, {\Phi}^{(\mathsf{g})}_2, \ldots)$ be a sequence of independent complex rotation invariant Gaussian random variables each with mean $0$ and complex variance $\mathbb{E}[ | {\Phi}^{(\mathsf{g})}_{k} |^2] = \mathsf{g}_{kk}^{-1}$.\end{definition}
\noindent In the literature, ${\Phi}^{(\mathsf{g})}$ is called the ``standard Gaussian in $\mathcal{H}^{(\mathsf{g})}(\mathbb{Z}_+)$'' even though it is almost surely not in $\mathcal{H}^{(\mathsf{g})}(\mathbb{Z}_+)$.  We recall this in the next Lemma [\ref{StandardGaussianLemma}], which follows by direct computation from $\mathbb{E}[| \Phi_k^{(\mathsf{g})} |^2] = \mathsf{g}_{kk}^{-1}$.  For background on Gaussians in Hilbert spaces, see Janson \cite{Janson}.
  \begin{lemma} \label{StandardGaussianLemma} $\Phi^{(\mathsf{g})}$ is almost surely not in $\mathcal{H}^{(\mathsf{g})}(\mathbb{Z}_+)$.  However, for any $\textcolor{black}{\phi} \in \mathcal{H}^{(\mathsf{g})}(\mathbb{Z}_+)$, the series \begin{equation} \label{WellDefinedGaussian} \mathsf{g}( \Phi^{(\mathsf{g})}, \textcolor{black}{\phi}) := \sum_{k=1}^{\infty} \mathsf{g}_{kk} \overline{\Phi^{(\mathsf{g})}_k} \textcolor{black}{\phi}_k \end{equation} is a well-defined rotation invariant Gaussian in $\mathbb{C}$ with mean $0$ and variance $\mathsf{g}(\textcolor{black}{\phi}, \textcolor{black}{\phi}) = || \textcolor{black}{\phi}||_{\mathsf{g}}^2\textcolor{black}{<\infty}$. \end{lemma}

\noindent We now compute the covariance of $\Phi^{(\mathsf{g})}$ using a generalization of $\mathcal{K}$ in (\ref{WeldingOperatorFormula}) which we call $\mathcal{K}^{(\mathsf{g})}$.
\begin{definition} \label{gWeldingOperatorDefinition} The $\mathsf{g}$-welding operator is the second-order differential operator 
\begin{equation}\label{gWeldingOperatorFormula} \mathcal{K}^{(\mathsf{g})}:= \sum_{k=1}^{\infty} \mathsf{g}_{kk}^{-1} \frac{ \partial^2}{\partial \overline{V_k}^{{(1)}} \partial V_k^{{(2)}}}. \end{equation} \end{definition}

\begin{lemma} \label{gWeldingCovarianceLemma} Fix $v \in \mathcal{H}^{(\mathsf{g})}(\mathbb{Z}_+)$ and two functions $O_1, O_2 : \mathcal{H}^{(\mathsf{g})}(\mathbb{Z}_+) \rightarrow \mathbb{R}$ \textcolor{black}{with dense domain of definition including $v$}. 
Let $(\nabla_{\mathsf{g}} O_i)|_v$ denote the gradient at $v$ computed with respect to the inner product $\mathsf{g}$ and regard $\Phi^{(\mathsf{g})}$ as the standard Gaussian in the tangent space $T_v \mathcal{H}^{(\mathsf{g})} (\mathbb{Z}_+)$.  Then for $\textnormal{Re}[a + \textnormal{\textbf{i}} b]=a$, \textcolor{black}{the Gaussian random variables in \textnormal{(\ref{WellDefinedGaussian})} with $\phi_i = ( \nabla_{\mathsf{g}} O_i ) |_v)$ have covariance} \begin{equation} \label{gCovarianceFormula}  \textnormal{Cov} \Bigg [ \mathsf{g} \Big (\Phi^{(\mathsf{g})}, (\nabla_{\mathsf{g}} O_1) \big |_v \Big )  , \mathsf{g} \Big (\Phi^{(\mathsf{g})},  (\nabla_{\mathsf{g}} \textcolor{black}{O_2}) \big |_v \Big ) \Bigg ] = 4 \textnormal{Re} \Bigg [ \mathcal{K}^{(\mathsf{g})} O_1(v^{(1)}) O_2(v^{(2)}) \Bigg |_{v^{(1)} = v^{(2)} = v} \Bigg ]. \end{equation}
\end{lemma}

\begin{itemize}
\item \textit{Proof of} Lemma [\ref{gWeldingCovarianceLemma}]: In $x_k = \textnormal{Re}[V_k], y_k = \textnormal{Im}[V_k]$, $(\nabla_{\mathsf{g}} O )|_v $ has components \begin{equation} (\nabla_{\mathsf{g}} O)|_v = \Big ( \mathsf{g}_{kk}^{-1} \frac{\partial O}{\partial x_k} \Big |_{v} ,  \mathsf{g}_{kk}^{-1} \frac{\partial O}{\partial y_k} \Big |_v \Big )_{k=1}^{\infty}  \end{equation} and defines a tangent vector in $T_v \mathcal{H}^{(\mathsf{g})} ( \mathbb{Z}_+).$ Since $\mathbb{A}$ is affine, $\Phi^{(\mathsf{g})}$ is also the standard Gaussian on $T_v\mathcal{H}^{(\mathsf{g})}(\mathbb{Z}_+)$, so by Lemma [\ref{StandardGaussianLemma}] the left-hand side of (\ref{gCovarianceFormula}) is

 \begin{equation} \label{WhatWeWantLateInTheGame} \mathsf{g}\big ( (\nabla_{\mathsf{g}} O_1)|_v, (\nabla_{\mathsf{g}} O_2)|_v \big ) =  \sum_{k=1}^{\infty} \mathsf{g}_{kk}^{-1}  \Bigg ( \frac{\partial O_1}{\partial x_k} \frac{\partial O_2}{\partial x_k} +  \frac{\partial O_1}{\partial y_k} \frac{\partial O_2}{\partial y_k} \Bigg ). \end{equation}  On the other hand, the right-hand side of (\ref{gCovarianceFormula}) is
 \begin{equation} 4 \textnormal{Re} \ \mathcal{K}^{(\mathsf{g})} O_1(v^{(1)}) O_2(v^{(2)}) \Bigg |_{v^{(1)} = v^{(2)} = v} = 4 \textnormal{Re} \  \sum_{k=1}^{\infty} \mathsf{g}_{kk}^{-1} \frac{\partial O_1}{\partial \overline{V_k}} \frac{\partial O_2}{\partial {V_k}} \end{equation} 
 
 \noindent Writing $\frac{\partial}{\partial V_k} = \frac{1}{2} \Big ( \frac{\partial}{\partial x_k} - \textbf{i} \frac{\partial}{\partial y_k} \Big ) $ and $\frac{\partial}{\partial \overline{V_k}} = \frac{1}{2} \Big ( \frac{\partial}{\partial x_k} + \textbf{i} \frac{\partial}{\partial y_k} \Big ) $ and taking $\textnormal{Re}$ gives (\ref{WhatWeWantLateInTheGame}). $\square$

\end{itemize}

\noindent \textcolor{black}{The random specialization} $\boldsymbol{\varphi} = \{\boldsymbol{\varphi}_k\}_{k=1}^{\infty}$ in Theorem [\ref{Theorem2CLT}] is the special case of the above.
 \begin{definition} \label{RelabelingDefinitionOneHalf} \textcolor{black}{In the case $\mathsf{g}_{kk} = k^{2s}$ for some $s \in \mathbb{R}$, write $\mathcal{H}^{{s}}( \mathbb{Z}_+)$ instead of $\mathcal{H}^{(\mathsf{g})}( \mathbb{Z}_+)$.}  In particular, for $s = -1/2$ so $\mathsf{g}_{kk} = k^{-1}$, write $\mathcal{H}^{{-1/2}}( \mathbb{Z}_+)$ and $( \cdot, \cdot )_{-1/2}$ instead of $\mathcal{H}^{(\mathsf{g})}( \mathbb{Z}_+)$ and $\mathsf{g}(\cdot, \cdot)$, write $\boldsymbol{\varphi} =  \{\boldsymbol{\varphi}_{k} \}_{k=1}^{\infty}$ instead of $\Phi^{(\mathsf{g})}= \{\Phi^{(\mathsf{g})}_{k} \}_{k=1}^{\infty}$, and write $\mathcal{K}$ instead of $\mathcal{K}^{(\mathsf{g})}$ as in \textnormal{(\ref{WeldingOperatorFormula})}.
\end{definition}

\begin{lemma} \label{JackItIsLemma} $v \in \mathcal{H}^{{-1/2}}( \mathbb{Z}_+)$ if and only if the Jack measure $M(v,v)$ is well-defined. \end{lemma}
\begin{itemize}
\item \textit{Proof of} Lemma [\ref{JackItIsLemma}]: Immediate from \textnormal{(\ref{DiagonalRegularityCondition})}. $\square$
\end{itemize}

\begin{lemma} \label{WhatIsThePhi} $\boldsymbol{\varphi}_k$ are independent Gaussians with $\mathbb{E}[\boldsymbol{\varphi}_k] = \mathbb{E}[\boldsymbol{\varphi}_k^2] = 0$ and $\mathbb{E}[| \boldsymbol{\varphi}_k|^2] = k$. \end{lemma} 
\begin{itemize}
\item \textit{Proof of} Lemma [\ref{WhatIsThePhi}]: \textcolor{black}{For $\mathsf{g}_{kk} = k^{2s}$ and $s = - \tfrac{1}{2}$,} $\mathbb{E}[ |\Phi_k^{(\mathsf{g})} |^2] = \mathsf{g}_{kk}^{-1} = (k^{-1})^{-1} = k$. $\square$
\end{itemize}

\noindent Lemma [\ref{gWeldingCovarianceLemma}] specializes to the following result involving the original welding operator $\mathcal{K}$.

\begin{lemma} \label{ThirdCovarianceLemma} Fix $v \in \mathcal{H}^{-1/2}(\mathbb{Z}_+)$ and two functions $O_1, O_2 : \mathcal{H}^{-1/2}(\mathbb{Z}_+) \rightarrow \mathbb{R}$ \textcolor{black}{with dense domain of definition including $v$}. 
Let $(\nabla_{-1/2} O_j)|_v$ denote the gradient at $v$ computed with respect to the inner product $(\cdot, \cdot)_{-1/2}$ and regard $\boldsymbol{\varphi}$ as the standard Gaussian in $T_v \mathcal{H}^{-1/2} (\mathbb{Z}_+)$.  Then \begin{equation} \label{SpecialgCovarianceFormula}  \textnormal{Cov} \Bigg [ \Big (\boldsymbol{\varphi}, (\nabla_{{-1/2}} O_1)\big |_v \Big )_{-1/2}  , \Big (\boldsymbol{\varphi},  (\nabla_{{-1/2}} O_1) \big |_v \Big )_{-1/2} \Bigg ] = 4\textnormal{Re} \Bigg [ \mathcal{K} O_1(v^{(1)}) O_2(v^{(2)}) \Bigg |_{v^{(1)} = v^{(2)} = v} \Bigg ]. \nonumber \end{equation}

\end{lemma}

\noindent The superscript $-1/2$ in Definition [\ref{RelabelingDefinitionOneHalf}] is carefully chosen notation.  Let $\mathbb{T} = \mathbb{R} / 2 \pi \mathbb{Z}$ be the unit circle.  The map $v \mapsto \sum_{k=1}^{\infty} (V_k e^{- \textbf{i} kx} + \overline{V_k} e^{+\textbf{i} kx})$ defines an isometry $\mathcal{H}^{-1/2}(\mathbb{Z}_+) \rightarrow \mathcal{H}_0^{-1/2} ( \mathbb{T}; \mathbb{R}) $ from our space of Jack measure specializations to the homogeneous $L^2$-Sobolev space of mean-zero real-valued functions on $\mathbb{T}$ of Sobolev regularity $s = - \frac{1}{2}$.  Across this isometry, $\boldsymbol{\varphi}$ is sent to \begin{equation} \label{FGFRelevantHere} \boldsymbol{\varphi} (x) = \sum_{k=1}^{\infty} \big ( \boldsymbol{\varphi}_k e^{- \textbf{i} kx} + \overline{\boldsymbol{\varphi}_k} e^{+\textbf{i} kx} \big ) \end{equation} the fractional Gaussian field $\mathsf{FGF}_{-1/2}( \mathbb{T}^1)$ on $\mathbb{T}$ for $s=-\tfrac{1}{2}$ with Hurst index $s - \tfrac{d}{2} = -\tfrac{1}{2} - \tfrac{1}{2} = -1$.  This $\boldsymbol{\varphi}(x)$ is almost surely of Sobolev regularity $\widetilde{s}$ for any $\widetilde{s}< -1$ but almost surely never of $\widetilde{s}=-1$.  For a survey of fractional Gaussian fields, see Lodhia-Sheffield-Sun-Watson \cite{FGFsurvey}.
 
\subsubsection{\textbf{Step 8: The Gaussian process in Regime I}} \label{FlexThisStep8GG} We now describe $\textcolor{black}{\mathbf{G}}(c| v; \ebar)$ in Regime I \textcolor{black}{as the {push-forward} of the Gaussian process $\boldsymbol{\varphi}$ along a linear map $d \mathbf{f}_{\ebar}|_v$.  Recall that if $\mathcal{L} : \mathcal{H}_1 \rightarrow \mathcal{H}_2$ is a linear map between real vector spaces and $\Phi \in \mathcal{H}_1$ is Gaussian, the \textit{push-forward $\mathcal{L}_{*}\Phi$ of $\Phi$ along $\mathcal{L}$} is the Gaussian in $\mathcal{H}_2$ satisfying $\upalpha ( \mathcal{L}_{*} \Phi) = (\upalpha \circ \mathcal{L}) ( \Phi)$ for all linear functionals $\upalpha: \mathcal{H}_2 \rightarrow \mathbb{R}$.}

\begin{proposition} \label{RegimeICLTDescription} For $\textcolor{black}{\ebar <0}$, let $\textcolor{black}{\mathbf{f}_{\ebar} }: \textcolor{black}{\widetilde{v}} \mapsto \mathbf{f}(\textcolor{black}{c} | \textcolor{black}{\widetilde{v}}; \ebar)$ be the map \textcolor{black}{which sends a specialization $\widetilde{v} = \{ \widetilde{V}_k \}_{k=1}^{\infty}$ with bounded symbol \textnormal{(\ref{Symbol})} to the dispersive action profile in \textnormal{Definition [\ref{DefinitionDispersiveActionProfiles}]}.  Let $\boldsymbol{\varphi} = \{\boldsymbol{\varphi}_k\}_{k=1}^{\infty}$ be the Gaussian random specialization from \textnormal{Definition [\ref{RelabelingDefinitionOneHalf}]}.}  Then \textcolor{black}{the Gaussian process defined implicitly by \textnormal{(\ref{WouldBeNiceToIBPThis})} in terms of $\mathbf{G}_p( v; \ebar)$ with covariance $\boldsymbol{\Sigma}_{p_1, p_2}(v; \ebar)$ in \textnormal{(\ref{SecondCovarianceFormula})} is} \begin{equation} \label{PushForwardPresentationIsItOk} \textcolor{black}{\mathbf{G}}(c | v; \ebar) =  \tfrac{1}{2} \big ( d\textcolor{black}{\mathbf{f}_{\ebar}}|_v\big )_* \boldsymbol{\varphi} \end{equation} the push-forward of $\tfrac{1}{2}\boldsymbol{\varphi}$ along the \textcolor{black}{linear map given by the} differential $d \textcolor{black}{\mathbf{f}_{\ebar}}|_{v}$ at $\widetilde{v}=v$.
\end{proposition}
\begin{itemize}
\item \textit{Proof of} Proposition [\ref{RegimeICLTDescription}]: \textcolor{black}{First, we argue that $(d \mathbf{f}_{\ebar} |_v )_{*} \boldsymbol{\varphi}$ in (\ref{PushForwardPresentationIsItOk}) is well-defined.  In \cite{Ke1}, Kerov defined a ``{profile}'' $f: \mathbb{R} \rightarrow \mathbb{R}$ to be any real-valued function of $c \in \mathbb{R}$ which is $1$-Lipshitz and has $\int_0^{\pm \infty} (1 \mp f'(c) ) \frac{dc}{1 + |c|} < \infty$.  Let $B$ denote the set of all profiles and $L^{\infty}(\mathbb{Z}_+)$ the set of all specializations with bounded symbol (\ref{Symbol}).  In \cite{Moll1}, we proved that \begin{equation} \label{SZeroMomentMap} \mathbf{f}_{\ebar} : L^{\infty}(\mathbb{Z}_+) \rightarrow B.\end{equation} However, $\boldsymbol{\varphi}$ is almost surely not in $L^{\infty} (\mathbb{Z}_+)$ but rather only in $\mathcal{H}^{\widetilde{s}} (\mathbb{Z}_+)$ for any $\widetilde{s} < -1$ in Definition [\ref{RelabelingDefinitionOneHalf}] as discussed in \textsection [\ref{SUBSUBSECFGF}].  Fortunately, the differential of (\ref{SZeroMomentMap}) at $\widetilde{v} = v$
\begin{equation} d \mathbf{f}_{\ebar} |_v : T_{v} L^{\infty}(\mathbb{Z}_+) \rightarrow T_{\mathbf{f}( c | v; \ebar)} B \end{equation} is
 {linear} and thus has a unique extension to a linear map on the larger domain $T_v \mathcal{H}^{\widetilde{s}} (\mathbb{Z}_+)$ for any $\widetilde{s} < -1$ since the embedding $T_v L^{\infty}(\mathbb{Z}_+) \subset T_v \mathcal{H}^{\widetilde{s}}(\mathbb{Z}_+)$ is dense, so $(d \mathbf{f}_{\ebar} |_v )_* \boldsymbol{\varphi}$ is well-defined.  To prove (\ref{PushForwardPresentationIsItOk}), set $\mathbf{O}_p(\widetilde{v}; \ebar) = \mathbf{O}_p (\widetilde{v}, \widetilde{v}; \ebar)$ in (\ref{WhatOhReallyIs}).  For any profile $b(c) \textcolor{black}{\in B}$, \begin{equation} \label{ThatLittleT} t_p(b) = \int_{-\infty}^{+\infty} c^p \tfrac{1}{2} b''(c) dc \end{equation} is a linear functional densely-defined on $B$ with respect to the convergence of profiles specified in Kerov \cite{Ke1}.  By the definition of $\mathbf{f}_{\ebar}$,} $\mathbf{O}_p( \widetilde{v} ; \ebar) = (t_p \circ \textcolor{black}{\mathbf{f}_{\ebar}}) ( \widetilde{v})$.  \textcolor{black}{Since $t_p$ and $\mathbf{f}_{\ebar}$ are densely-defined as discussed, $\mathbf{O}_p (\widetilde{v}; \ebar)$ has dense domain of definition in $\mathcal{H}^{\widetilde{s}}(\mathbb{Z}_+)$ including $v$ for all $\widetilde{s} < -1$.}  By the chain rule \textcolor{black}{and unique extensions of linear functionals on dense subsets}, the differential $d \mathbf{O}_p (\cdot ; \ebar)|_v : T_v \mathcal{H}^{\widetilde{s}} (\mathbb{Z}_+) \rightarrow \mathbb{R}$ is \textcolor{black}{the composition} \begin{equation} \label{TheCompositionYaDig} d \mathbf{O}_p (\cdot ; \ebar)|_v = ( d t_p \circ d \mathbf{f}_{\ebar} \big |_v) ( \cdot) .\end{equation} \noindent \textcolor{black}{Since $t_p(b)$ is linear in $b$, $d t_p = t_p$ and hence (\ref{WouldBeNiceToIBPThis}) reads \begin{equation} \label{LookAtWhereWeAre}  d t_p \big ( \mathbf{G}( c | v; \ebar) \big )  = \mathbf{G}_p (v ; \ebar) .\end{equation} On the other hand, by (\ref{TheCompositionYaDig}) and the definition of gradients in terms of differentials, \begin{equation} \label{LookAtWhereWeAre2}  ( d t_p \circ d \mathbf{f}_{\ebar}|_v ) ( \boldsymbol{\varphi}) = \Big ( \boldsymbol{\varphi} , \big ( \nabla_{-1/2} \mathbf{O}_p ( \cdot ; \ebar ) \big ) |_v \Big )_{-1/2}. \end{equation} \noindent Finally, Lemma [\ref{SecondCovarianceLemma}], the reality of covariances, $O_i = \mathbf{O}_{p_i}( \cdot ; \ebar)$ in Lemma [\ref{ThirdCovarianceLemma}] for $i=1,2$, and formulas (\ref{LookAtWhereWeAre}) and (\ref{LookAtWhereWeAre2}) imply the equality in joint law \begin{equation}  dt_p ( \mathbf{G}(c | v; \ebar) ) = (dt_p \circ d \mathbf{f}_{\ebar}|_v ) ( \tfrac{1}{2} \boldsymbol{\varphi}) \end{equation} \noindent of the real centered Gaussian random variables indexed by $p=1,2,3,\ldots$ with covariance $\boldsymbol{\Sigma}_{p_1, p_2} ( v; \ebar)$.  If $\upalpha$ is any linear combination of $dt_p$, $\upalpha( \mathbf{G} ( c |v ; \ebar) ) = ( \upalpha \circ d \mathbf{f}_{\ebar}|_v ) ( \tfrac{1}{2} \boldsymbol{\varphi})$.}  $\square$

\end{itemize} 

\pagebreak

\subsubsection{\textbf{Step 9: The Gaussian process in Regime II}} We now describe $\mathbf{G}(c| v; 0)$ in Regime II.

\begin{proposition} \label{RegimeIICLTDescription} For $\ebar =0$, let $\textcolor{black}{\mathbf{f}_0}: \textcolor{black}{\widetilde{v}} \mapsto \mathbf{f} (\textcolor{black}{c} | \textcolor{black}{\tilde{v}}; 0)$ be the map \textcolor{black}{which sends a specialization $\widetilde{v} = \{ \widetilde{V}_k\}_{k=1}^{\infty}$ with bounded symbol \textnormal{(\ref{Symbol})} to the convex action profile in \textnormal{Definition [\ref{DefinitionConvexActionProfiles}]}}.  Then \begin{equation} \mathbf{G}(c | v; 0) = \tfrac{1}{2}\big ( \textcolor{black}{d \mathbf{f}_0 |_v}\big )_* \boldsymbol{\varphi} \end{equation} is the push-forward of $\tfrac{1}{2}\boldsymbol{\varphi}$ along the differential $\textcolor{black}{d \mathbf{f}_0 |_{v}}$ at $v$. \end{proposition}

\begin{itemize}
\item \textit{Proof of} Proposition [\ref{RegimeIICLTDescription}]: \textcolor{black}{Repeat the proof} of Proposition [\ref{RegimeICLTDescription}] \textcolor{black}{with $\ebar=0$}. $\square$
\end{itemize}

\noindent In this case $\ebar=0$, we know from (\ref{PushForwardFormulas}) and (\ref{WhatOhReallyIs}) that \begin{equation} \mathbf{O}_p (v; 0) = \int_0^{2 \pi}  v(x)^p \tfrac{dx}{2\pi} \end{equation} where $v(x)$ is defined by (\ref{Symbol}).  \noindent By Lemma [\ref{SecondCovarianceLemma}], the covariance of $\mathbf{G}(c| v; 0)$ is thus \begin{eqnarray} \boldsymbol{\Sigma}_{p_1, p_2}(v; 0) &=& \mathcal{K} \mathbf{O}_{p_1}(v^{(1)}; 0) \mathbf{O}_{p_2} (v^{(2)} ; 0) \Big |_{v^{(1)} = v^{(2)} = v}  \nonumber \\ &=& \sum_{k=1}^{\infty} k \frac{\partial}{\partial \overline{V_k^{(1)}} } \Bigg ( \int_{0}^{2\pi} \big ( v^{(1)}(x_1) \big )^{p_1}  \tfrac{dx_1}{2\pi} \Bigg ) \frac{\partial}{\partial V_k^{(2)}}\Bigg ( \int_{0}^{2\pi} \big ( v^{(2)}(x_2) \big )^{p_2}  \tfrac{dx_2}{2\pi} \Bigg ) \Bigg |_{v^{(1)} = v^{(2)} = v} \\ &=& \label{NiceBDCovarianceMatch} \sum_{k=1}^{\infty} k \Bigg ( \int_{0}^{2\pi} p_1 v(x_1)^{p_1-1} e^{+ \textbf{i} kx_1 } \tfrac{dx_1}{2\pi} \Bigg ) \Bigg (  \int_{0}^{2\pi} p_2 v(x_2)^{p_2-1} e^{- \textbf{ i} k x_2}   \tfrac{dx_2}{2\pi}  \Bigg )  \end{eqnarray} \noindent which coincides with the covariance of the Gaussian process in Breuer-Duits \cite{BreuerDuits}.  Finally, let $\hbar \rightarrow 0$ and $\ebar \rightarrow 0$ at a comparable rate $\frac{\ebar^2}{\hbar} = \gamma \in \mathbb{R}$.  \textcolor{black}{By Lemma [\ref{LinearStatisticsCumulantExpansionLemma}], in Regime II} the fluctuations at scale $\hbar^{1/2}$ are asymptotically Gaussian with covariance independent of $\gamma$ \textcolor{black}{and a mean shift determined by $\mathbf{X}_p(v;0) = Y_{1,0,1}^{\textnormal{dec}}(p | v,v) = W_{1,0,1}^{\textnormal{dec}}(p | v,v)$, the weighted enumeration of} $1$-decorated ribbon paths on $N=1$ site with $g=0$ pairings and $m=1$ slide since at $m=1$, $\ebar^m = \ebar^1 \sim \gamma^{1/2} {\hbar}^{1/2}$ is the order of the fluctuations.  In conclusion, the limit is now $\mathbf{G}(c | v; 0) + \sqrt{\gamma} \cdot \mathbf{X}(c|v;0)$ for $\mathbf{X}(c| v;0)$ satisfying \begin{equation} \label{FormulaForRegimeIIGaussianFieldMeanShift} \int_{-\infty}^{+\infty} c^p \tfrac{1}{2} \mathbf{X}''(c|v;0)dc   = \textcolor{black}{\frac{d}{d \ebar}} \Bigg ( \int_{-\infty}^{+\infty} c^{p} \tfrac{1}{2} \mathbf{f}''( c | v; \ebar) dc \Bigg ) \Bigg |_{\ebar = 0}. \end{equation}

\subsubsection{\textbf{Step 10: The Gaussian process in Regime III}} We now describe $\mathbf{G}(c| v; \ebar)$ in Regime III.

\begin{proposition} \label{RegimeIIICLTDescription} For $\textcolor{black}{\ebar>0}$, $\mathbf{G}(c| v; \ebar)$ is determined by $\mathbf{G}(c | -v; - \ebar)$ by the non-local relation \begin{equation} \label{NonLocalRelationForCLT} \boldsymbol{\Sigma}_{p_1, p_2} (-v; - \ebar) = (-1)^{p_1 + p_2} \boldsymbol{\Sigma}_{p_1, p_2} ( v; \ebar) .\end{equation} 
\end{proposition}

\begin{itemize}
\item \textit{Proof of} Proposition [\ref{RegimeIIICLTDescription}]. By (\ref{KMKRelationOT}), $\textnormal{KMK}_p(T_1^{\uparrow}, \cdots, T_{p}^{\uparrow})$ contain the term $T_{\ell_1}^{\uparrow} \cdots T_{\ell_n}^{\uparrow}$ only if $\ell_1 + \cdots + \ell_n=p$.  Combining (\ref{NewConstantPolynomialRelation}) and (\ref{SymmetryRelation}) implies for all $p$ the relation $\mathbf{O}_p(-v; -\ebar) = (-1)^{p} \mathbf{O}_p(v; \ebar)$ equivalent to (\ref{NonLocalRelationOneThree}) which is non-local in $v$.  The desired relationship between the covariances now follows from (\ref{SecondCovarianceFormula}). $\square$
\end{itemize}
\subsection{Example: Gaussian fluctuations for Poissonized Jack-Plancherel measures} \label{SUBSECJackPlancherelGaussianFluctuations}

\noindent In \textsection [\ref{SUBSECJackPlancherelLimitShapes}], we saw Theorem [\ref{Theorem1LLN}] gives limit shapes for Poissonized Jack-Plancherel measures $M(v_{\textnormal{PL}}, v_{\textnormal{PL}})$. We now illustrate Theorem [\ref{Theorem2CLT}] for these same measures $M(v_{\textnormal{PL}}, v_{\textnormal{PL}})$ in Regime II.  In \textsection{[\ref{APPENDIXdePoissonization}]}, we discuss the corresponding fluctuation results for Jack-Plancherel measures $M_d^{\textnormal{PL}}$ \textcolor{black}{with law} \textcolor{black}{\textnormal{(\ref{JackPlancherelDoubleHookLaw})}}.
\subsubsection{Covariance of Poissonized Jack-Plancherel fluctuations in Regime II} \label{ContactIvanovOlshanski} For $v = v_{\textnormal{PL}}$, (\ref{Symbol}) specializes to $v_{\textnormal{PL}}(x) = 2 \cos x$ and so (\ref{NiceBDCovarianceMatch}) implies {\small \begin{equation} \label{ThisIsNiceAreYouDone} \boldsymbol{\Sigma}_{p_1, p_2} ( v_{\textnormal{PL}}; 0)  = \sum_{k=1}^{\infty} k \Bigg ( \int_{0}^{2\pi} p_1 \big ( 2\cos x_1\big )^{p_1-1} e^{+ \textbf{i} kx_1 } \tfrac{dx_1}{2\pi} \Bigg ) \Bigg (  \int_{0}^{2\pi} p_2 \big ( 2 \cos x_2 \big )^{p_2-1} e^{- \textbf{ i} k x_2}   \tfrac{dx_2}{2\pi}  \Bigg ) .\end{equation} } \noindent Integrating by parts in $c$ in (\ref{WouldBeNiceToIBPThis}), \begin{equation} \mathbf{G}_p(v; \ebar) = - \int_{-\infty}^{+\infty} pc^{p-1}  \tfrac{1}{2}\mathbf{G}' (c | v; \ebar) dc \end{equation} and so by (\ref{CovarianceFirstPresented}) and (\ref{ThisIsNiceAreYouDone}), for any polynomials $Q_1(c), Q_2(c)$ and $\pm_j = (-1)^{j-1}$, {\small 
\begin{equation} \label{DontYouWantToPlugInChebychev} \textnormal{Cov} \Bigg ( \Bigg \{ \int_{-\infty}^{+\infty} Q_j(c_j) \tfrac{1}{2} \mathbf{G}'(c_j | v_{\textnormal{PL}}, 0) dc_j \Bigg \}_{j=1}^2 \Bigg ) = \sum_{k=1}^{\infty} k  \prod_{j=1}^2 \Bigg ( \int_{0}^{2\pi} Q_j ( 2\cos x_j\big ) e^{\pm_j \textbf{i} kx_j } \tfrac{dx_j}{2\pi} \Bigg ).\end{equation}}\noindent As a consequence, let $k_1, k_2 \in \{1,2,3,\ldots\}$ be independent of the index $k$ in (\ref{DontYouWantToPlugInChebychev}) and choose $Q_j(c)= \mathsf{T}_{k_j}(\tfrac{c}{2})$ for $j=1,2$ where $\mathsf{T}$ is the Chebychev polynomial of the first kind \begin{equation} \mathsf{T}_{k_j} (\cos x_j) = \cos (k_j x_j) .\end{equation} \noindent In this case, (\ref{DontYouWantToPlugInChebychev}) implies that the Gaussian random variables \begin{equation} \label{LoveItHereTheWorldIsRound}- \int_{-\infty}^{+\infty} \mathsf{T}_k (\tfrac{c}{2}) \tfrac{1}{2} \mathbf{G} '(c | v_{\textnormal{PL}}; 0) dc \end{equation} \noindent are independent real Gaussian random variables with variance $k$.  Using the relation \begin{equation} 2 \frac{d}{dc} \mathsf{T}_k( \tfrac{c}{2}) = k \mathsf{U}_{k-1} ( \tfrac{c}{2}) \end{equation}
\noindent to the Chebychev polynomials of the second kind, (\ref{LoveItHereTheWorldIsRound}) and integration by parts imply \begin{equation} \label{PoissonizedIvanovOlshanski} \int_{-\infty}^{+\infty} \mathsf{U}_{k-1} (\tfrac{c}{2} ) \mathbf{G}(c | v_{\textnormal{PL}}, 0) dc \end{equation} are independent real Gaussian random variables with variance $1/k$.  For $k \geq 2$, this is coincides exactly with the description of Kerov's Gaussian fluctuation result for Plancherel measures \cite{Ke5} in Theorem 7.1 of Ivanov-Olshanski \cite{IvOl}.  We account for this difference at $k=1$ in \textsection [\ref{APPENDIXdePoissonization}].

\subsubsection{Mean shift of Poissonized Jack-Plancherel fluctuations in Regime II} For $v = v_{\textnormal{PL}}$, we calculate the mean shift $\mathbf{X}(c | v_{\textnormal{PL}}; 0)$ present if $\ebar, \hbar \rightarrow 0$ so that $\ebar^2 \sim \gamma \hbar$.  For $u \in \mathbb{C} \setminus \mathbb{R}$, define \begin{equation} S_{\pm}(u) = \frac{u \pm \sqrt{ 4 - u^2}}{2} \end{equation}  At $\ebar=0$, recall that the Stieltjes transform of Wigner's semi-circle in (\ref{CatalanQuadraticTing}) is \begin{equation} \label{FirstThingToUseForMeanShift} \mathbf{T}^{\uparrow}(u |  v_{\textnormal{PL}}; 0) = S_-(u). \end{equation} \begin{lemma} \label{OneSlideLemma} For varying $\ebar$, $\mathbf{T}^{\uparrow}(u | v_{\textnormal{PL}}; \ebar)$ in \textnormal{(\ref{MotzkinQuadraticTing})} satisfies \begin{equation} \label{SecondThingToUseForMeanShift} \textcolor{black}{\frac{d}{d \ebar}} \mathbf{T}^{\uparrow}(u | v_{\textnormal{PL}}; \ebar) \Bigg |_{\ebar = 0} = \frac{S_-(u)^2}{ 4 - u^2}. \end{equation} 
\end{lemma}
\begin{itemize}
\item \textit{Proof of} Lemma [\ref{OneSlideLemma}]: Since $\mathbf{T}^{\uparrow}(u | v; \ebar)$ is equal to both sides of (\ref{DefinitionDispersiveActionProfilesFORMULA}), the left side of (\ref{SecondThingToUseForMeanShift}) extracts the weighted count of \textcolor{black}{certain} sliding paths $\gamma \in W_{1,0,1}(\ell)$ on $n=1$ site of length $\ell$ and exactly $m=1$ slide.  \textcolor{black}{For $v = v_{\textnormal{PL}}$ with $V_k^{\textnormal{PL}} = \delta(k-1)$, the only} sliding path\textcolor{black}{s} $\gamma$ \textcolor{black}{which appear} are the concatenation of a simple walk in $\mathbb{Z}_{\bullet}^{2}$ from the origin to height $j$, a slide at height $j$, and a simple walk back down to $(\ell,0)$.  To quantify this observation, for $j=0,1,2,\ldots$ let $R_{0,j}({\ell})$ be the number of lattice paths in $\mathbb{Z}_{\bullet}^2$ which start at $(0,0)$, end at $(\ell, j)$, and take only jumps $\textbf{\textit{e}}=(1, \pm 1)$ of degree $\pm 1$.  $R_{0,0}(\ell)$ are the Catalan numbers $W_{1,0,0}(\ell | v_{\textnormal{PL}}, v_{\textnormal{PL}})$ and $R_{0,j}(0) = \delta(j-0)$.  The generating function is easily shown to be \begin{equation} \label{NiceGeneratingFunction202020} \sum_{\ell=0}^{\infty} R_{0,j}(\ell) u^{-\ell-1} = S_-(u)^{j+1}. \end{equation}  \noindent By (\ref{DefinitionDispersiveActionProfilesFORMULA}) and the decomposition of $\gamma \in \Gamma_{1,0,1}(\ell)$ above, the left-hand side of (\ref{SecondThingToUseForMeanShift}) is \begin{equation} \sum_{\ell=0}^{\infty} u^{-\ell-1} \sum_{\gamma \in \Gamma_{1,0,1}(\ell)} \mathfrak{W}( \gamma | v_{\textnormal{PL}}, v_{\textnormal{PL}}) = \sum_{j=1}^{\infty} j S_{-}(u)^{2j+2} = \frac{S_-(u)^2}{ (S_-(u)^{-1} - S_-(u) )^2}  \end{equation} \noindent which proves (\ref{SecondThingToUseForMeanShift}) since $S_{\pm}(u)^{-1} = S_{\mp}(u)$ and $(S_+(u) - S_-(u))^2 = {4-u^2}$.  $\square$
\end{itemize}
\noindent By (\ref{FormulaForRegimeIIGaussianFieldMeanShift}), (\ref{DefinitionDispersiveActionProfilesFORMULA}), (\ref{FirstThingToUseForMeanShift}), and (\ref{SecondThingToUseForMeanShift}), 
\begin{eqnarray}  \int_{-\infty}^{+\infty}\log \Bigg [ \frac{1}{u-c}\Bigg ] \tfrac{1}{2} \mathbf{X}''(c | v_{\textnormal{PL}}; 0) dc &=& \frac{d}{d \ebar} \Bigg ( \int_{-\infty}^{+\infty} \log \Bigg [ \frac{1}{u-c} \Bigg ]  \tfrac{1}{2} \mathbf{f}''(c | v_{\textnormal{PL}}; \ebar) dc \Bigg ) \Bigg |_{\ebar = 0} \\  &=& \frac{d}{d \ebar} \Bigg ( \log \mathbf{T}^{\uparrow} (u | v_{\textnormal{PL}}; \ebar) \Bigg ) \Bigg |_{\ebar = 0} \\ &=& \label{DoPartialFractionsYourself} \frac{S_-(u)}{4 - u^2}  .
\end{eqnarray}

\noindent As a consequence, the mean shift for Poissonized Jack-Plancherel measures in Regime II coincides with the shift for Jack-Plancherel measures discovered by Do{\l}{\k{e}}ga-F{\'e}ray \cite{DoFe}:\begin{equation} \label{DFMeanShift} \mathbf{X}(c | v_{\textnormal{PL}};0) = \begin{cases} - \tfrac{1}{2\pi} \arcsin \big ( \tfrac{c}{2} \big ), \ \ \ \ \ \ \ \ \ |c| \leq 2 \\ 0 \ \ \ \ \ \ \ \ \ \ \ \  \ \ \ \ \ \ \  \ \ \ \ \ \    \ \ \ \ \ \ \ \ |c| > 2\end{cases}\end{equation} We discuss this coincidence in \textsection [\ref{APPENDIXdePoissonization}].  Note that (\ref{DFMeanShift}) is discontinuous in $c$, so its distributional derivatives produce delta functions at $c = \pm 2$ reflected in the partial fraction expansion of (\ref{DoPartialFractionsYourself}).



\section{Comments and comparison with previous results} \label{SECComments}

\subsection{Jack measures} \label{SUBSECJackMeasureComments} \textcolor{black}{We first discuss our results for Jack measures with generic specializations $v^{\textnormal{out}} = \{ V_k^{\textnormal{out}} \}_{k=1}^{\infty}, v^{\textnormal{in}} =  \{ V_k^{\textnormal{in}} \}_{k=1}^{\infty}$ and generic Jack parameter $\alpha >0$ encoded by $\ebar \in \mathbb{R}$ as in \textsection [\ref{SUBSECJackPolynomials}].}
 
\subsubsection{Previous work} \label{SUBSUBSECPreviousWork} The Regime II cases of our main results in this paper, Theorem [\ref{Theorem1LLN}] and Theorem [\ref{Theorem2CLT}], as well as our technical Theorem [\ref{Theorem3AOE}] relating Jack measures and ribbon paths, first appeared in the year 2015 in versions 1 and 2 of the author's unpublished arXiv preprint \cite{Moll0}.  Later in the year 2017 in version 3 of \cite{Moll0}, we proved these results in all three Regimes I, II, III as a corollary of much more general limit shape and Gaussian fluctuation results for stochastic processes that arise in geometric quantization and the semi-classical analysis of coherent states.  These general results in version 3 of \cite{Moll0} do not assume classical or quantum integrability and will appear in \cite{ChangMoll4}.  In this \textcolor{black}{current} paper, we prove Theorems [\ref{Theorem1LLN}] and [\ref{Theorem2CLT}] in Regimes I, II, III by the combinatorial method of moments as we originally did for Regime II in versions 1, 2 of \cite{Moll0}. 

\subsubsection{Distinct specializations}  \textcolor{black}{Although Theorems [\ref{Theorem1LLN}] and [\ref{Theorem2CLT}] are stated for Jack measures with $v^{\textnormal{out}} = v^{\textnormal{in}} = v$, many of our results above in \textsection [\ref{SECRibbonAnisotropic}], \textsection [\ref{SECProof1}], \textsection [\ref{SECProof2}] are proven for $v^{\textnormal{out}}, v^{\textnormal{in}}$ possibly distinct.  This level of generality in Theorem [\ref{Theorem3AOE}] is necessary to establish a connection to the work of Chapuy-Do{\l}{\k{e}}ga \cite{ChapuyDolega2020} in \textsection [\ref{APPENDIXRibbonPaths}] as discussed in \textsection [\ref{SUBSECRibbonPathsComments}].  In our proofs of Theorems [\ref{Theorem1LLN}] and [\ref{Theorem2CLT}], the assumption $v^{\textnormal{out}} = v^{\textnormal{in}} = v$ is used to identify limit shapes with the dispersive and convex action profiles from \cite{GerardKappeler2019, Moll1}.  In \cite{GerardKappeler2019, Moll1}, these profiles are defined through the spectral theory of generalized Toeplitz operators with symbols $\sum_{k=1}^{\infty} \overline{V_k^{\textnormal{out}} } e^{\textbf{i} kx} + V_k^{\textnormal{in}} e^{- \textbf{i} kx}$ which are self-adjoint if and only if the symbol is real-valued as in (\ref{Symbol}), i.e. if $v^{\textnormal{out}} = v^{\textnormal{in}} = v$.  In this way, the problem in \textsection [\ref{SUBSUBSECWhereLimitShapesAt}] of finding limit shapes $\mathbf{f}( c | v^{\textnormal{out}}, v^{\textnormal{in}}, \ebar)$ for Jack measures with distinct specializations reduces to the spectral theory of the generalized Toeplitz operators in \cite{GerardKappeler2019, Moll1} with complex-valued symbols.  For example, such analysis might be applied to describe the known limit shapes of Schur-Weyl measures of Biane \cite{Bi3}, since these measures on partitions can be defined as the marginals of the Schur measures $M(v_{\textnormal{PL}}, v_{\star})$ with $\ebar=0$ defined by distinct Plancherel and principal specializations.}

\subsubsection{Jack non-negative specializations} \textcolor{black}{In \cite{KeOkOl}, Kerov-Okounkov-Olshanski classified the set of all specializations $v^{\dagger}= \{V^{\dagger}_k\}_{k=1}^{\infty}$ of the ring of symmetric functions which are non-negative on Jack polynomials $P_{\lambda}(V^{\dagger}_1, V^{\dagger}_2, \ldots | \ebar, \hbar) \geq 0$ for every partition $\lambda$.  In the Definition [\ref{DefinitionJackMeasures}] of Jack measures $M(v^{\textnormal{out}}, v^{\textnormal{in}})$, it is not necessary that either of the two specializations $v^{\textnormal{out}}, v^{\textnormal{in}}$ be chosen from the Jack non-negative specializations $v^{\dagger}$ from \cite{KeOkOl}.  In fact, the $V_k^{\dagger}$ in \cite{KeOkOl} depend on both $\ebar$ and $\hbar$, while in the current paper we only consider specializations with $V_k$ independent of $\ebar$ and $\hbar$.  It would be interesting to adapt the methods in this paper to study the extremal Gibbs measures on the Young graph with Jack edge multiplicities in \cite{KeOkOl}, since these measures are marginals of the Jack measures $M(v_{\textnormal{PL}}, v^{\dagger})$ defined by distinct Plancherel and Jack non-negative specializations.}

\subsubsection{Limit shapes for Jack measures} \textcolor{black}{For the limit shapes $\mathbf{f} ( c | v ; \ebar)$ derived in our Theorem [\ref{Theorem1LLN}], $\mathbf{f} (c | v; \ebar) - | c|$ has compact support in two cases: in Regime II for generic $v$ (as proven above) or in Regimes I and III for non-generic ``multi-phase'' specializations $v$ (as proven in \cite{Moll6}; see \textsection [\ref{SUBSECzMeasureComments}]).  If $\mathbf{f}( c| v; \ebar) - |c|$ has compact support, the convergence in Theorem [\ref{Theorem1LLN}] could be improved as follows.  The $T^{\uparrow}_{\ell} (\ebar, \hbar)|_{\lambda}$ in (\ref{TransitionStatistics}) are the $\ell^{\textnormal{th}}$ moments of the random transition measures $\tau_{\lambda}^{\uparrow}(c | \ebar, \hbar)$ of Kerov \cite{Ke4}.  Likewise, the limit shapes $\mathbf{f}( c | v; \ebar)$ have deterministic transition measures $\boldsymbol{\tau}^{\uparrow}( c | v; \ebar)$ defined in \cite{Moll1} following Kerov \cite{Ke1}.  If $\mathbf{f} (c | v; \ebar) - |c|$ has compact support, then so does $\boldsymbol{\tau}^{\uparrow}( c | v; \ebar)$ and Weierstrass approximation implies weak convergence $\tau_{\lambda}^{\uparrow}(c | \ebar, \hbar) \rightarrow \boldsymbol{\tau}^{\uparrow}( c | v; \ebar)$ in probability.  If we also knew $f_{\lambda}(c | \ebar, \hbar) - |c|$ was compactly supported near the support of $\mathbf{f} (c | v; \ebar) - |c|$ with high probability as $\hbar \rightarrow 0$, Kerov's Markov-Kre\u{\i}n correspondence \cite{Ke1} as stated in \textsection 2.2.2 of Sodin \cite{SodinKMK} would imply uniform convergence $f_{\lambda}(c | \ebar, \hbar) \rightarrow \mathbf{f}( c | v; \ebar)$ in probability.  To control the random supports of $f_{\lambda}(c | \ebar, \hbar) - |c|$, one strategy would be to prove \textit{large deviations principles} for Jack measures in Regimes I, II, or III.  Our proof of Theorem [\ref{Theorem1LLN}] does not rely on any large deviations principles for Jack measures.  To our knowledge, even candidates for the rate functions generalizing the hook integral of Vershik-Kerov \cite{KeVe} and Logan-Shepp \cite{LoSh} are not known for Schur measures with generic specializations.  For details regarding the hook integral and uniform convergence in the case of Plancherel measures, see Romik \cite{RomikBOOK} and Sodin \cite{SodinKMK}.}

\subsubsection{Gaussian fluctuations for Jack measures} \textcolor{black}{Our proof of Theorem [\ref{Theorem2CLT}] -- most notably the labeling of the quantities $W_{n,g,m}$ and construction of a ``welding operator'' $\mathcal{K}$ in \textsection [\ref{SUBSUBSECWelding}] and \textsection [\ref{SUBSUBSECFGF}] -- is heavily inspired by the loop equation analysis of continuous $\beta$-ensembles in Borot-Guionnet \cite{BtGu1}.  For background, see Guionnet \cite{GuionnetBOOK}.    In \textsection [\ref{FlexThisStep8GG}], $\mathcal{K}$ allows us to prove that the $v$ and $\ebar$ dependence of the Jack measure global Gaussian fluctuations $\mathbf{G}(c | v; \ebar)$ are completely determined by a deterministic linear map $d \mathbf{f}_{\ebar} |_v$ according to (\ref{PushForwardPresentationIsItOk}).  In (\ref{PushForwardPresentationIsItOk}), $\boldsymbol{\varphi}$ is independent of $v$ and $\ebar$ with symbol (\ref{FGFRelevantHere}) given by the spatial gradient of the log-correlated Gaussian field on the unit circle $\mathbb{T}$ which is itself the restriction of any 2D Gaussian free field to a small circle.  For this reason, our Theorem [\ref{Theorem2CLT}] has the same form as several global fluctuation results for continuous $\beta$-ensembles.  For a survey, see Borodin \cite{BoCMI}.  In fact, our results in \cite{ChangMoll4, Moll1} reveal that it is not a coincidence that our Theorem [\ref{Theorem2CLT}] can be stated through (\ref{FGFRelevantHere}), since the Hilbert space completion of the ring of symmetric functions with its Hall inner product (\ref{NormIntro}) can be defined without completion by the Segal-Bargmann construction \textit{starting from} the fractional Gaussian field (\ref{FGFRelevantHere}) as in Janson \cite{Janson}.}

\subsection{Schur measures} \label{SUBSECSchurMeasureComments} At $\ebar=0$, Jack measures recover the Schur measures of Okounkov \cite{Ok1}.  
\subsubsection{Limit shapes for Schur measures} The $\ebar=0$ case of Regime II of our Theorem [\ref{Theorem1LLN}] was originally proven in \textsection 4.2.1 of Okounkov \cite{Ok2}.  The term ``convex action profile'' was later introduced in \cite{Moll1}.  The proof in \cite{Ok2} uses the realization \cite{Ok1} of Schur measures as determinantal point processes and the method of steepest descent.  Our method of proof in this $\ebar=0$ case is new.

\subsubsection{Gaussian fluctuations for Schur measures} \noindent The Gaussian process $\mathbf{G}(c | v;0)$ in the $\ebar=0$, $\gamma=0$ case of Regime II of our Theorem [\ref{Theorem2CLT}] has also appeared in work of Breuer-Duits \textnormal{\cite{BreuerDuits}} on Borodin's biorthogonal ensembles \textnormal{\cite{Bo1}}.  To see this, identify our $v(x)$ in (\ref{NiceBDCovarianceMatch}) with the symbol of the Laurent matrix in the main Theorem 2.1 in \cite{BreuerDuits}.  By contrast, $v(x)$ in (\ref{Symbol}) is the symbol of the Toeplitz matrix in Szeg\H{o}'s First Theorem which we used to define our $\mathbf{f}(c | v; 0)$ in \textsection [\ref{SECSzegoConvex}].

\subsection{Plancherel measures} \label{SUBSECPlancherelMeasureComments} At $\ebar =0$ and $v_{\textnormal{PL}}=\{1,0,0,\ldots\}$, Jack measures $M(v_{\textnormal{PL}},v_{\textnormal{PL}})$ are the Poissonized Plancherel measures, a mixture of the {Plancherel measures} of the symmetric group $S(d)$ on $d$ elements by a Poisson distribution of intensity $\frac{1}{\hbar}$.  \textcolor{black}{For their law, set $- \varepsilon_2 = \varepsilon_1 = \hbar^{1/2}$ in (\ref{JackPlancherelDoubleHookLaw}).}  For background on $d \rightarrow \infty$ asymptotics of Plancherel measures, see Baik-Deift-Sudan \cite{BaikDeiftSuidanBOOK} and Romik \cite{RomikBOOK}.  Note that $v_{\textnormal{PL}}$ are called ``exponential specializations'' in \textsection 7.8 of Stanley \cite{StanleyVol2}.

\subsubsection{Limit shapes for Plancherel measures} As $d \rightarrow \infty$, a well-known result of Vershik-Kerov \cite{KeVe} and Logan-Shepp \cite{LoSh} is that the profile of a typical partition $\lambda$ sampled from the Plancherel measures approaches the convex limit shape (\ref{FormulaVKLS}).  In \textsection [\ref{SUBSECJackPlancherelLimitShapes}], we saw that the case $v=v_{\textnormal{PL}}$, $\ebar=0$ of our Theorem [\ref{Theorem1LLN}] in Regime II proves the same result for the Poissonized Plancherel measures.  In \textsection [\ref{APPENDIXdePoissonization}], we ``dePoissonize'' our methods above to give a new proof of the result in \cite{KeVe, LoSh}.

\subsubsection{Gaussian fluctuations for Plancherel measures} In \textsection [\ref{ContactIvanovOlshanski}], we evaluated the covariance of the mean zero Gaussian process $\mathbf{G}(c | v_{\textnormal{PL}}; 0)$ in Regime II of our Theorem [\ref{Theorem2CLT}] assuming $\ebar=0$ and $v= v_{\textnormal{PL}}$.  We showed that these Gaussian fluctuations of Poissonized Plancherel measures are slightly different from the Gaussian fluctuations of Plancherel measures found by Kerov \cite{Ke5} as described in Theorem 7.1 of Ivanov-Olshanski \cite{IvOl}.  In \textsection [\ref{APPENDIXdePoissonization}], we explain this discrepancy and give a new proof of the result in \cite{Ke5, IvOl}.  The Chebychev polynomials in \cite{IvOl} emerge naturally in our above argument in \textsection [\ref{ContactIvanovOlshanski}] from the assumption $v_{\textnormal{PL}}(x) = 2 \cos x$ on the symbol $v(x)$ in (\ref{Symbol}).

\subsection{Jack-Plancherel measures} \label{SUBSECJackPlancherelMeasureComments} \noindent At a Plancherel specialization $v_{\textnormal{PL}}=\{1,0,0,\ldots\}$, Jack measures $M(v_{\textnormal{PL}},v_{\textnormal{PL}})$ are the Poissonized Jack-Plancherel measures, a mixture of {Jack-Plancherel measures} $M_d^{\textnormal{PL}}$ of Kerov \cite{Ke4} \textcolor{black}{with law} \textcolor{black}{\textnormal{(\ref{JackPlancherelDoubleHookLaw})}} by a Poisson distribution of intensity $\frac{1}{\hbar}$.  \textcolor{black}{To verify this claim, one may use \textsection VI.10 in Macdonald \cite{Mac}, take the $z,z' \rightarrow \infty$ limit of (1.1) in Borodin-Olshanski \cite{BoOl1}, or set $(\varepsilon_2, \varepsilon_1) = (- \epsilon_2, \epsilon_1)$ in (3.4) in the survey by Okounkov \cite{Ok2}}.  The Jack-Plancherel measures $M_d^{\textnormal{PL}}$ are the simplest case of the Nekrasov-Okounkov measures on partitions \cite{NekOk}; see \textsection 3.2 of \cite{NekOk}.  In this paper we use the conventions (\ref{Convention1}), (\ref{Convention2}) from \cite{NekOk}.  For work independent of \cite{NekOk} on Jack-Plancherel measures, see Fulman \cite{Fulman}, Hora-Obata \cite{HoOb}, Matsumoto \cite{MatsumotoJack, MatsumotoJack2}, Do{\l}{\k{e}}ga-F{\'e}ray \cite{DoFe}, Do{\l}{\k{e}}ga-{\'S}niady \cite{DoSni}, Borodin-Gorin-Guionnet \cite{BoGoGu}, and Guionnet-Huang \cite{GuionnetHuang}.

\subsubsection{Limit shapes for Jack-Plancherel measures in Regime II} In \textsection [\ref{SUBSECJackPlancherelLimitShapes}], we applied our first Theorem [\ref{Theorem1LLN}] to prove that if $\hbar \rightarrow 0$ and $\ebar \rightarrow 0$ at a comparable rate $\gamma = \big ( \sqrt{\alpha} - \frac{1}{\sqrt{\alpha}} \big )^2$ with fixed Jack parameter $\alpha$ in (\ref{Convention2}), the Poissonized Jack-Plancherel measures form a limit shape $\mathbf{f}( c | v_{\textnormal{PL}};0)$ independent of $\alpha$ that coincides with the limit shape (\ref{FormulaVKLS}) found by Vershik-Kerov \cite{KeVe} and Logan-Shepp \cite{LoSh}.  For Jack-Plancherel measures $M_d^{\textnormal{PL}}$, this result is due to Do{\l}{\k{e}}ga-F{\'e}ray \cite{DoFe}.  We give a new proof of their result by a short dePoissonization argument in \textsection [\ref{APPENDIXdePoissonization}].  

\subsubsection{Limit shapes for Jack-Plancherel measures in Regimes I and III} In \textsection [\ref{SUBSECJackPlancherelLimitShapes}], we applied our first Theorem [\ref{Theorem1LLN}] to prove that if $\hbar \rightarrow 0$ with $\ebar \neq 0$ fixed so that $\textcolor{black}{\alpha \rightarrow 0}$ or $\textcolor{black}{\alpha \rightarrow \infty}$ in (\ref{Convention2}), the Poissonized Jack-Plancherel measures $M(v_{\textnormal{PL}}, v_{\textnormal{PL}})$ form a limit shape $\mathbf{f}(c | v_{\textnormal{PL}}; \ebar)$ that coincides with a Nekrasov-Pestun-Shatashvili profile in \cite{NekPesSha} as we discuss in \textsection [\ref{SUBSECMotzkinNPS}].  While $M(v_{\textnormal{PL}}, v_{\textnormal{PL}})$ appear in Nekrasov-Okounkov \cite{NekOk}, these limit shapes $\mathbf{f}(c | v_{\textnormal{PL}}; \ebar)$ do not appear in \cite{NekOk} since in this paper the authors consider only the Seiberg-Witten limit $\varepsilon_1, \varepsilon_2 \rightarrow 0$ and assume $\varepsilon_1 + \varepsilon_2 = 0$ which is $\ebar = 0$ in (\ref{OmegaVariables}), (\ref{Ebar}), (\ref{Hbar}).  For Jack-Plancherel measures, this same limit shape was independently discovered by Do{\l}{\k{e}}ga-{\'S}niady \cite{DoSni}.  We give a new proof of this result from \cite{NekPesSha, DoSni} in \textsection [\ref{APPENDIXdePoissonization}].

\subsubsection{Gaussian fluctuations for Jack-Plancherel measures in Regime II} In \textsection [\ref{SUBSECJackPlancherelGaussianFluctuations}], we applied our second Theorem [\ref{Theorem2CLT}] to prove that if $\hbar \rightarrow 0$ and $\ebar \rightarrow 0$ at a comparable rate $\gamma = \big ( \sqrt{\alpha} - \frac{1}{\sqrt{\alpha}} \big )^2$ with fixed Jack parameter $\alpha$ in (\ref{Convention2}), the fluctuations of the Poissonized Jack-Plancherel measures at scale $\hbar^{1/2}$ converge to a Gaussian process $\mathbf{G}(c | v_{\textnormal{PL}}; 0) +\big ( \sqrt{\alpha} - \frac{1}{\sqrt{\alpha}} \big )\mathbf{X}(c | v_{\textnormal{PL}};0)$ where $\mathbf{G}(c| v_{\textnormal{PL};0})$ is independent of $\alpha$ and $\mathbf{X}(c | v_{\textnormal{PL}};0)$ is a deterministic mean shift.  This result is nearly identical to the Gaussian fluctuations of the ordinary Jack-Plancherel measures found by Do{\l}{\k{e}}ga-F{\'e}ray \cite{DoFe}.  We give a new proof of their result by a short dePoissonization argument in \textsection [\ref{APPENDIXdePoissonization}].  

\subsubsection{Gaussian fluctuations for Jack-Plancherel measures in Regimes I and III} For Poissonized Jack-Plancherel measures, our Theorem [\ref{Theorem2CLT}] implies that the fluctuations of the random partition profile around the piecewise-linear $\mathbf{f}( c |v_{\textnormal{PL}}; \ebar)$ in \cite{NekPesSha} converge to a Gaussian process $\mathbf{G}(c | v_{\textnormal{PL}}; \ebar)$ whose covariance can be computed in formula (\ref{SecondCovarianceFormula}) directly from the limit shapes.  For the ordinary Jack-Plancherel measures, such Gaussian fluctuations around piecewise-linear limit shapes were first discovered by Do{\l}{\k{e}}ga-{\'S}niady \cite{DoSni}.  We give a new proof of their result in this case by a dePoissonization argument in \textsection [\ref{APPENDIXdePoissonization}].  However, note that Do{\l}{\k{e}}ga-{\'S}niady establish limit shapes and Gaussian fluctuations for a large class of random partitions in \cite{DoSni} by delicate analysis of asymptotics of Jack characters.  For background, see {\'S}niady \cite{Sniady, SniadyJackCharacters}.  Relating our methods for Jack measures to the models and techniques in Do{\l}{\k{e}}ga-{\'S}niady \cite{DoSni} would be especially exciting in light of the connection between our ribbon paths and the work of Chapuy-Do{\l}{\k{e}}ga \cite{ChapuyDolega2020} in \textsection [\ref{APPENDIXRibbonPaths}].

\subsection{z-measures} \label{SUBSECzMeasureComments} Our original motivation for studying Jack measures $M(v^{\textnormal{out}}, v^{\textnormal{in}})$ \textcolor{black}{with generic complex-valued specializations $v^{\textnormal{out}}, v^{\textnormal{in}}$} was to gain insight into properties of the \textit{$z$-measures} on partitions introduced by Borodin-Olshanski \cite{BoOl1}.  As noted in \textsection 3 of \cite{BoOl1}, the principal series of these $z$-measures are marginals of Jack measures $M(v_{\star}, v_{\star})$ defined by two identical principal specializations $v_{\star}$ \textcolor{black}{which are complex-valued}.  In a subsequent work \cite{Moll6}, we combine the results in this paper with the finite gap analysis of the $n$-phase solutions of the classical Benjamin-Ono equation in \cite{GerardKappeler2019, Moll2, Moll1} to define and analyze \textcolor{black}{\textit{$n$-phase specializations} of the ring of symmetric functions and corresponding} \textit{$n$-phase $z$-measures} for any $n=1,2,3,\ldots$ which generalize the $z$-measures in \cite{BoOl1} at $n=1$.  \textcolor{black}{For $n$-phase specializations $v$, $n_{*}(v)=n$ in Proposition [\ref{PropositionSlidingPaths}] and the dispersive action profile $\mathbf{f}(c | v; \ebar)$ has $2n+1$ interlacing local extrema}.  Remarkably, the original $z$-measures in \cite{BoOl1} correspond directly to the $1$-phase periodic traveling wave solution of the Benjamin-Ono equation.  \textcolor{black}{As we will show in \cite{Moll7}, the coherent state quantization of this classical periodic traveling wave on the circle is a superposition of interacting quantum periodic traveling waves with wavenumbers $j=1,2,3,\ldots$ whose quantized actions \cite{Moll2} are discrete random variables $2\pi \hbar ( \lambda_j - \lambda_{j+1})$ governed by the mixed $z$-measures on partitions $\lambda = (\lambda_j)_{j=1,2,3,\ldots}$ in \cite{BoOl1}.}

\subsection{Related measures}  \label{SUBSECRelatedMeasuresComments} Jack measures are a special case of the Macdonald measures introduced by Borodin-Corwin \cite{BoCo}.  For generalizations of Jack-Plancherel measures that are a priori distinct from Jack measures, see Borodin-Gorin-Guionnet \cite{BoGoGu} and Do{\l}{\k{e}}ga-{\'S}niady \cite{DoSni}.  In \cite{HuangJack}, Huang extended the framework developed by Bufetov-Gorin \cite{BufetovGorinSchur} and found a large class of measures on partitions that are distinct from Jack measures and yet are also amenable to asymptotic analysis by means of the Nazarov-Sklyanin operators \cite{NaSk2} as in our \textsection [\ref{SECRibbonAnisotropic}], proving both limit shape and Gaussian fluctuation results and giving several applications in \cite{HuangJack}.  Finally, Bufetov \cite{BufetovKMK} and Sodin \cite{SodinKMK} gave several probabilistic applications of Kerov's transition measures which we also analyzed in \textsection [\ref{SECRibbonAnisotropic}].

\subsection{Ribbon paths} \label{SUBSECRibbonPathsComments} \noindent Our Definition [\ref{DefinitionRibbonPaths}] of ribbon paths $\boldsymbol{\vec{\gamma}}$, together with our Theorem [\ref{Theorem3AOE}] relating Jack measures and ribbon paths, first appeared in the year 2015 in version 1 of the author's unpublished arXiv preprint \cite{Moll0}.  Our choice of the name ``ribbon paths'' and notation $W_{n,g,m}$ was meant to highlight similarities between our results and the theory of \textit{ribbon graphs} on non-oriented real surfaces.  Precisely, we noted in \cite{Moll0} that our Theorem [\ref{Theorem3AOE}] relating Jack measures $M(v^{\textnormal{out}}, v^{\textnormal{in}})$ to $W_{n,g,m}(\ell_1, \ldots, \ell_n |v^{\textnormal{out}}, v^{\textnormal{in}})$ in the $(v^{\textnormal{out}}, v^{\textnormal{in}})$-weighted enumeration Problem [\ref{ProblemRibbonPathEnumeration}] for connected ribbon {paths} on $n$ sites of lengths $\ell_1, \ldots, \ell_n$ with $n-1+g$ pairings and $m$ slides {has the same form as the refined topological expansion} for $\beta$-ensembles in a one-cut potential $V$ over ribbon {graphs} on $n$ vertices with $m$ twists of genus $g$ in Chekhov-Eynard \cite{Ch, ChEy1, ChEy2} and Borot-Guionnet \cite{BtGu1}.  That being said, in light of the breakthrough work of Chapuy-Do{\l}{\k{e}}ga \cite{ChapuyDolega2020}, we can now establish a definite connection between ribbon paths and ribbon graphs on non-oriented surfaces in \textsection [\ref{APPENDIXRibbonPaths}].  For further comments on the constructions in Chapuy-Do{\l}{\k{e}}ga \cite{ChapuyDolega2020} and their larger context, see \textsection [\ref{APPENDIXRibbonPaths}].
\appendix
\section{dePoissonization for Poissonized Jack-Plancherel measures} \label{APPENDIXdePoissonization}

\noindent In \textsection [\ref{SUBSECJackPlancherelLimitShapes}] and \textsection [\ref{SUBSECJackPlancherelGaussianFluctuations}], we illustrated our Theorems [\ref{Theorem1LLN}] and [\ref{Theorem2CLT}] in the case of the Poissonized Jack-Plancherel measures $M(v_{\textnormal{PL}}, v_{\textnormal{PL}})$.  We now show that a small modification of our moment method in \textsection [\ref{SECRibbonAnisotropic}] applies to the case of Jack-Plancherel measures $M_d^{\textnormal{PL}}$ on partitions of size $d$ \textcolor{black}{with law} \textcolor{black}{\textnormal{(\ref{JackPlancherelDoubleHookLaw})}}.  In doing so, we give new proofs of results in \cite{KeVe, LoSh, IvOl, DoFe, DoSni} as discussed in \textsection [\ref{SECComments}].

\begin{proposition} \label{dePoissonizationProposition} For $\hbar = \tfrac{1}{d}$, consider the $d \rightarrow \infty$ behavior of the random profile ${f}_{\lambda}( c | \ebar, \hbar)$ for $\lambda$ sampled from the Jack-Plancherel measures $M_d^{\textnormal{PL}}$ \textcolor{black}{with law} \textcolor{black}{\textnormal{(\ref{JackPlancherelDoubleHookLaw})}}.  Then in Regimes I, II, and III,
\begin{enumerate}
\item \textnormal{Theorem [\ref{Theorem1LLN}]} still holds with the same limit shapes $\mathbf{f}(c | v_{\textnormal{PL}}; \ebar)$.
\item \textnormal{Theorem [\ref{Theorem2CLT}]} still holds except the Gaussian process is different due to the conditioning $| \lambda|=d$ before $d \rightarrow \infty$.  In Regime II, the mean shift $\mathbf{X}(c | v_{\textnormal{PL}};0)$ is the same and the Gaussian process is $\mathbf{G}(c | v_{\textnormal{PL}}; 0)$ conditioned on the event $\int_{-\infty}^{+\infty} \mathbf{G}(c | v_{\textnormal{PL}}; 0) dc= 0$.
\end{enumerate}

\end{proposition} 

\noindent \textit{Proof of} Proposition [\ref{dePoissonizationProposition}]: Let $M_d(v^{\textnormal{out}}, v^{\textnormal{in}})$ be the Jack measure $M(v^{\textnormal{out}}, v^{\textnormal{in}})$ conditioned on the event $| \lambda|=d$.  \textcolor{black}{For Jack-Plancherel measures, $M_d^{\textnormal{PL}} = M_d(v_{\textnormal{PL}}, v_{\textnormal{PL}})$}.  Recall $\rho_{\mu}$ from \textsection [\ref{SUBSECJackPolynomials}].  Define \begin{equation} \label{TruncatedChocolate} \Upsilon^{(d)}_{v}(\cdot | \hbar) := \sum_{| \mu|=d} \frac{ \overline{V_{\mu}} \rho_{\mu}}{ || \rho_{\mu}||_{\hbar}^2 } \end{equation} so that the reproducing kernel in \textsection [\ref{SUBSUBSECStep1}] decomposes as $\Upsilon_v( \cdot | \hbar) = \sum_{d=0}^{\infty} \Upsilon_v^{(d)} ( \cdot |  \hbar)$.  Since the Nazarov-Sklyanin hierarchy $\widehat{T}^{\uparrow}_{\ell} (\ebar, \hbar)$ includes $\widehat{T}^{\uparrow}_2(\ebar, \hbar) = \sum_{k=1}^{\infty} \widehat{\rho}_{k} \widehat{\rho}_{-k}$ with eigenspaces $\mathcal{F}_d$ being the span of $\rho_{\mu}$ with $|\mu|=d$, the argument in \textsection [\ref{SUBSUBSECStep2}] may be repeated to prove that joint moments \begin{equation} \mathbb{E}_d[ T_{\ell_1}^{\uparrow} ( \ebar, \hbar) \cdots T_{\ell_n}^{\uparrow} (\ebar, \hbar) ] = \frac{ \langle \Upsilon_{v^{\textnormal{in}}}^{(d)} ( \cdot | \hbar),  \widehat{T}^{\uparrow}_{\ell_1}(\ebar, \hbar) \cdots \widehat{T}^{\uparrow}_{\ell_n} (\ebar, \hbar) \cdot  \Upsilon_{v^{\textnormal{out}}}^{(d)} (\cdot | \hbar) \rangle_{\hbar}}{\langle \Upsilon_{v^{\textnormal{in}}}^{(d)} ( \cdot | \hbar), \Upsilon_{v^{\textnormal{out}}}^{(d)} (\cdot | \hbar) \rangle_{\hbar}}  \end{equation} \noindent of $T_{\ell}^{\uparrow}( \ebar, \hbar)|_{\lambda}$ for $\lambda$ sampled from $M_d(v^{\textnormal{out}}, v^{\textnormal{in}})$ can be computed from (\ref{TruncatedChocolate}).  For $v=v_{\textnormal{PL}}$, (\ref{TruncatedChocolate}) is
\begin{equation} \Upsilon_{v_{\textnormal{PL}}}^{(d)} ( \cdot | \hbar) = \frac{ \overline{V_1}^d \rho_1^d }{ \hbar^d d!}. \end{equation} \noindent In this case, one has an analog of the eigenvalue relation $\widehat{\rho}_{-k} \Upsilon_v( \cdot | \hbar) = \overline{V_k} \Upsilon_v(\cdot | \hbar)$ from (\ref{HaveToLoveThisEigenvalueRelation}): 
\begin{equation} \label{CoolShiftsBro} \widehat{\rho}_{-k} \Upsilon_{v_{\textnormal{PL}}}^{(d)} = \overline{V_k} \Upsilon_{v_{\textnormal{PL}}}^{(d-1)} (\cdot | \hbar) \delta(k-1) \end{equation} 
 \noindent which shifts $d \mapsto d-1$.   Repeating the expansion of the path operators, Theorem [\ref{Theorem3AOE}] implies \begin{equation} \label{ExpansionMegaCool} \mathbb{E}_d [ T_{\ell_1}^{\uparrow} ( \ebar, \hbar) \cdots T_{\ell_n}^{\uparrow} (\ebar, \hbar) ] \sim \mathbb{E}[T_{\ell_1}^{\uparrow} ( \ebar, \hbar) \cdots T_{\ell_n}^{\uparrow} (\ebar, \hbar)] \cdot \mathcal{C}_d( \tfrac{1}{2} ||\vec{\ell}||_1; \hbar) \end{equation} \noindent to leading order as $\hbar \rightarrow 0$ with correction due to $\textcolor{black}{\eta=}\tfrac{1}{2} || \vec{\ell}||_1 = \tfrac{1}{2}(\ell_1+ \cdots + \ell_n)$ shifts in (\ref{CoolShiftsBro}) with \begin{equation} \label{ThisCThisCThisC} \mathcal{C}_d(\eta; \hbar) = \frac{ {\langle \Upsilon_{v^{\textnormal{out}}}^{(d - \eta)} ( \cdot | \hbar), \Upsilon_{v^{\textnormal{in}}}^{(d - \eta)} (\cdot | \hbar) \rangle_{\hbar}} }{ {\langle \Upsilon_{v^{\textnormal{out}}}^{(d)} ( \cdot | \hbar), \Upsilon_{v^{\textnormal{in}}}^{(d)} (\cdot | \hbar) \rangle_{\hbar}}} = \frac{ d(d-1) (d-2) \cdots (d- \eta +1) }{ \hbar^{- \eta} }. \end{equation} 
 
\noindent Specialize $\hbar = \frac{1}{d}$ so $f_{\lambda}(c | \ebar, \hbar)$ encloses area $1$ independent of $d$.  By (\ref{LeibnizLove}), pairings now cost a factor of $\tfrac{1}{d}$ and so to leading order as $d \rightarrow \infty$, only ribbon paths with no pairings contribute to (\ref{ExpansionMegaCool}).  Since $\hbar = \frac{1}{d}$, the correction in (\ref{ExpansionMegaCool}) is negligible since $\mathcal{C}_d (\eta; \tfrac{1}{d}) \rightarrow 1$ \textcolor{black}{in (\ref{ThisCThisCThisC})} as $d \rightarrow \infty$.  The asymptotic factorization (\ref{AsymptoticFactorizationLemmaFormula}) applied to (\ref{ExpansionMegaCool}) proves (1).  For (2), use (\ref{ExpansionMegaCool}) to write \begin{eqnarray} \ \ \ \ \ \ \ \  \textnormal{Cov}_d [ T_{\ell_1}^{\uparrow}(\ebar, \hbar), T_{\ell_2}^{\uparrow}(\ebar, \hbar) ]  &:=& \mathbb{E}_d [ T_{\ell_1}^{\uparrow}(\ebar, \hbar)  T_{\ell_2}^{\uparrow}(\ebar, \hbar) ]  - \mathbb{E}_d [ T_{\ell_1}^{\uparrow}(\ebar, \hbar) ] \mathbb{E}_d[ T_{\ell_2}^{\uparrow}(\ebar, \hbar) ] \\ \label{ThisLineIsItAnyway} & \sim &  + \mathbb{E} [ T_{\ell_1}^{\uparrow}(\ebar, \hbar)  T_{\ell_2}^{\uparrow}(\ebar, \hbar) ] \cdot \mathcal{C}_d (\ell_1 + \ell_2; \tfrac{1}{d}) \\ \label{SubtractThisAnyway} & \ &   - \mathbb{E} [ T_{\ell_1}^{\uparrow}(\ebar, \hbar) ] \mathbb{E}[ T_{\ell_2}^{\uparrow}(\ebar, \hbar) ] \cdot \mathcal{C}_d (\ell_1; \tfrac{1}{d}) \mathcal{C}_d ( \ell_2; \tfrac{1}{d})  \\ \label{VerifiedThisCheckTheStatistacs} & \sim & \textnormal{Cov}[ T_{\ell_1}^{\uparrow}(\ebar, \hbar), T_{\ell_2}^{\uparrow}(\ebar, \hbar)] - \frac{1}{d} \cdot  \ell_1 \ell_2 \cdot \mathbb{E}[ T_{\ell_1}^{\uparrow}(\ebar, \hbar) ] \mathbb{E}[ T_{\ell_2}^{\uparrow}(\ebar, \hbar) ].\end{eqnarray}

\noindent The terms (\ref{VerifiedThisCheckTheStatistacs}) come from those in (\ref{ThisLineIsItAnyway}) - (\ref{SubtractThisAnyway}) with (i) all $\mathcal{C}_d \sim 1$ or (ii) no pairings but $\mathcal{C}_d(\ell_1 + \ell_2; \tfrac{1}{d}) - \mathcal{C}_d(\ell_1; \tfrac{1}{d}) \mathcal{C}_d(\ell_2; \tfrac{1}{d}) \sim - \tfrac{1}{d} \cdot \ell_1, \ell_2$.  The same (\ref{ExpansionMegaCool}) implies that $n$th cumulants scaled by $d^{n/2}$ still vanish as $d \rightarrow \infty$ for $n \geq 3$ in any Regime I, II, III.  In Regime II when $\ebar \sim \hbar^{1/2} \sim d^{-1/2}$, since $\mathcal{C}_d$ produce only integer powers of $\tfrac{1}{d}$, again by (\ref{ExpansionMegaCool}), the mean shift $\mathbf{X}(c | v_{\textnormal{PL}}; 0)$ must be the same as (\ref{DFMeanShift}) and the covariance independent of $\alpha$.  Since $\mathbb{E}[ T_{\ell}^{\uparrow}(\ebar, \hbar)] \sim W_{1,0,0}(\ell | v_{\textnormal{PL}}, v_{\textnormal{PL}})$ are the Catalan numbers, the second term in (\ref{VerifiedThisCheckTheStatistacs}) cancels the contribution of the $k=1$ case of (\ref{PoissonizedIvanovOlshanski}) to the Poissonized covariance.  Conditioning that this variable is $0$ produces the same effect. $\square$

\section{Ribbon paths and ribbon graphs on non-oriented real surfaces} \label{APPENDIXRibbonPaths}

\noindent In \textsection [\ref{AppendixB1}] we express the $b$-deformed $2D$-Toda $\tau$-functions $\widehat{F}^G$ introduced by Chapuy-Do{\l}{\k{e}}ga \cite{ChapuyDolega2020} with arbitrary weights $G$ in terms of Jack measures.  In \textsection [\ref{AppendixB2}] we introduce a class of $b$-dependent weight functions $G_n^{\uparrow}$ and show that $\widehat{F}^{G_n^{\uparrow}}$ is determined by the $n$th joint moments of the random variables $T_{\ell}^{\uparrow} (\ebar, \hbar)|_{\lambda}$ which we analyzed in \textsection [\ref{SECRibbonAnisotropic}].  In \textsection [\ref{AppendixB3}], we combine our Theorem [\ref{Theorem3AOE}] with Theorem 6.2 in \cite{ChapuyDolega2020} to relate our ribbon paths to the ribbon graphs on non-oriented surfaces in \cite{ChapuyDolega2020}.  

\subsection{Jack measures and $b$-deformations of 2D-Toda $\tau$-functions} \label{AppendixB1}

\noindent In this section, we show that $b$-deformations of the $2D$-Toda $\tau$-functions introduced by Chapuy-Do{\l}{\k{e}}ga \cite{ChapuyDolega2020} may be expressed in terms of Jack measures $M(v^{\textnormal{out}}, v^{\textnormal{in}})$.  This connection is immediate from the definition above of Jack measures and is only a matter of changing variables and conventions appropriately.  Recall that the scaled content of a box $\square \in \lambda$ in a Young diagram in the $i$th row and $j$th column is \begin{equation} \label{NekOkContent} c(\square | \ebar, \hbar) = \textcolor{black}{\varepsilon_2 (i-1) + \varepsilon_1 ( j-1)} \end{equation} \noindent where $\varepsilon_2<0< \varepsilon_1$ come from $\ebar, \hbar$ by (\ref{OmegaVariables}), (\ref{Ebar}), (\ref{Hbar}) according to the conventions in \cite{NekOk}.  In the anisotropic profile $f_{\lambda}(c | \ebar, \hbar)$, the scaled content (\ref{NekOkContent}) is the $c$-coordinate of the lowest corner of $\square$.  
\begin{proposition} \label{ConvertFGToJackMeasureProposition} Let $\ebar \in \mathbb{R}$ and $\hbar>0$.  Let $v^{\textnormal{out}} = \{V_k^{\textnormal{out}}\}_{k=1}^{\infty} ,v^{\textnormal{in}} = \{V_k^{\textnormal{in}}\}_{k=1}^{\infty}$ be two arbitrary specializations.  For any indeterminate $t$ and any $G: \mathbb{R} \rightarrow \mathbb{C}$, the $b$-deformed $2D$-Toda $\tau$-function $\widehat{F}^G$ in \textnormal{\textsection 6.2} of \textnormal{Chapuy-Do{\l}{\k{e}}ga \cite{ChapuyDolega2020}} with higher times $\mathbf{p}_k, \mathbf{q}_k$ may be expressed as \begin{equation} \label{ConvertFGToJackMeasure}\widehat{F}^{G} = \textnormal{exp} \Bigg ( \frac{1}{\hbar} \sum_{k=1}^{\infty}  \frac{ \overline{V_k^{\textnormal{out}}} V_k^{\textnormal{in}} }{k} \Bigg ) \cdot \mathbb{E} \Bigg [ t^{|\lambda|} \prod_{\square \in \lambda} G\big ( c(\square | \ebar, \hbar) \big ) \Bigg ] \end{equation}
\noindent the product of Stanley's Cauchy kernel \textnormal{(\ref{CauchyIdentity})} with an expectation $\mathbb{E}$ against the Jack measure $M(v^{\textnormal{out}}, v^{\textnormal{in}})$ after the change of variables \begin{eqnarray} \label{WebbieFlowOut} \mathbf{p}_k &=& \overline{V_k^{\textnormal{out}}} / (- \varepsilon_2) \\ \label{WebbieFlowIn} \mathbf{q}_k &=& V_k^{\textnormal{in}} / (- \varepsilon_2)  \\ \label{WhereThebTho} b &=& \  \ebar \ / (- \varepsilon_2). \end{eqnarray} 

\end{proposition}

\begin{itemize}

\item \textit{Proof of} Proposition [\ref{ConvertFGToJackMeasureProposition}]: By (\ref{Convention2}), the $\alpha$-content $c_{\alpha} (\square) = -\textcolor{black}{(i-1)} + \alpha \textcolor{black}{(j-1)}$ in \cite{Mac, Stanley} is related to (\ref{NekOkContent}) by $c(\square | \ebar, \hbar) = (- \varepsilon_2) c_{\alpha}(\square)$.  As written, $G$ in (\ref{ConvertFGToJackMeasure}) is a function of $c(\square | \ebar, \hbar)$, not $c_{\alpha}(\square)$.  On the other hand, formula (63) in \cite{ChapuyDolega2020} features a weight $G$ evaluated on $c_{\alpha}(\square)$.  To resolve this discrepancy, by Remark 14 of \cite{ChapuyDolega2020}, the authors discuss a distinguished extra variable that may be inserted into their weight function $G$ to clarify the expansion of $\log \widehat{F}^G$ as a topological expansion.  We choose this variable to be the $(- \varepsilon_2)$ introduced by Nekrasov-Okounkov \cite{NekOk}.  After making use of their Remark 14, their weight function $G$ in their formula (63) is now a function of (\ref{NekOkContent}).  The desired claim now follows from formula (63) in \cite{ChapuyDolega2020} using Definition [\ref{DefinitionJackMeasures}], $b= \alpha-1$, (\ref{Convention2}), and (\ref{Convention1}). $\square$

\end{itemize}

\noindent At $\ebar=0$ so that $\alpha=1$ and $b=0$, (\ref{ConvertFGToJackMeasure}) recovers the relationship between the $2D$-Toda $\tau$-functions and the Schur measures of Okounkov \cite{Ok1}.  For recent analysis of (\ref{ConvertFGToJackMeasure}) for $\ebar = 0$ and arbitrary $G$, see Guay-Paquet-Harnad \cite{GuayPaquetHarnad} and Alexandrov-Chapuy-Eynard-Harnad \cite{AlexandrovChapuyEynardHarnad2018} and references therein.  In this case $\ebar =0$, (\ref{ConvertFGToJackMeasure}) is an expression for these $\tau$-functions in \textit{content-product form}, a name which refers directly to contents (\ref{NekOkContent}) of $\square \in \lambda$ which appear in the product over boxes in (\ref{ConvertFGToJackMeasure}).

\subsection{Transition measure joint moments from $b$-dependent weight $G_n^{\uparrow}$} \label{AppendixB2}

\noindent In sections \textsection 1 - \textsection 5 of \cite{ChapuyDolega2020}, Chapuy-Do{\l}{\k{e}}ga first analyze their $b$-deformed $\tau$-functions (\ref{ConvertFGToJackMeasure}) for a large class of weights $G_{\tilde{k}} (c ) = \prod_{a=1}^{\tilde{k}}  (  1+ \tilde{u}_a   c  )$ defined by an arbitrary number $\widetilde{k}$ of indeterminates $\tilde{u}_{1}, \ldots, \tilde{u}_{\tilde{k}}$.  In \textsection 6 of \cite{ChapuyDolega2020}, they use their results for $\widehat{F}^{G_{\tilde{k}}}$ over all $\widetilde{k}$ to derive results for arbitrary weights $G(c)$ that are not necessarily polynomials in $c$.  We now introduce a different class of weights $G_n^{\uparrow}$ for which the $b$-deformed $\tau$-function (\ref{ConvertFGToJackMeasure}) recovers the joint moments of $T_{\ell}^{\uparrow}(\ebar, \hbar)|_{\lambda}$ we analyzed above in \textsection [\ref{SECRibbonAnisotropic}].    These weights $G_n^{\uparrow}$ arise in defining bulk $\mathcal{Y}$ observables in Nekrasov \cite{NekYI} and hence $qq$-characters and discrete loop equations in \cite{NekYI} independently studied by Borodin-Gorin-Guionnet \cite{BoGoGu}. \begin{definition} Fix $\ebar \in \mathbb{R}$, $\hbar>0$, and $\varepsilon_2 < 0 < \varepsilon_1$ determined from $\ebar, \hbar$ by \textnormal{(\ref{OmegaVariables}), (\ref{Ebar}), (\ref{Hbar})}. For any $n=1,2,3,\ldots$ and $u_1, \ldots, u_n \in \mathbb{C} \setminus \mathbb{R}$, define $G_n^{\uparrow}: \mathbb{R} \rightarrow \mathbb{C}$ by \begin{equation} \label{GUpWeights} G_n^{\uparrow}(c) = \prod_{a=1}^n \frac{ (u_a - c) ( u_a - c - \varepsilon_2 -  \varepsilon_1)}{ (u_a - c - \varepsilon_2) (u_a - c - \varepsilon_1)}. \end{equation} \noindent  We may write $G_n^{\uparrow}( c ) = G_n^{\uparrow}( c | u_1, \ldots, u_n ; \ebar, \hbar)$ to emphasize the variables defining $G_n^{\uparrow}$.
\end{definition} \noindent Unlike the weights $G_{\tilde{k}}$ from \cite{ChapuyDolega2020}, the weights (\ref{GUpWeights}) depend on both $\ebar$ and $\hbar$, hence on $b$ by (\ref{WhereThebTho}).  To relate (\ref{GUpWeights}) to our results above, for any $u \in \mathbb{C} \setminus \mathbb{R}$, consider the generating function \begin{equation} \label{TUpObservableAppendix} T^{\uparrow}(u | \ebar, \hbar)|_{\lambda} = \sum_{\ell=0}^{\infty} u^{- \ell-1} T_{\ell}^{\uparrow}(\ebar, \hbar) |_{\lambda} \end{equation} \noindent for the $T_{\ell}^{\uparrow} (\ebar, \hbar)|_{\lambda}$ defined in (\ref{TransitionStatistics}).  In Kerov \cite{Ke1, Ke4}, (\ref{TUpObservableAppendix}) is the Stieltjes transform of the transition measure associated to the anisotropic partition profile $f_{\lambda}(c | \ebar, \hbar)$.  

\begin{proposition} \label{ConvertFGUPToJackMeasureProposition} Fix $\ebar \in \mathbb{R}$, $\hbar>0$, and $u_1, \ldots, u_n \in \mathbb{C} \setminus \mathbb{R}$.  For the $b$-dependent weights $G_n^{\uparrow}$ in \textnormal{(\ref{GUpWeights})}, the $b$-deformed $2D$-Toda $\tau$-function $\widehat{F}^{G_n^{\uparrow}}$ in \textnormal{(\ref{ConvertFGToJackMeasure})} is determined by the $n$th joint moments of the random variables $T^{\uparrow}_{\ell}(\ebar, \hbar)|_{\lambda}$ analyzed above in \textnormal{Theorem [\ref{Theorem3AOE}]} except for $\lambda$ sampled from the $t$-rescaled Jack measure $M_t(v^{\textnormal{out}}, v^{\textnormal{in}})$ with specializations $\{t^{k/2} V_k^{\textnormal{out}}\}_{k=1}^{\infty}$, $\{t^{k/2} V_k^{\textnormal{in}}\}_{k=1}^{\infty}$ by
\begin{equation} \label{ConvertFGUPToJackMeasure} \widehat{F}^{G_n^{\uparrow}} = \textnormal{exp} \Bigg ( \frac{1}{\hbar} \sum_{k=1}^{\infty}t^k  \frac{\overline{V_k^{\textnormal{out}}} V_k^{\textnormal{in}}}{k} \Bigg ) \sum_{\ell_1, \ldots, \ell_n=0}^{\infty} u_1^{-\ell_1} \cdots u_n^{-\ell_n} \mathbb{E}_t \Big [ T^{\uparrow}_{\ell_1} (\ebar, \hbar)|_{\lambda} \cdots T_{\ell_n}^{\uparrow}(\ebar, \hbar)|_{\lambda} \Big ]  . \end{equation}
\end{proposition} 
\begin{itemize}\item  \textit{Proof of} Proposition [\ref{ConvertFGUPToJackMeasureProposition}]: The definition of (\ref{TUpObservableAppendix}) in terms of $f_{\lambda}''(c | \ebar, \hbar)$ in (\ref{TransitionStatistics}) implies \begin{equation} \label{TakingProductOfThisnTimes} u T^{\uparrow}(u  | \ebar, \hbar)|_{\lambda} = \prod_{\square \in \lambda} \frac{ \big (u - c(\square | \ebar, \hbar) \big ) \big ( u - c(\square | \ebar, \hbar) - \varepsilon_2 - \varepsilon_1\big )}{ \big (u - c(\square | \ebar, \hbar) - \varepsilon_2 \big ) \big ( u - c(\square | \ebar, \hbar) - \varepsilon_1 \big )} .\end{equation} \noindent Taking a product of (\ref{TakingProductOfThisnTimes}) $n$ times, once for each $u_a$ for $a=1,2,\ldots, n$, gives \begin{equation} \label{TUpToGUp} \prod_{a=1}^n u_a T^{\uparrow}(u_a | \ebar, \hbar) |_{\lambda} = \prod_{\square \in \lambda} G_n^{\uparrow} \big ( c( \square | \ebar, \hbar) \big ) \end{equation} \noindent for $G_n^{\uparrow}$ in (\ref{GUpWeights}).  For any function $O|_{\lambda}$ of $\lambda$, if $\mathbb{E}$ and $\mathbb{E}_t$ are expectations against $M(v^{\textnormal{out}}, v^{\textnormal{in}})$ and its $t$-rescaling $M_t(v^{\textnormal{out}}, v^{\textnormal{in}})$, respectively, we always have
\begin{equation} \label{BalancedByTheT} \textnormal{exp} \Bigg ( \frac{1}{\hbar} \sum_{k=1}^{\infty} \frac{\overline{V_k^{\textnormal{out}}} V_k^{\textnormal{in}}}{k} \Bigg )\mathbb{E} [ t^{|\lambda|} O|_{\lambda}] = \textnormal{exp} \Bigg ( \frac{1}{\hbar} \sum_{k=1}^{\infty} t^k \frac{\overline{V_k^{\textnormal{out}}} V_k^{\textnormal{in}}}{k} \Bigg )  \mathbb{E}_t[ O |_{\lambda}]. \end{equation}
\noindent The case $G = G_n^{\uparrow}$ of Proposition [\ref{ConvertFGToJackMeasureProposition}], (\ref{TakingProductOfThisnTimes}), (\ref{TUpToGUp}), and (\ref{BalancedByTheT}) imply (\ref{ConvertFGUPToJackMeasure}). $\square$ 
\end{itemize}

\subsection{Ribbon paths and ribbon graphs with $b$-dependent weight $G^{\uparrow}$}\label{AppendixB3}
\noindent In Proposition [\ref{RibbonPathsInChapuyDolega}] below, we relate our weighted enumeration of ribbon paths above in our Theorem [\ref{Theorem3AOE}] to the weighted enumeration of ribbon graphs on non-oriented surfaces in Chapuy-Do{\l}{\k{e}}ga \cite{ChapuyDolega2020}.  Let $\Delta_n$ be the set of ribbon paths $\boldsymbol{\vec{\gamma}}$ on $n$ sites \textcolor{black}{as in Definition [\ref{DefinitionRibbonPaths}]}.  For any $\boldsymbol{\vec{\gamma}}$, let
\begin{itemize}
\item $\mathbf{d} [\boldsymbol{\vec{\gamma}}]$ be the size of $\boldsymbol{\vec{\gamma}}$ \textcolor{black}{as in Definition [\ref{DefinitionSizeRibbonPath}]}
\item $\mathbf{q} [\boldsymbol{\vec{\gamma}}]$ be the number of pairings in $\boldsymbol{\vec{\gamma}}$
\item $\mathbf{m} [\boldsymbol{\vec{\gamma}}]$ be the number of slides in $\boldsymbol{\vec{\gamma}}$
\item $\boldsymbol{\ell}_1[\boldsymbol{\vec{\gamma}}], \ldots, \boldsymbol{\ell}_n[\boldsymbol{\vec{\gamma}}]$ be the lengths of the $n$ sites in $\boldsymbol{\vec{\gamma}}$.
\end{itemize}
\noindent Recall the length of a partition $\mu$ is ${\ell}(\mu) = \# \{i : \mu_i  > 0\}$.  By Definition [\ref{DefinitionUnpairedJumpProfiles}], for any $\boldsymbol{\vec{\gamma}} \in \Delta_n$, \begin{equation} \boldsymbol{\ell}_1[\boldsymbol{\vec{\gamma}}] +  \cdots + \boldsymbol{\ell}_n[\boldsymbol{\vec{\gamma}}]  = {\ell}(\boldsymbol{\mu}^{{+}}[ \boldsymbol{\vec{\gamma}}]) + {\ell} (\boldsymbol{\mu}^{{-}}[ \boldsymbol{\vec{\gamma}}]) + 2\mathbf{q} [\boldsymbol{\vec{\gamma}}] + \mathbf{m} [\boldsymbol{\vec{\gamma}}]. \end{equation}  For any $t, \hbar, \ebar, u_1, \ldots, u_n$ and any specializations $v^{\textnormal{out}} = \{V_k^{\textnormal{out}} \}_{k=1}^{\infty}$, $v^{\textnormal{in}} = \{V_k^{\textnormal{in}}\}_{k=1}^{\infty}$, define a sum \begin{equation} \label{MySum} I^{\uparrow}_n := \sum_{\boldsymbol{\vec{\gamma}} \in \Delta_n} t^{\mathbf{d} [ \boldsymbol{\vec{\gamma}}]} \ \hbar^{\mathbf{q} [ \boldsymbol{\vec{\gamma}}]}
\  \textcolor{black}{\ebar}^{\mathbf{m} [ \boldsymbol{\vec{\gamma}}]} \  u_1^{-\boldsymbol{\ell}_1 [ \boldsymbol{\vec{\gamma}}]} \cdots  u_n^{-\boldsymbol{\ell}_n [ \boldsymbol{\vec{\gamma}}]}  \cdot \mathfrak{W}( \boldsymbol{\vec{\gamma}} | v^{\textnormal{out}}, v^{\textnormal{in}}) \end{equation} where $\mathfrak{W}$ is the weight in Definition [\ref{DefinitionRibbonPathsWeight}] which we recall is \begin{equation} \label{SeriouslyLookHere} \mathfrak{W}( \boldsymbol{\vec{\gamma}} | v^{\textnormal{out}}, v^{\textnormal{in}})  =  \prod_{ \textbf{p} \in \boldsymbol{\vec{\gamma}}}  \textnormal{size}(\textbf{p}) \prod_{\textit{\textbf{e}} \in \textnormal{\textbf{E}}_S(\boldsymbol{\vec{\gamma}})} \textnormal{height}(\textit{\textbf{e}}) \prod_{\textbf{\textit{e}} \in \mathbf{E}_{UJ}^-(\boldsymbol{\vec{\gamma}})} \overline{V_{-\deg ( \textbf{\textit{e}})}^{\textnormal{out}}} \prod_{\textbf{\textit{e}} \in \mathbf{E}_{UJ}^+(\boldsymbol{\vec{\gamma}})} {V_{\deg ( \textbf{\textit{e}})}^{\textnormal{in}}}. \end{equation} 
\noindent Next, let $\Gamma_{\textnormal{CD}}$ be the set of ribbon graphs of maps on non-oriented surfaces which are connected rooted normal infinite constellations $(\mathbf{M}, \mathbf{c})$ of any size from \textnormal{\textsection 6.2} of \textnormal{Chapuy-Do{\l}{\k{e}}ga \cite{ChapuyDolega2020}}.  For any $(\mathbf{M},\mathbf{c}) \in \Gamma_{CD}$, any $G: \mathbb{R} \rightarrow \C$, and any measure of non-orientability $\rho$ in \cite{ChapuyDolega2020}, let \begin{itemize}
\item $|\mathbf{M}|$ be the size of the constellation $\mathbf{M}$
\item $\mathbf{m}[ \mathbf{M},\mathbf{c}]$ be the exponent of the $b$-weight of $(\mathbf{M}, \mathbf{c})$, denoted $\nu_{\rho}(\mathbf{M}, \mathbf{c})$ in \textnormal{\cite{ChapuyDolega2020}}
\item $\boldsymbol{\chi}[\mathbf{M}]$ be the Euler characteristic of $\mathbf{M}$
\item $\mathbf{F}[\mathbf{M}]$ be the set of faces $\mathbf{f}$ in $\mathbf{M}$ of any degree $\deg (\mathbf{f}) \in \{1,2,\ldots\}$
\item $\mathbf{V}_0 [ \mathbf{M}]$ be the set of vertices $\mathbf{v}$ in $\mathbf{M}$ of colour $0$ of any degree $\deg ( \mathbf{v}) \in \{1,2,3,\ldots\}$
\item $\mathbf{U}^G [ \mathbf{M}]$ be the $G$-weight of $\mathbf{M}$, denoted $f_{\mathfrak{v}[\mathbf{M}]}(g_1, \ldots)$ in \textnormal{\cite{ChapuyDolega2020}}
\end{itemize}

\noindent In addition, let $\boldsymbol{\chi}_0 [ \mathbf{M}, \mathbf{c}]$ denote the modified Euler characteristic \begin{equation} \label{ModifiedEulerCharacteristic} \boldsymbol{\chi}_0 [ \mathbf{M}, \mathbf{c}] = \boldsymbol{\chi}[\mathbf{M}] + \mathbf{m}[ \mathbf{M},\mathbf{c}] \end{equation}

\noindent which agrees with the usual Euler characteristic if $\mathbf{M}$ is orientable.\\
\\
\noindent For any $t, \ebar, \varepsilon_2$ satisfying (\ref{OmegaVariables}), (\ref{Ebar}), (\ref{Hbar}) and any specializations $v^{\textnormal{out}} = \{V_k^{\textnormal{out}} \}_{k=1}^{\infty}$, $v^{\textnormal{in}} = \{V_k^{\textnormal{in}}\}_{k=1}^{\infty}$, and any weight $G : \mathbb{R} \rightarrow \mathbb{C}$, consider the weighted sum over ribbon graphs  \begin{equation} \label{TheirSum} S_{CD}^{G} := \sum_{(\mathbf{M}, \mathbf{c}) \in \Delta_{CD}} t^{| \mathbf{M}|} \ \ebar^{\mathbf{m} [ \mathbf{M}, \mathbf{c}]} (-\varepsilon_2)^{-\boldsymbol{\chi}_0[\mathbf{M}, \mathbf{c}]}   \prod_{\mathbf{f} \in \mathbf{F}[ \mathbf{M}]} \overline{V_{\deg ( \mathbf{f})}^{\textnormal{out}}} \prod_{\mathbf{v} \in \mathbf{V}_0 [\mathbf{M}]} V_{\deg (\mathbf{v})}^{\textnormal{in}} \cdot \mathbf{U}^{G}[\mathbf{M}].
\end{equation}

\begin{proposition} \label{RibbonPathsInChapuyDolega} 
\noindent For the weight $G =G_n^{\uparrow}$ in \textnormal{(\ref{GUpWeights})}, the weighted sum of ribbon graphs $S_{CD}^G$ in \textnormal{(\ref{TheirSum})} from \textnormal{Chapuy-Do{\l}{\k{e}}ga \cite{ChapuyDolega2020}} and our weighted sum $I_n^{\uparrow}$ of ribbon paths in \textnormal{(\ref{MySum})} are related by
\begin{equation} \label{TheirSumToMySum} S_{CD}^{G^{\uparrow}_n} = \frac{1+b}{\hbar} \sum_{k=1}^{\infty} t^k \overline{V_k^{\textnormal{out}}} V_k^{\textnormal{in}} +  (1 + b) t \frac{\partial}{\partial t} \log I_n^{\uparrow} .\end{equation}
 \end{proposition}

\begin{itemize}
\item \textit{Proof of \textnormal{Proposition [\ref{RibbonPathsInChapuyDolega}]}}: For any $G$, $S_{CD}^G$ in (\ref{TheirSum}), and $\widehat{F}^G$ in (\ref{ConvertFGToJackMeasure}), \textcolor{black}{we claim} \begin{equation} \label{Their62Assertion} S_{CD}^G = ( 1 + b ) t \frac{\partial }{\partial t} \log \widehat{F}^{G}. \end{equation} \noindent This claim (\ref{Their62Assertion}) is precisely the statement of Theorem 6.2 in \cite{ChapuyDolega2020} after applying Remark 14 from \cite{ChapuyDolega2020} with parameter $(-\varepsilon_2)$ from \cite{NekOk} as in our proof of Proposition [\ref{ConvertFGToJackMeasureProposition}].  On the one hand, as pointed out in Remark 14 of \cite{ChapuyDolega2020}, this scaling introduces a factor of $(-\varepsilon_2)^{\ell(\mathbf{\mu}^{-}[\mathbf{M}]) + \ell(\mathbf{\mu}^{+} [ \mathbf{M}]) - \boldsymbol{\chi}[\mathbf{M}]}$ where $\mathbf{\mu}^{\pm}[ \mathbf{M}]$ are the ramification profiles of the generalized branched cover $\mathcal{S} \rightarrow \mathbb{S}^2$ associated to the constellation $\mathbf{M}$ over two distinguished points $z^{\pm} \in \mathbb{S}^2$.  In the constellation $M$, $\mathbf{\mu}^{\pm}[ \mathbf{M}]$ keep track of the number of faces and vertices of colour $0$ of each degree, respectively.  On the other hand, the change of variables (\ref{WebbieFlowOut}), (\ref{WebbieFlowIn}) introduces a factor of $(-\varepsilon_2)^{-\ell(\mathbf{\mu}^{-}[\mathbf{M}]) - \ell(\mathbf{\mu}^{+} [ \mathbf{M}])}$.  The lengths of the ramification profiles in the exponents cancel, leaving $(- \varepsilon_2)^{- \boldsymbol{\chi}[\mathbf{M}]}$.  Finally, changing variables from $b$ to $\ebar$ by (\ref{WhereThebTho}) introduces a factor $(- \varepsilon_2)^{- \mathbf{m}[ \mathbf{M}, \mathbf{c}]}$ which we can absorb into the Euler characteristic by means of the modified Euler characteristic (\ref{ModifiedEulerCharacteristic}).  This verifies the agreement between (\ref{Their62Assertion}) and Theorem 6.2 in \cite{ChapuyDolega2020}.  Next, apply Part (1) of Theorem [\ref{Theorem3AOE}] in the case of the $t$-rescaled Jack measure $M(v^{\textnormal{out}}, v^{\textnormal{in}})$ with specializations $\{t^{k/2} V_k^{\textnormal{out}}\}_{k=1}^{\infty} , \{ t^{k/2} V_k^{\textnormal{in}}\}_{k=1}^{\infty}$.  The resulting $\ebar$ and $\hbar$ expansion of the $n$th joint moments of the random variables $T_{\ell}^{\uparrow}(\ebar, \hbar)|_{\lambda}$ is precisely the weighted sum $I_n^{\uparrow}$ in (\ref{MySum}): \begin{equation} \label{ScaleThatTHollaAtchaWodie} \mathbb{E}_t[T_{\ell_1}^{\uparrow}(\ebar, \hbar)|_{\lambda} \cdots T_{\ell_n}^{\uparrow}(\ebar, \hbar) |_{\lambda}] = I_n^{\uparrow} \end{equation}
\noindent Indeed, we know the result at $t=1$ and \textcolor{black}{how} specializations $\overline{V^{\textnormal{out}}_{- \textnormal{deg}(\mathbf{e})}}$, $V^{\textnormal{in}}_{\deg(\mathbf{e})}$ \textcolor{black}{weight} unpaired jumps $\mathbf{e} \in \mathbf{E}_{UJ}^{\mp} [ \boldsymbol{\vec{\gamma}}]$ of negative and positive degree, respectfully.  The scaling $V_k \mapsto t^{k/2} V_k$ in each specialization contributes an overall factor of $t^{ \mathbf{d}[ \boldsymbol{\vec{\gamma}}]}$ where $\mathbf{d}[ \boldsymbol{\vec{\gamma}}]$ is the size of the ribbon path in Definition [\ref{DefinitionSizeRibbonPath}].  Finally, substituting our $b$-dependent weight $G = G_n^{\uparrow}$ in (\ref{Their62Assertion}) and applying Proposition [\ref{ConvertFGUPToJackMeasureProposition}] and (\ref{ScaleThatTHollaAtchaWodie}) yields (\ref{TheirSumToMySum}). $\square$
\end{itemize}
\noindent The relation (\ref{Their62Assertion}) in Theorem 6.2 of \cite{ChapuyDolega2020} is a highly non-trivial connection between ribbon graphs on surfaces and Jack measures justified by several novel techniques and insights of Chapuy-Do{\l}{\k{e}}ga.  The transparent parallels above between ribbon graphs and ribbon paths in (\ref{SeriouslyLookHere}) and (\ref{TheirSum}), as well as the geometric implications of our main Theorems [\ref{Theorem1LLN}] and [\ref{Theorem2CLT}], remain to be explored.\\
\\
\noindent \textcolor{black}{\textbf{Acknowledgements.} The author would like to thank Alexei Borodin, Ga\"{e}tan Borot, Percy Deift, Vadim Gorin, Ryan Mickler, and Jonathan Novak for helpful discussions and the referees for their careful reading and thoughtful comments on the manuscript.}
{\footnotesize

\bibliographystyle{plain}
\bibliographystyle{amsalpha}

\bibliography{Bbib2021june19final}}

\end{document}